\documentclass[review,12pt]{elsarticle}

\usepackage{setspace}
\usepackage{subcaption}
\usepackage{multirow}
\usepackage{multicol}
\usepackage{graphicx}
\RequirePackage{amsmath}
\usepackage{bm}
\usepackage{xfrac}
\usepackage{float}
\usepackage{arydshln} 

\usepackage{amsmath,esint,amsfonts}

\usepackage[most]{tcolorbox}
\usepackage[colorinlistoftodos]{todonotes}

\usepackage{lineno,hyperref}
\usepackage{cleveref} 

\journal{arXiv}

\biboptions{sort&compress} 

\begin{document}
\begin{frontmatter}

	\title{A Semi-Implicit Meshless Method for Incompressible Flows in Complex Geometries}

	\author{Shantanu Shahane\fnref{Corresponding Author}}
	\author{Surya Pratap Vanka}
	\address{Department of Mechanical Science and Engineering,\\
		University of Illinois at Urbana-Champaign, Urbana, Illinois 61801}
	\fntext[Corresponding Author]{\vspace{0.3cm}Corresponding Author Email Address: \url{shahaneshantanu@gmail.com}}


	%

	\begin{abstract}
		We present an exponentially convergent semi-implicit meshless algorithm for the solution of Navier-Stokes equations in complex domains. The algorithm discretizes partial derivatives at scattered points using radial basis functions as interpolants. Higher-order polynomials are appended to polyharmonic splines (PHS-RBF) and a collocation method is used to derive the interpolation coefficients. The interpolating kernels are then differentiated and the partial-differential equations are satisfied by collocation at the scattered points. The PHS-RBF interpolation is shown to be exponentially convergent with discretization errors decreasing as a high power of a representative distance between points. We present here a semi-implicit algorithm for time-dependent and steady state fluid flows in complex domains. At each time step, several iterations are performed to converge the momentum and continuity equations. A Poisson equation for pressure corrections is formulated by imposing divergence free condition on the iterated velocity field. At each time step, the momentum and pressure correction equations are repeatedly solved until the velocities and pressure converge to a pre-specified tolerance. We have demonstrated the convergence and discretization accuracy of the algorithm for two model problems and simulated three other complex problems. In all cases, the algorithm is stable for Courant numbers in excess of ten. The algorithm has the potential to accurately and efficiently solve many fluid flow and heat transfer problems in complex domains. An open source code Meshless Multi-Physics Software (MeMPhyS) \cite{shahanememphys} is available for interested users of the algorithm.
	\end{abstract}


\end{frontmatter}


\section{Introduction}
Numerical algorithms to solve the coupled Navier-Stokes equations for incompressible flows can be categorized as: i) time marching pressure-projection methods \cite{harlow1965numerical,chorin1968numerical,temam1969approximation,kim1985application,rosenfeld1991fractional,shahane2019finite}, ii) semi-implicit equation-by-equation decoupled iterative schemes \cite{issa1986solution,van1984enhancements,patankar1983calculation,patankar1983prediction,patankar2018numerical}, or fully-coupled block-implicit methods using Newton methods for the coupled nonlinear discrete equations \cite{vanka1985block,braaten1990fully,darwish2009coupled,vanka1983fully,chen2010coupled}. The appropriate technique for a given flow depends on the flow characteristics and the desired temporal information from the numerical solution. The explicit fractional step (pressure-projection) methods are ideal for problems with needs to accurately calculate the flow instabilities and turbulent fluctuations, such as in Direct and Large-eddy Simulations (DNS and LES) of turbulence \cite{sohankar2004and,manhart2004zonal,ha2018gpu,kalitzin2003dns,you2000modified}. In such explicit methods for incompressible Navier-Stokes equations, the time step is based on a combination of velocity and diffusion Courant numbers or only on the velocity Courant number, if the diffusion terms are treated implicitly. A number of variants of pressure-projection methods with second-order temporal accuracy have been developed in combination with finite-difference, finite-volume, spectral and spectral-element methods \cite{robichaux1999three,nek5000,spalart1991spectral,tafti1991numerical,tafti1991three}.
\par If the temporal variations are not fast (slow transients), or only the steady state flow fields are desired, one can use a semi-implicit method in which the time-dependent equations are solved with larger time steps (several times the velocity Courant number) with iterations over the coupled equations at each time step. For steady state problems, the temporally evolving solution is of interest only in passing, and therefore, moderate temporal errors can be accepted. An alternative to time marching is a direct relaxation method of the coupled equations using an iterative correction strategy. Such a strategy is followed in the well-known SIMPLE series \cite{patankar1983calculation, patankar1983prediction} of algorithms where, a pressure-correction equation is derived from truncated discrete momentum equations to satisfy the divergence-free condition on the velocities. A number of different pressure-correction strategies \cite{patankar2018numerical, van1984enhancements, issa1986solution} have been proposed as an improvement to the original SIMPLE algorithm. Since the momentum equations are nonlinear, an iterative method is always needed. Several works \cite{vanka1983fully, darwish2009coupled, chen2010coupled, braaten1990fully, vanka1985block} have combined Newton type methods with coupled direct or iterative solvers to efficiently solve the discrete equations. If the initial guess is reasonably good, the coupled solvers converge rapidly in a small number of iterations. However, coupled direct solvers are not attractive for large numbers of nonlinear equations usually required for industrial flow simulations. Iterative methods for such coupled equations combining Krylov subspace methods \cite{pawlowski2006globalization, knoll2004jacobian} and multigrid techniques \cite{vanka1986blockcmame, vanka1986blockjcp, john2000numerical} have also been pursued.
\par Recently, there has been growing interest in solving the Navier-Stokes equations in complex domains using meshless methods that discretize the governing partial differential equations using only scattered points. Among the several types of meshless methods \cite{monaghan2012smoothed, perrone1975general, liu1995reproducing, belytschko1994element, liszka1996hp, babuvska1997partition, onate1996finite}, interpolation with radial basis functions provides a novel way to construct discrete forms of the Navier Stokes equations \cite{shu2003local, ding2006numerical, sanyasiraju2008local, vidal2016direct, zamolo2019solution, shahane2020high, bayona2017role_II, kosec2008solution, kosec2020radial}. Several RBFs have been previously used for such interpolation of values at scattered locations \cite{flyer2016onrole_I}. These include:
\begin{equation}
\begin{aligned}
\text{Multiquadrics (MQ): } \phi(r)=&(r^2 + \epsilon ^2)^{1/2}\\
\text{Inverse Multiquadrics (IMQ): } \phi(r)=&(r^2 + \epsilon ^2)^{-1/2}\\
\text{Gaussian: } \phi(r)=&\exp\left(\frac{-r^2}{\epsilon ^2}\right)\\
\text{Polyharmonic Splines (PHS): } \phi(r)=&r^{2a+1},\hspace{0.1cm} a \in \mathbb{N}\\
\text{Thin Plate Splines (TPS): } \phi(r)=&r^{2a} log(r),\hspace{0.1cm} a \in \mathbb{N}\\
\end{aligned}
\label{Eq:RBF_list}
\end{equation}
\par The PHS-RBF interpolation has the advantage over MQ, IMQ and Gaussians of not requiring a shape factor. The value of this shape factor is not known precisely, but the accuracy of interpolation as well as the condition number of the matrix required to determine the interpolation coefficients strongly depend on the shape factor. The RBF interpolations can also be appended with polynomials to get high accuracy. It has been demonstrated earlier \cite{flyer2016onrole_I, shahane2020high, radhakrishnan2021non} that when the RBFs are appended with polynomials and the derivatives of a smooth function are evaluated, the discretization errors converge exponentially per the degree of the appended polynomial. RBFs have been originally used as global interpolants of scattered data \cite{hardy1971multiquadric} and to solve partial differential equations \cite{kansa1990multiquadrics_I, kansa1990multiquadrics_II, kansa2000circumventing}, but the size of the problem has been limited because of the condition number of the solution matrix. Cloud based methods \cite{shu2003local, ding2006numerical}, on the other hand, are attractive because of the lower condition number and sparsity of the coefficient matrix. PHS-RBF methods with appended polynomials have been used in the past to solve several partial differential equations, including the Navier Stokes equations. Recently, we presented \cite{shahane2020high} a pressure-projection based algorithm using PHS-RBF to solve the time-dependent Navier Stokes equations of an incompressible flow. We presented the spatial discretization accuracy of the method in several problems by systematically varying the appended polynomial and the numbers of scattered points. That algorithm used an explicit formulation of the advection and diffusion terms using an Adams-Bashforth discretization in time. The divergence-free condition on the velocities was imposed by solving a pressure Poisson equation. Such explicit methods are well-suited for problems in which rapid transients are to be resolved, such as in simulations of turbulence. However, they are not optimal to compute steady state flows, and time-varying flows which permit large timesteps from accuracy viewpoint. For such needs, semi-implicit methods that permit several multiples of the velocity Courant number can be faster in solution time. Semi-implicit methods iterate between the momentum and continuity equations until desired convergence is obtained at each time step. Although the computational cost per timestep increases due to the iterations, overall this approach is efficient because the number of timesteps to reach a stationary solution is reduced significantly.
\par In this paper, we describe a semi-implicit PHS-RBF meshless method which is attractive to compute steady as well as temporally evolving flows. The method can use time steps significantly larger than the time step restricted by an explicit algorithm. We first present the convergence characteristics of the method in several model flows. Subsequently, the algorithm is applied to compute three complex steady and unsteady recirculating flows.

\section{The PHS-RBF Discretization} \label{Sec:The PHS-RBF Method}
The Navier-Stokes equations consist of differential operators of first and second order. In this work, PHS-RBF expansions collocated over a cloud of neighboring points are used to numerically estimate these differential operators. Monomials upto a maximum degree $k$ are appended for stabilization and accuracy of the PHS functions. A continuous scalar variable $s$ can be interpolated over $q$ cloud points as follows:
\begin{equation}
s(\textbf{x}) = \sum_{i=1}^{q} \lambda_i \phi_i (||\bm{x} - \bm{x_i}||_2) + \sum_{i=1}^{m} \gamma_i P_i (\bm{x})
\label{Eq:RBF_interp}
\end{equation}
where, $\phi(r)=r^{2a+1},\hspace{0.1cm} a \in \mathbb{N}$. In this work, we use the PHS-RBF $\phi(r)=r^3$. There are total $q+m$ unknowns ($\lambda_i$ and $\gamma_i$) to be estimated. Collocating \cref{Eq:RBF_interp} at the $q$ cloud points gives $q$ conditions. Additional $m$ equations are obtained by imposing the following constraints on the polynomials \cite{flyer2016onrole_I}:
\begin{equation}
\sum_{i=1}^{q} \lambda_i P_j(\bm{x_i}) =0 \hspace{0.5cm} \text{for } 1 \leq j \leq m
\label{Eq:RBF_constraint}
\end{equation}
The collocation conditions together with \cref{Eq:RBF_constraint} can be written in a matrix vector form,
\begin{equation}
\begin{bmatrix}
\bm{\Phi} & \bm{P}  \\
\bm{P}^T & \bm{0} \\
\end{bmatrix}
\begin{bmatrix}
\bm{\lambda}  \\
\bm{\gamma} \\
\end{bmatrix} =
\begin{bmatrix}
\bm{A}
\end{bmatrix}
\begin{bmatrix}
\bm{\lambda}  \\
\bm{\gamma} \\
\end{bmatrix} =
\begin{bmatrix}
\bm{s}  \\
\bm{0} \\
\end{bmatrix}
\label{Eq:RBF_interp_mat_vec}
\end{equation}
where, the superscript $T$ denotes the transpose, $\bm{\lambda} = [\lambda_1,...,\lambda_q]^T$, $\bm{\gamma} = [\gamma_1,...,\gamma_m]^T$, $\bm{s} = [s(\bm{x_1}),...,s(\bm{x_q})]^T$ and $\bm{0}$ is the matrix of zeros. Sizes of the submatrices $\bm{P}$ and $\bm{\Phi}$ are $q\times m$ and $q\times q$ respectively. The submatrix $\bm{\Phi}$ is given by:
\begin{equation}
\bm{\Phi} =
\begin{bmatrix}
\phi \left(||\bm{x_1} - \bm{x_1}||_2\right) & \dots  & \phi \left(||\bm{x_1} - \bm{x_q}||_2\right) \\
\vdots & \ddots & \vdots \\
\phi \left(||\bm{x_q} - \bm{x_1}||_2\right) & \dots  & \phi \left(||\bm{x_q} - \bm{x_q}||_2\right) \\
\end{bmatrix}
\label{Eq:RBF_interp_phi}
\end{equation}
For two dimensional problem ($d=2$) with maximum degree of appended polynomial set to 2 ($k=2$), there are $m=\binom{k+d}{k}=\binom{2+2}{2}=6$ polynomial terms: $[1, x, y, x^2, xy, y^2]$. Thus, the submatrix $\bm{P}$ is obtained by evaluating the polynomial terms at the $q$ cloud points.
\begin{equation}
\bm{P} =
\begin{bmatrix}
1 & x_1  & y_1 & x_1^2 & x_1 y_1 & y_1^2 \\
\vdots & \vdots & \vdots & \vdots & \vdots & \vdots \\
1 & x_q  & y_q & x_q^2 & x_q y_q & y_q^2 \\
\end{bmatrix}
\label{Eq:RBF_interp_poly}
\end{equation}
Unknown interpolation coefficients $\lambda_i$ and $\gamma_i$ are obtained by solving the \cref{Eq:RBF_interp_mat_vec}:
\begin{equation}
\begin{bmatrix}
\bm{\lambda}  \\
\bm{\gamma} \\
\end{bmatrix} =
\begin{bmatrix}
\bm{A}
\end{bmatrix} ^{-1}
\begin{bmatrix}
\bm{s}  \\
\bm{0} \\
\end{bmatrix}
\label{Eq:RBF_interp_mat_vec_solve}
\end{equation}
In \cref{Eq:RBF_interp_mat_vec_solve}, the inverse of matrix $\bm{A}$ is used as a notation. Practically, the equation is solved by a direct linear solver.
\par The differential operators in the Navier-Stokes equation are estimated as a weighted linear summation of the function values at the neighboring cloud points. Applying any scalar linear operator $\mathcal{L}$ such as $\frac{\partial}{\partial x}$ or the Laplacian $\nabla ^2$ on \cref{Eq:RBF_interp} gives:
\begin{equation}
\mathcal{L} [s(\textbf{x})] = \sum_{i=1}^{q} \lambda_i \mathcal{L} [\phi_i (\bm{||\bm{x} - \bm{x_i}||_2})] + \sum_{i=1}^{m} \gamma_i \mathcal{L}[P_i (\bm{x})]
\label{Eq:RBF_interp_L}
\end{equation}
\Cref{Eq:RBF_interp_L} is collocated by evaluating $\mathcal{L}[\bm{s}]$ at the $q$ cloud points. This gives a rectangular matrix vector system:
\begin{equation}
\mathcal{L}[\bm{s}] =
\begin{bmatrix}
\mathcal{L}[\bm{\Phi}] & \mathcal{L}[\bm{P}]  \\
\end{bmatrix}
\begin{bmatrix}
\bm{\lambda}  \\
\bm{\gamma} \\
\end{bmatrix}
\label{Eq:RBF_interp_mat_vec_L}
\end{equation}
where, $\mathcal{L}[\bm{P}]$ and $\mathcal{L}[\bm{\Phi}]$ are matrices of sizes $q\times m$ and $q\times q$ respectively. Substituting \cref{Eq:RBF_interp_mat_vec_solve} in \cref{Eq:RBF_interp_mat_vec_L} and simplifying, we get:
\begin{equation}
\begin{aligned}
\mathcal{L}[\bm{s}] &=
\left(\begin{bmatrix}
\mathcal{L}[\bm{\Phi}] & \mathcal{L}[\bm{P}]  \\
\end{bmatrix}
\begin{bmatrix}
\bm{A}
\end{bmatrix} ^{-1}\right)
\begin{bmatrix}
\bm{s}  \\
\bm{0} \\
\end{bmatrix}
=
\begin{bmatrix}
\bm{B}
\end{bmatrix}
\begin{bmatrix}
\bm{s}  \\
\bm{0} \\
\end{bmatrix}\\
&=
\begin{bmatrix}
\bm{B_1} & \bm{B_2}
\end{bmatrix}
\begin{bmatrix}
\bm{s}  \\
\bm{0} \\
\end{bmatrix}
= [\bm{B_1}] [\bm{s}] + [\bm{B_2}] [\bm{0}]
= [\bm{B_1}] [\bm{s}]
\end{aligned}
\label{Eq:RBF_interp_mat_vec_L_solve}
\end{equation}
Thus, $\mathcal{L}[\bm{s}]$ at any point is estimated using the weighted summation of the values of $\bm{s}$ at its neighboring cloud points. The matrix $[\bm{B_1}]$ in \cref{Eq:RBF_interp_mat_vec_L_solve} is precomputed using the coordinates of the cloud points and stored for use throughout the solution steps. Properties of the PHS-RBF expansion are reported in our previous publication \cite{shahane2020high}.

\section{Semi-Implicit Solution Algorithm} \label{Sec:Semi-Implicit Algorithm}
We have previously reported a fully explicit fractional step algorithm with the PHS-RBF method to solve the incompressible Navier-Stokes equations \cite{shahane2020high}. This explicit method is computationally expensive since the time step is limited by the condition that the convection and diffusion Courant number be less than unity. In this work, we describe a semi-implicit iterative algorithm which is shown to be stable for Courant numbers much higher than unity. For temporal discretization, we use a second-order accurate backward differentiation formula (BDF2) \cite{suli2003introduction} which expresses the momentum equations as

\begin{equation}
\frac{\rho (\alpha_1 u^{n+1} + \alpha_2 u^n + \alpha_3 u^{n-1})}{\Delta t} = -\rho \bm{u}^{n+1} \bullet (\nabla u^{n+1}) + \mu \nabla^2 u^{n+1} - \left(\frac{\partial p}{\partial x}\right)^{n+1}
\label{Eq:u:n+1}
\end{equation}
\begin{equation}
\frac{\rho (\alpha_1 v^{n+1} + \alpha_2 v^n + \alpha_3 v^{n-1})}{\Delta t} = -\rho \bm{u}^{n+1} \bullet (\nabla v^{n+1}) + \mu \nabla^2 v^{n+1} - \left(\frac{\partial p}{\partial y}\right)^{n+1}
\label{Eq:v:n+1}
\end{equation}
where, $n$ is timestep number, $\rho$ is density, $\mu$ is dynamic viscosity, $\bm{u}=[u,v]$ is velocity vector, $p$ is pressure and $\alpha_1=1.5$, $\alpha_2=-2$ and $\alpha_3=0.5$ are the BDF2 coefficients \cite{suli2003introduction}. The incompressible continuity equation is evaluated at the new timestep $n+1$ as
\begin{equation}
\left(\frac{\partial u}{\partial x}\right)^{n+1} + \left(\frac{\partial v}{\partial y}\right)^{n+1} = 0
\label{Eq:continuity:n+1}
\end{equation}
\Cref{Eq:u:n+1,Eq:v:n+1} have to be solved iteratively since the convection terms are nonlinear and $p^{n+1}$ is unknown. Given the $[u^n, v^n, p^n]$ and $[u^{n-1}, v^{n-1}, p^{n-1}]$ values at the previous two timesteps, the algorithm iterates to solve for $[u^{n+1}, v^{n+1}, p^{n+1}]$ as follows. First, the advection terms are linearized with previous iterate values and the diffusion terms are expressed implicitly at the current iteration. Together with the old iterate pressure gradients, we solve: 
\begin{equation}
\frac{\rho (\alpha_1 \widetilde{u} + \alpha_2 u^n + \alpha_3 u^{n-1})}{\Delta t} = -\rho \bm{u}^r \bullet (\nabla \widetilde{u}) + \mu \nabla^2 \widetilde{u} - \left(\frac{\partial p}{\partial x}\right)^r
\label{Eq:u:tilde}
\end{equation}
\begin{equation}
\frac{\rho (\alpha_1 \widetilde{v} + \alpha_2 v^n + \alpha_3 v^{n-1})}{\Delta t} = -\rho \bm{u}^r \bullet (\nabla \widetilde{v}) + \mu \nabla^2 \widetilde{v} - \left(\frac{\partial p}{\partial y}\right)^r
\label{Eq:v:tilde}
\end{equation}
where, $\widetilde{u}$ and $\widetilde{v}$ denote intermediate velocity fields. The spatial derivatives are expressed as a weighted linear combination of the velocities and pressure at the cloud points as discussed in \cref{Sec:The PHS-RBF Method}. \Cref{Eq:u:tilde,Eq:v:tilde} are solved for $\widetilde{u}$ and $\widetilde{v}$ using a preconditioned BiCGSTAB algorithm. We then correct the pressure field by approximately relating the velocity corrections to the corrections of pressure gradients.
\begin{equation}
\frac{\rho \alpha_1 (u^{r+1} - \widetilde{u})}{\Delta t} \approx - \frac{\partial (p^{r+1} - p^r)}{\partial x} = - \frac{\partial p'}{\partial x}
\label{Eq:u:p'}
\end{equation}
\begin{equation}
\frac{\rho \alpha_1 (v^{r+1} - \widetilde{v})}{\Delta t} \approx - \frac{\partial (p^{r+1} - p^r)}{\partial y} = - \frac{\partial p'}{\partial y}
\label{Eq:v:p'}
\end{equation}
where, $p'$ denotes pressure corrections. Note that in the above approximation, we have neglected the changes in the advection and diffusion terms. Imposing the divergence-free condition on the velocity field at the new iteration ($r+1$) gives the Poisson equation for pressure corrections:
\begin{equation}
\nabla^2 p' = \frac{\rho \alpha_1}{\Delta t}\left(\frac{\partial \widetilde{u}}{\partial x} + \frac{\partial \widetilde{v}}{\partial y}\right)
\label{Eq:pressure poisson}
\end{equation}
The pressure-correction equation is also discretized with high accuracy using the PHS-RBF discretization scheme and solved with a preconditioned BiCGSTAB algorithm to a pre-specified tolerance. After the pressure-corrections are determined, the velocities and pressures are corrected as
\begin{equation}
u^{r+1} = \widetilde{u} - \frac{\Delta t}{\rho \alpha_1} \left(\frac{\partial p'}{\partial x}\right)
\label{Eq:u:correction}
\end{equation}
\begin{equation}
v^{r+1} = \widetilde{v} - \frac{\Delta t}{\rho \alpha_1} \left(\frac{\partial p'}{\partial y}\right)
\label{Eq:v:correction}
\end{equation}
Similarly, the pressure is corrected as:
\begin{equation}
p^{r+1} = p^r + p'
\label{Eq:p:correction}
\end{equation}
Note that the above algorithm has similarities to the SIMPLE series of algorithms \cite{patankar2018numerical}, but here the pressure correction equation is formulated only through the time-dependent term in the momentum equations whereas, SIMPLE formulates the $p’$ equation from truncated discretized momentum equations with only the diagonal terms. In situations where time-dependent problems are solved with SIMPLE, the contributions of the time-dependent terms are added to the diagonal coefficients. The present algorithm is more straight forward in deriving the pressure-correction equation. Also, we do not use any under-relaxation factors in any of the equations. The time-dependent term provides the needed damping to the nonlinearities and the coupling terms.
\subsection{Boundary Conditions}\label{Sec:Boundary Conditions}
The present algorithm is compatible with most of the boundary conditions used in practical flow conditions. These include Dirichlet, Neumann and periodic boundary conditions. At inlet and wall boundaries, the velocities are specified to be the known values. These velocities are either fixed in time or can be varied depending on the problem. The values of these velocities are also imposed on the $\widetilde{u}$ and $\widetilde{v}$ velocities of \cref{Eq:u:tilde,Eq:v:tilde}. At symmetry boundaries, the normal velocities are prescribed to be zero and the tangential velocities at the boundaries are determined by implicitly imposing zero normal derivatives. Derivatives at the boundary points are discretized over a local cloud surrounding the boundary points. The pressures at the boundary points are determined by imposing the normal momentum equation given by
\begin{equation}
\nabla p \bullet \widehat{\bm{N}} = \left[-\rho\left(\frac{\partial \bm{u}}{\partial t}\right) -\rho(\bm{u}\bullet \nabla)\bm{u} + \mu \nabla^2 \bm{u}\right] \bullet \widehat{\bm{N}}
\label{Eq:normal momentum}
\end{equation}
These boundary pressures are required for evaluating the pressure-gradients appearing in the momentum equations. In the solution of the pressure correction equation however, homogeneous Neumann conditions ($\partial p'/ \partial \widehat{\bm{N}} = 0$) are imposed at the walls, inlet, and symmetry boundaries. For boundaries with prescribed pressure, the pressure corrections are specified to be zero. The outlet boundary can be prescribed to be either fixed pressure or have zero normal derivatives on the velocities. The relevant derivatives are computed in Cartesian coordinates and resolved appropriately in the direction  normal to the boundary. Note that only one-sided cloud points are used at all boundaries.
\par
\begin{figure}[H]
	\centering
	\includegraphics[width=0.75\textwidth]{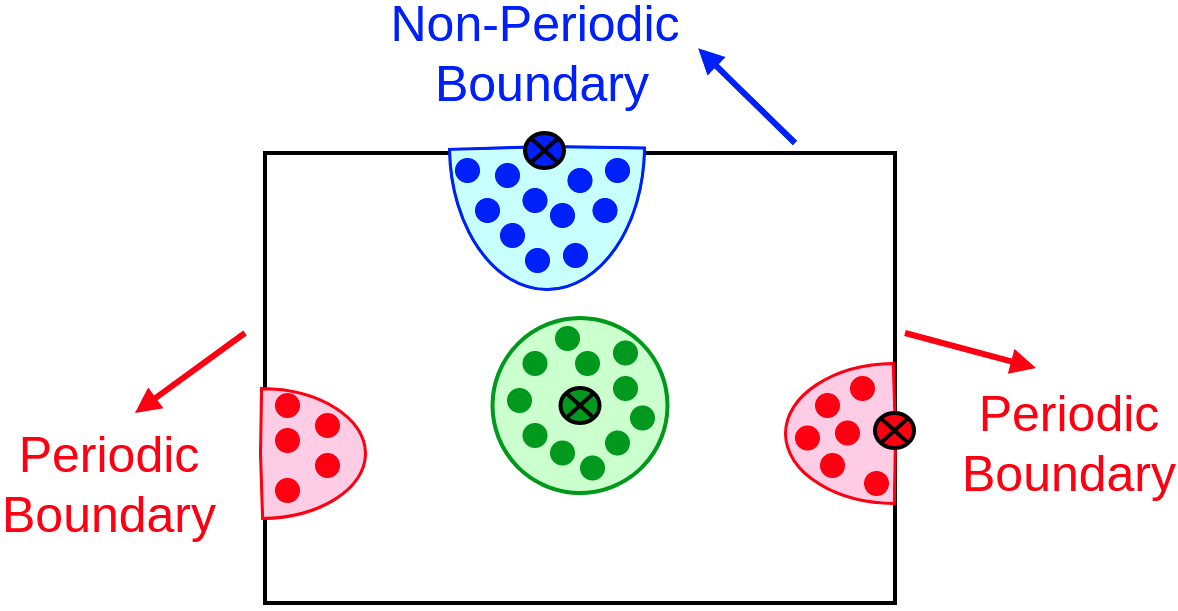}
	\caption{Clouds for Interior (Green), Periodic (Red) and Non-Periodic (Blue) Boundaries}
	\label{Fig:schematic cloud}
\end{figure}

 Lastly, we have also implemented periodicity boundary conditions in any one or both the spatial directions. Periodicity implies repeated domain with same velocities and pressures. Thus, the variables are wrapped around in the periodic direction. \Cref{Fig:schematic cloud} shows a schematic of the cloud definitions at the periodic and non-periodic boundaries. The periodic boundary is imposed implicitly by defining clouds for near boundary points that include some points adjacent to the boundary from the other end of the domain. Thus, for periodic boundary conditions, the clouds are much like any interior point shown in the center of \cref{Fig:schematic cloud}. This contrasts with the asymmetric cloud of points at a non-periodic boundary point. For periodic flows, such as a fully-developed flow in a channel, a mean pressure gradient is provided as an additional forcing term in the direction of periodicity. The governing equations are solved on the points near a periodic boundary and thus, no boundary condition is required. For this case, the two sides of the domain are coupled with each other. 

\subsection{Summary of Solution Steps}
The steps in the present solution algorithm for time marching from timesteps $n-1$ and $n$ to $n+1$ can be summarized as follows:
\begin{enumerate}
	\item Initialize the values at the iteration number $r=0$ as $[u^r, v^r, p^r] = [u^n, v^n, p^n]$
	\item Solve \cref{Eq:u:tilde,Eq:v:tilde} for intermediate velocities ($\widetilde{u}$ and $\widetilde{v}$) with the appropriate boundary condition  (\cref{Sec:Boundary Conditions}) using any linear solver such as a preconditioned BiCGSTAB solver
	\item Solve the Poisson \cref{Eq:pressure poisson} for pressure correction ($p'$)
	\item Correct the pressure to estimate $p^{r+1}$ at the interior points (\cref{Eq:p:correction})
	\item Set the pressure at the boundary points using the normal momentum \cref{Eq:normal momentum}
	\item Update the velocities at the new iteration ($u^{r+1}$ and $v^{r+1}$) using the velocity correction \cref{Eq:u:correction,Eq:v:correction}
	\item Compute the change during the current iteration: $\Delta = ||u^{r+1}-u^r|| + ||v^{r+1}-v^r||$
	\item If $\Delta$ is less than the prescribed iterative tolerance, set $[u^{n+1}, v^{n+1}, p^{n+1}] = [u^{r+1}, v^{r+1}, p^{r+1}]$ and proceed to a new timestep. Otherwise, go back to step 2
\end{enumerate}
We calculate the Courant number from the eigenvalues of the matrix corresponding to the discrete convection and diffusion operators. The Courant number given by $\lambda_{rmax} \Delta t/2$ is kept to be less than unity for first order explicit Euler method where, $\lambda_{rmax}$ denotes the eigenvalue with the maximum real component \cite{leveque2007finite}. In the semi-implicit method with BDF2, we can have Courant numbers much greater than unity. We have successfully tested the semi-implicit algorithm with BDF2 for time integration upto a Courant number of 12. For the problems studied in this paper, we observed that at each time step typically 2-3 iterations have been necessary to obtain a good convergent solution. However, the first few timesteps (around 20) may require 5 or 10 outer iterations. We initially set an iterative tolerance of 1E--4. For subsequent timesteps, the tolerance is updated as: $\min[1\text{E--4}, \Delta t \epsilon_s]$ where, $\epsilon_s=\frac{\phi^{n+1} - \phi^n}{\Delta t}$ is the rate of change in the velocity fields at consecutive timesteps denoted by $\phi^{n+1}$ and $\phi^n$. As steady state is approached $\epsilon_s$ reduces and thus, the term $\Delta t \epsilon_s$ tightens the iterative tolerance. 

\section{Convergence Characteristics in Model Problems}
\subsection{Steady Cylindrical Couette Flow} \label{Sec:couette flow}
\begin{figure}[H]
	\centering
	\begin{subfigure}[t]{0.45\textwidth}
		\includegraphics[width=\textwidth]{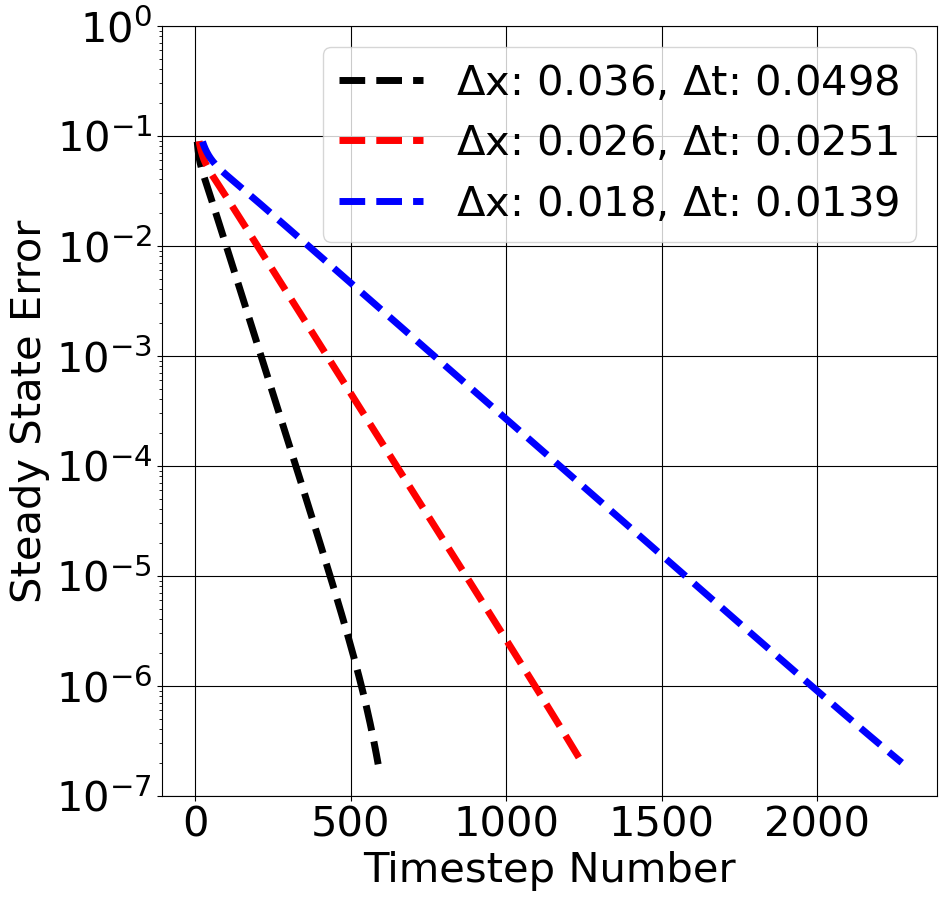}
		\caption{Polynomial Degree 3}
	\end{subfigure}
	\hspace{0.05\textwidth}
	\begin{subfigure}[t]{0.45\textwidth}
		\includegraphics[width=\textwidth]{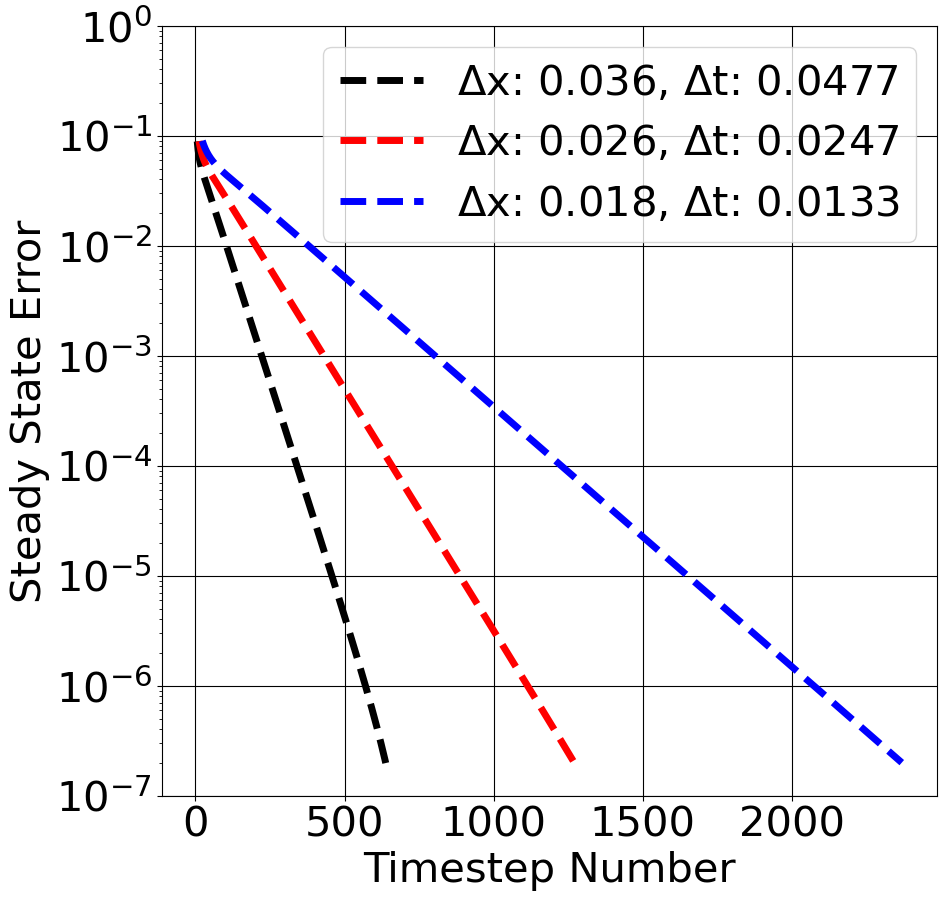}
		\caption{Polynomial Degree 4} \vspace{0.25cm}
	\end{subfigure}
	\begin{subfigure}[t]{0.45\textwidth}
		\includegraphics[width=\textwidth]{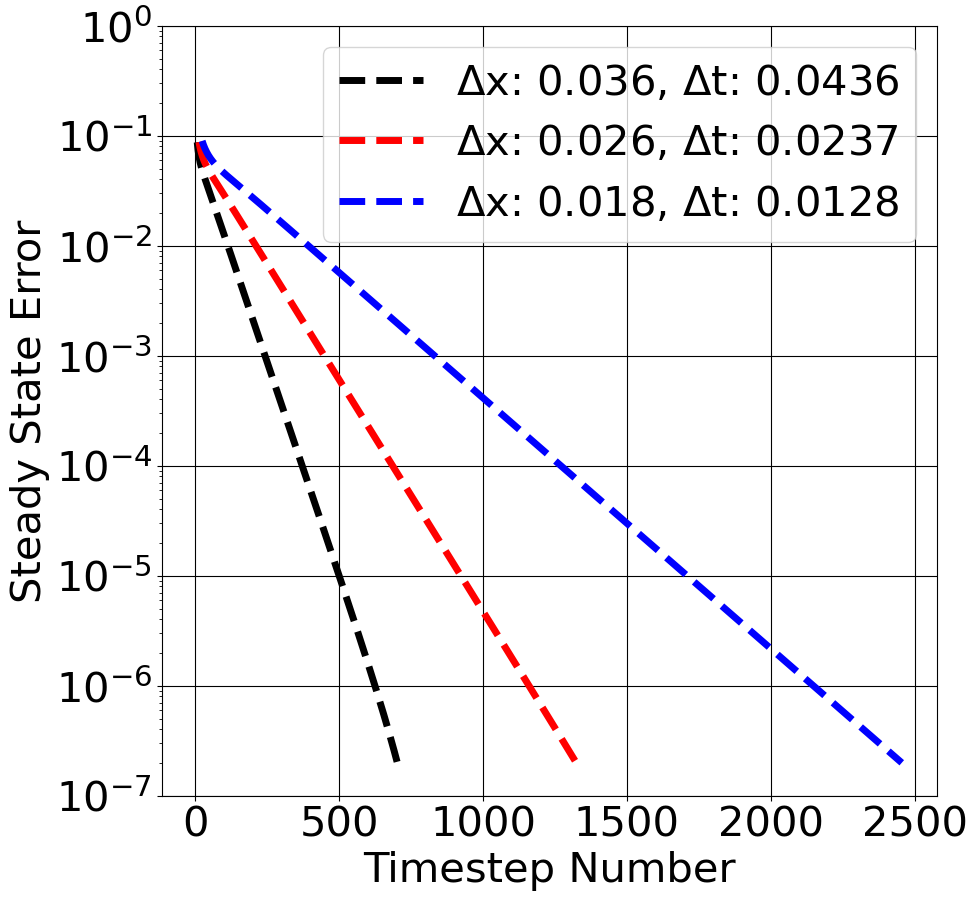}
		\caption{Polynomial Degree 5}
	\end{subfigure}
	\hspace{0.05\textwidth}
	\begin{subfigure}[t]{0.45\textwidth}
		\includegraphics[width=\textwidth]{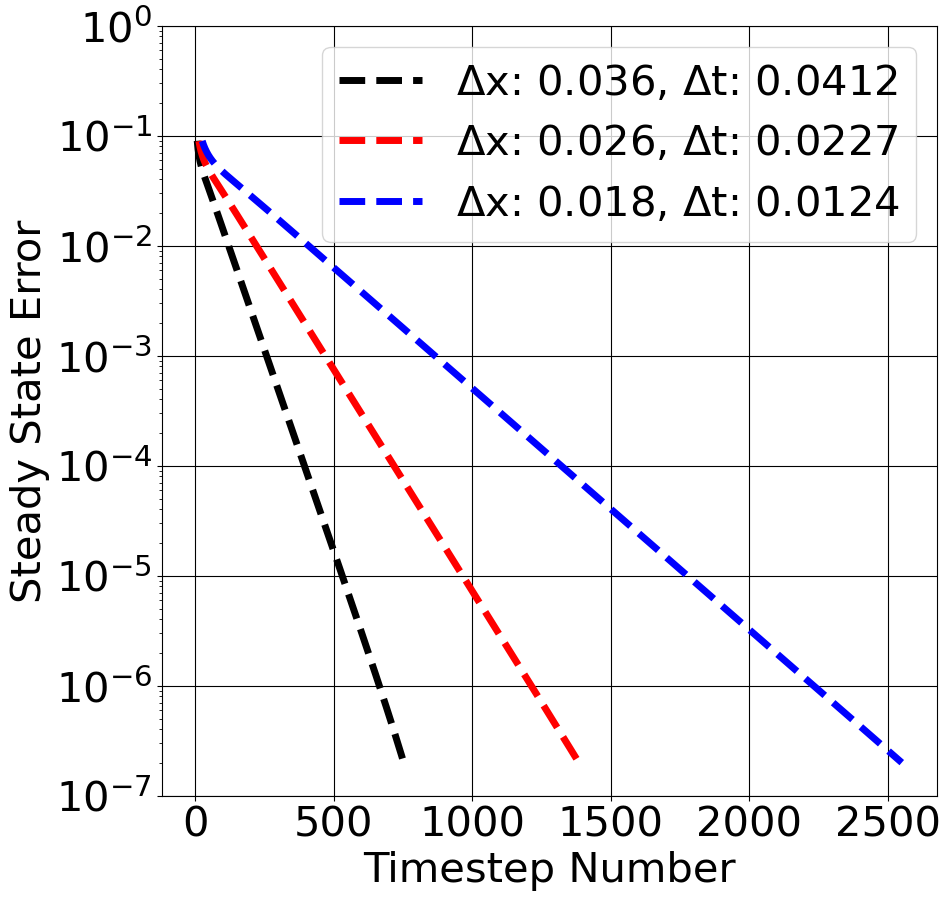}
		\caption{Polynomial Degree 6}
	\end{subfigure}
	\caption{Convergence of the Semi-Implicit Algorithm: Courant No.: 8}
	\label{Fig:couette steady error Co 8}
\end{figure}

\begin{figure}[H]
	\centering
	\begin{subfigure}[t]{0.45\textwidth}
		\includegraphics[width=\textwidth]{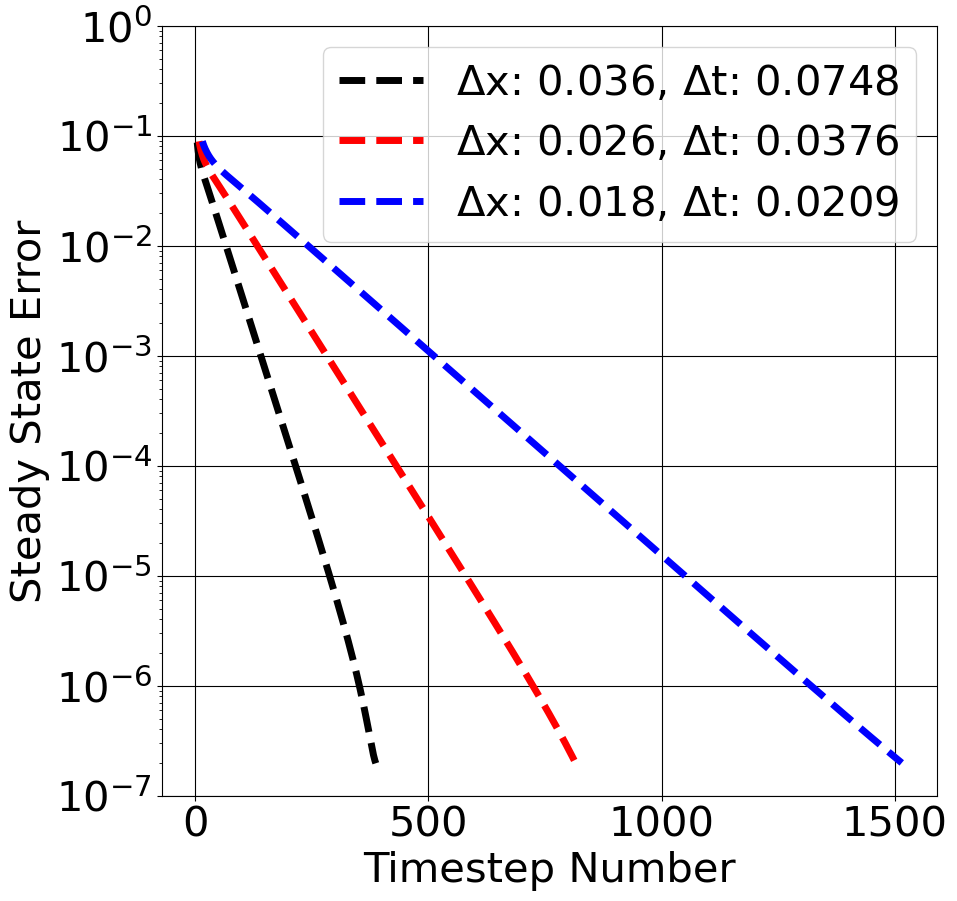}
		\caption{Polynomial Degree 3}
	\end{subfigure}
	\hspace{0.05\textwidth}
	\begin{subfigure}[t]{0.45\textwidth}
		\includegraphics[width=\textwidth]{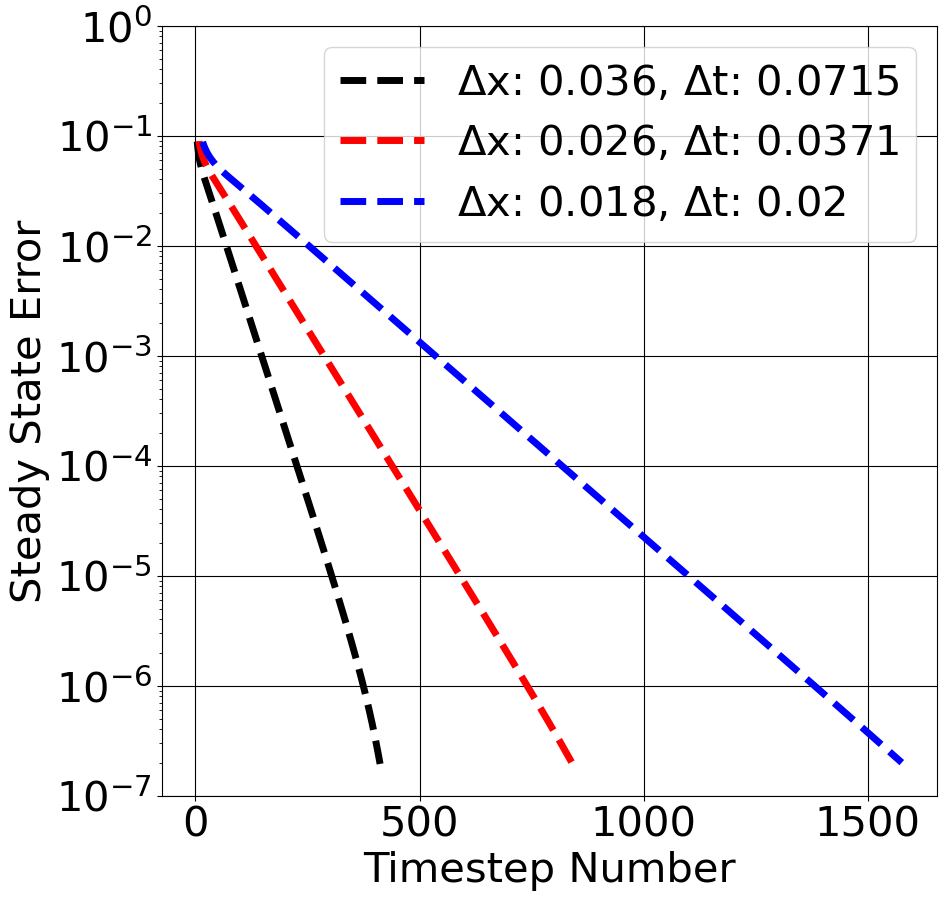}
		\caption{Polynomial Degree 4} \vspace{0.25cm}
	\end{subfigure}
	\begin{subfigure}[t]{0.45\textwidth}
		\includegraphics[width=\textwidth]{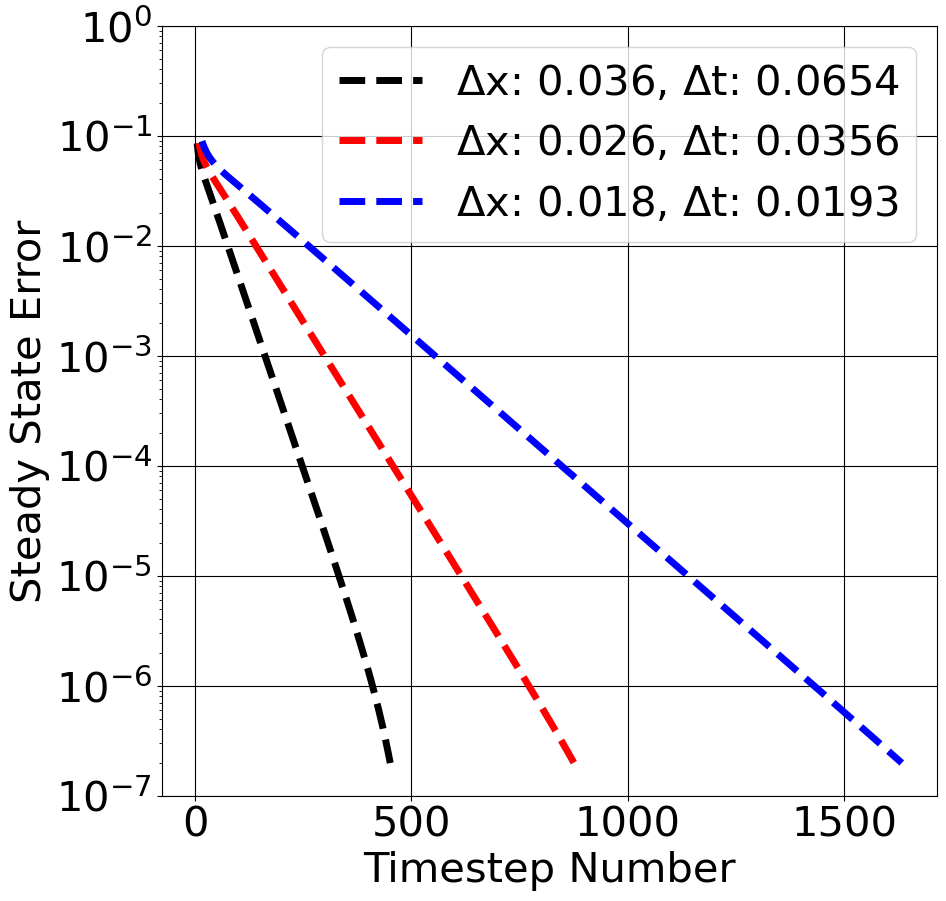}
		\caption{Polynomial Degree 5}
	\end{subfigure}
	\hspace{0.05\textwidth}
	\begin{subfigure}[t]{0.45\textwidth}
		\includegraphics[width=\textwidth]{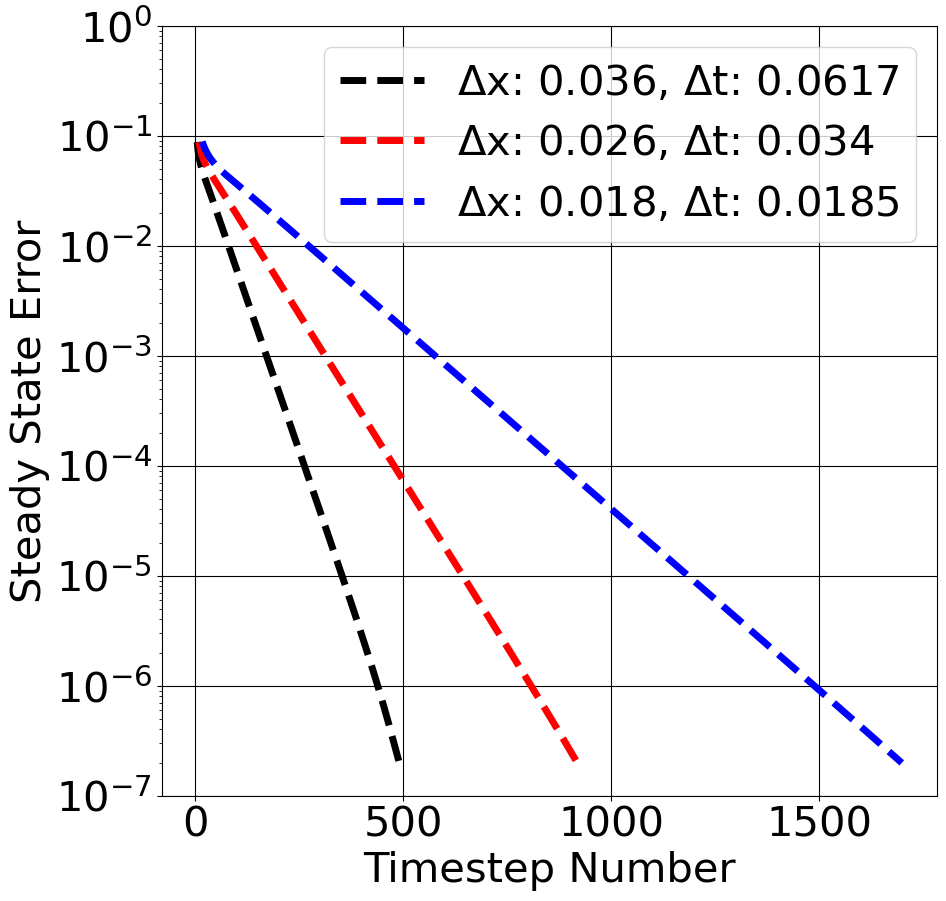}
		\caption{Polynomial Degree 6}
	\end{subfigure}
	\caption{Convergence of the Semi-Implicit Algorithm: Courant No.: 12}
	\label{Fig:couette steady error Co 12}
\end{figure}

We have first studied the convergence characteristics in cylindrical Couette flow to test the convergence of the algorithm as well as the discretization errors at steady state. Cylindrical Couette flow is a one dimensional problem in cylindrical polar coordinates but here, it is solved as a two dimensional problem in Cartesian coordinates. The annular region between the two cylinders is filled with fluid with inner cylinder rotating and outer cylinder stationary. The tangential velocity as a function of radial coordinate ($r$) is analytically given by:
\begin{equation}
v_{\theta}(r) = r_1 \omega \frac{r_1 r_2}{r_2^2 - r_1^2} \left(\frac{r_2}{r} - \frac{r}{r_2}\right)
\label{Eq:couette v_theta}
\end{equation}
where, $\omega$, $r_1$ and $r_2$ are rotational velocity of the inner cylinder, inner and outer radii respectively. This tangential velocity is decomposed into velocity fields in the Cartesian X and Y directions. At low rotational Reynolds number, the flow does not develop into longitudinal Taylor vortices and it remains steady and laminar. The flow Reynolds number defined with respect to inner cylinder's diameter and its tangential velocity is set to 100 although, the velocity field is independent of the Reynolds number. We have prescribed a rotational velocity of $\omega=2$, $r_1=0.5$ and $r_2=1$. Three different points distributions having 2080, 4033 and 8078 points are used with four increasing degrees of appended polynomials 3, 4, 5 and 6. These point distributions correspond to $\Delta x$ of 0.036, 0.026 and 0.018 respectively. $\Delta x$ is defined as the average value of the distance between each point and its nearest neighbor. Steady state solution is obtained by integrating in time from null velocity and pressure fields. Steady state error is defined as $\frac{\phi^{n+1} - \phi^n}{\Delta t}$ where, $\phi$ denotes X and Y components of velocity evaluated at consecutive timesteps $n$ and $n+1$. When the spatial average of the above error reaches a tolerance of 2E--7 for X and Y velocity components, the time integration is stopped.
\begin{figure}[H]
	\centering
	\begin{subfigure}[t]{0.45\textwidth}
		\includegraphics[width=\textwidth]{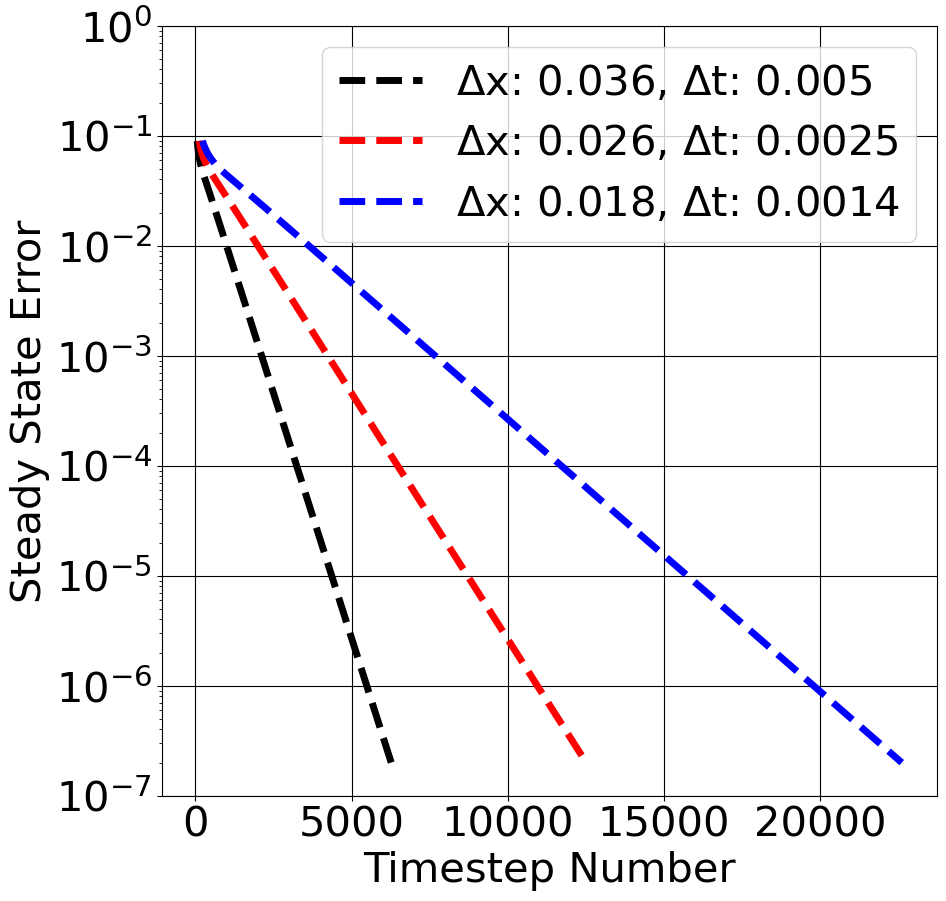}
		\caption{Polynomial Degree 3}
	\end{subfigure}
	\hspace{0.05\textwidth}
	\begin{subfigure}[t]{0.45\textwidth}
		\includegraphics[width=\textwidth]{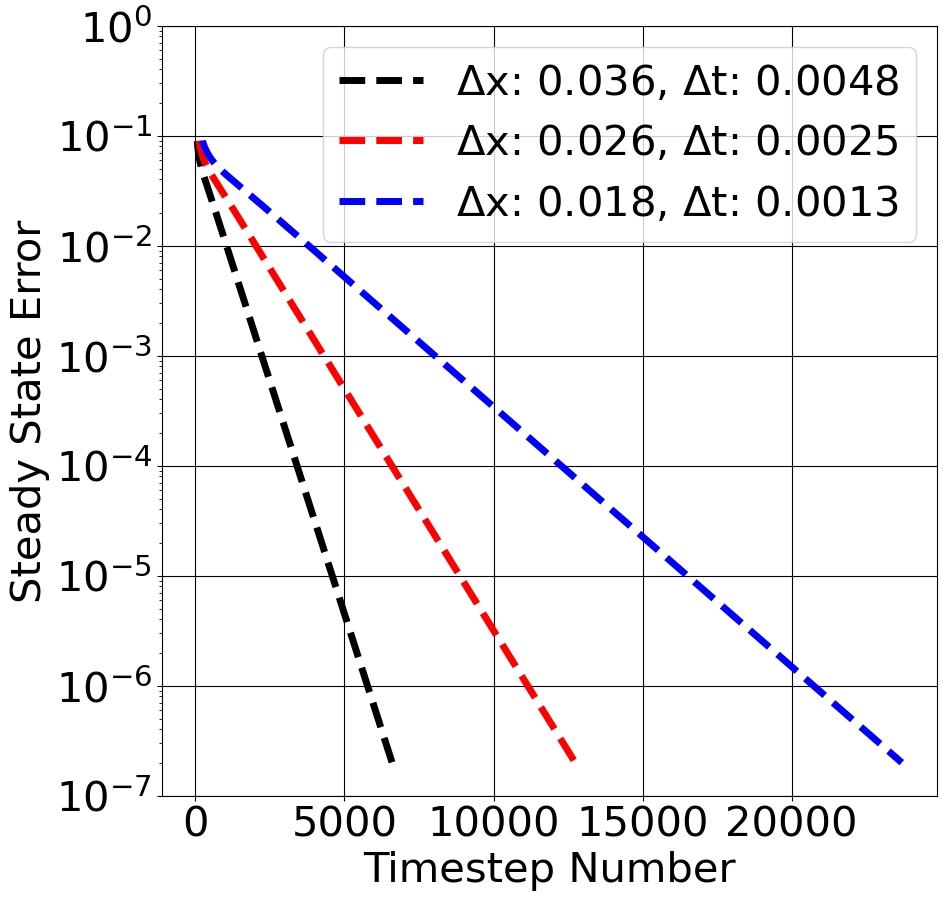}
		\caption{Polynomial Degree 4} \vspace{0.25cm}
	\end{subfigure}
	\begin{subfigure}[t]{0.45\textwidth}
		\includegraphics[width=\textwidth]{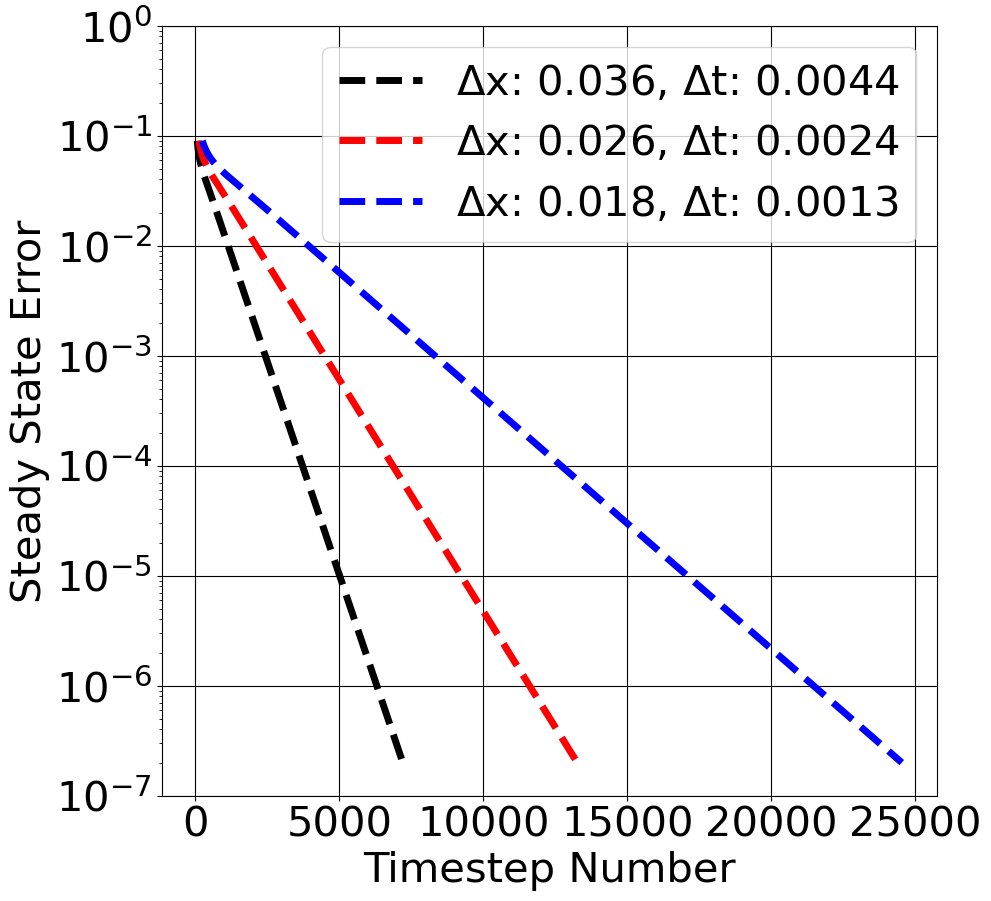}
		\caption{Polynomial Degree 5}
	\end{subfigure}
	\hspace{0.05\textwidth}
	\begin{subfigure}[t]{0.45\textwidth}
		\includegraphics[width=\textwidth]{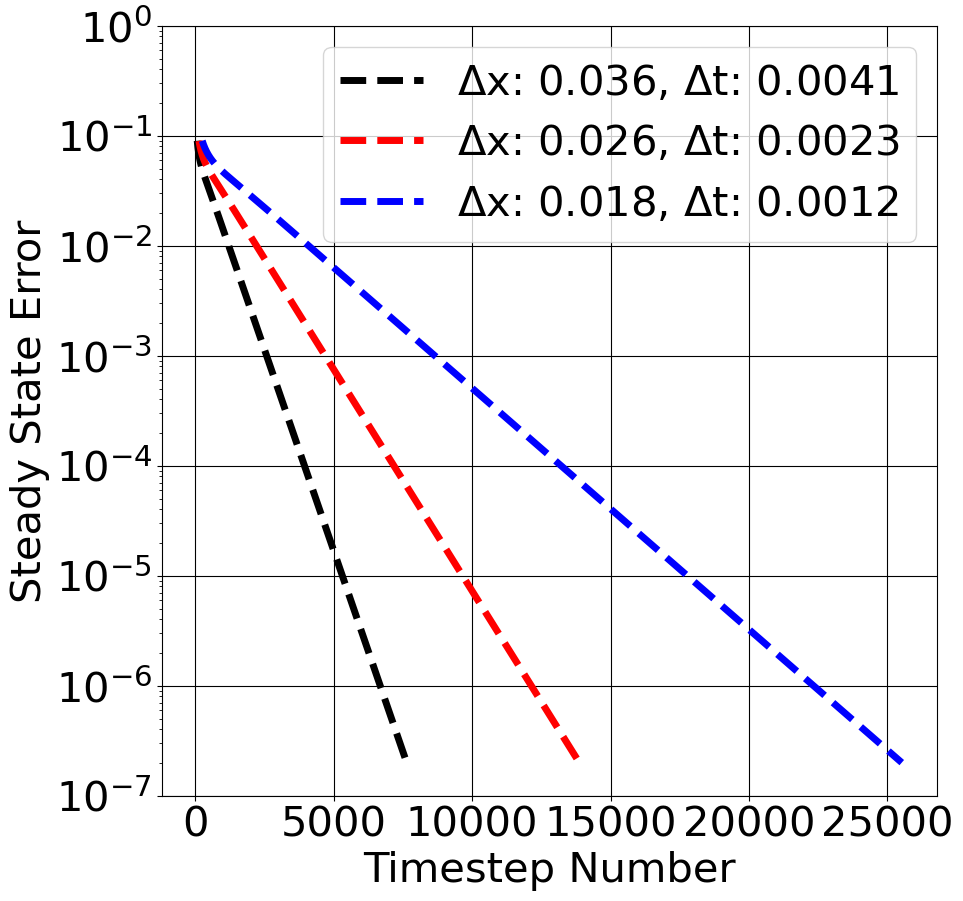}
		\caption{Polynomial Degree 6}
	\end{subfigure}
	\caption{Convergence of the Fracctional Step Algorithm \cite{shahane2020high}: Courant No.: 0.8}
	\label{Fig:couette steady error Co 0.8}
\end{figure}
For comparison, first we simulate using the fractional step method at a Courant number of 0.8 \cite{shahane2020high}. Since the stability is a function of the weighted sum of the eigenvalues of the discrete convection and diffusion operators, we see that when $\Delta x$ reduces by a factor of 2 from 0.036 to 0.018, the timestep $\Delta t$ reduces approximately by a factor of 3.5. The timestep only marginally reduces when the polynomial degree is increased for the same $\Delta x$. For problems having steady state, we have used first order explicit Euler method for time integration in the fractional step algorithm \cite{shahane2020high}. \Cref{Fig:couette steady error Co 0.8} plots the temporal evolution of the steady state error using the fractional step method. The steady state error drops at a constant rate on a logarithmic scale with timestep. The number of timesteps needed to reach steady state are inversely proportional to $\Delta t$ thus, coarser point distributions converge faster. \Cref{Fig:couette steady error Co 8,Fig:couette steady error Co 12} plot errors for the current semi-implicit algorithm with Courant numbers of 8 and 12 respectively. This case shows a steeper drop in the steady state errors compared to \cref{Fig:couette steady error Co 0.8}. Thus, the steady state is reached in less number of timesteps as expected. These plots show the computational advantage of the semi-implicit algorithm over the previous explicit fractional step method.

\par
\begin{figure}[H]
	\centering
	\includegraphics[width=0.7\textwidth]{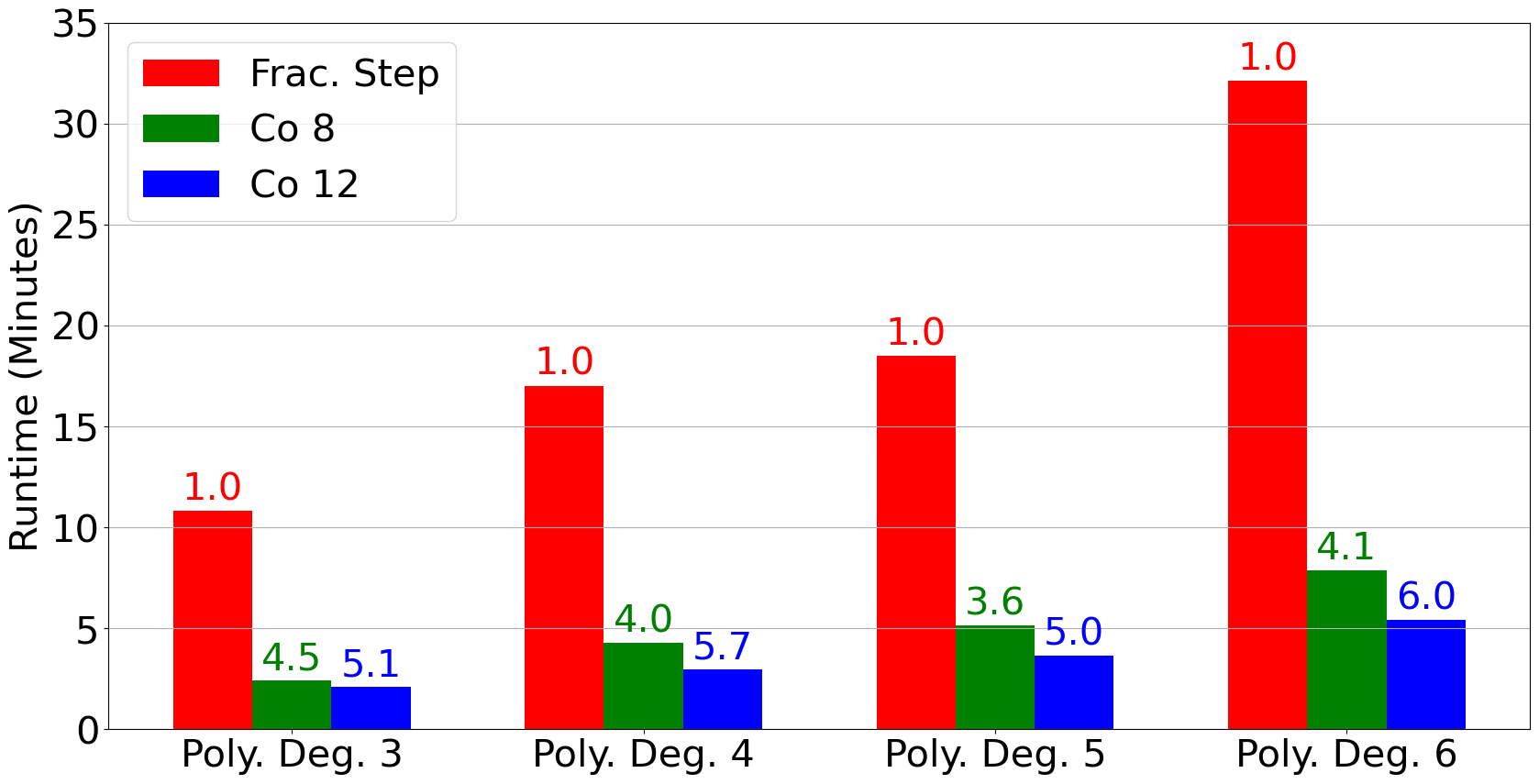}
	\caption{CPU Runtime with Speedup Factors for 8078 points with $\Delta x=0.018$}
	\label{Fig:couette runtimes}
\end{figure}

 \Cref{Fig:couette runtimes} compares CPU runtimes for the finest point distribution. The red colored bar plots runtime for the fractional step algorithm. The green and blue bars are for the semi-implicit algorithm with Courant numbers of 8 and 12 respectively. All calculations used the same compilers and computer hardware. For each polynomial degree, the speedup is defined as the ratio of runtime of the fractional step to the semi-implicit algorithm. Thus, the speedup for the fractional step is unity whereas, for the semi-implicit method, it is in the range [3.6, 4.5] and [5.0, 6.0] for the Courant numbers of 8 and 12 respectively.
\par For transient cases, we find that each timestep of the semi-implicit algorithm converges in 3 to 4 outer iterations for a tolerance of 1E--4. The fractional step algorithm with second-order accurate Adams Bashforth differencing requires a Courant number of 0.5 or less. Hence with a Courant number of 12, the semi-implicit algorithm is 6 to 8 times faster over the fractional step algorithm when transient problems are of interest. However, because of the larger time step in the semi-implicit algorithm the temporal errors can be larger than those in the fractional step algorithm, depending on the transient. Note that although we have used a maximum Courant number of 12 in this paper, it is not a hard limit.
We have observed the algorithm to be stable up to Courant numbers of 50 in some problems. However, at large Courant numbers there may be some minor oscillations in convergence. Hence, a value of Courant number around 12 is used for both steady and transient problems considered in this paper.
\par
\begin{figure}[H]
	\centering
	\begin{subfigure}[t]{0.45\textwidth}
		\includegraphics[width=\textwidth]{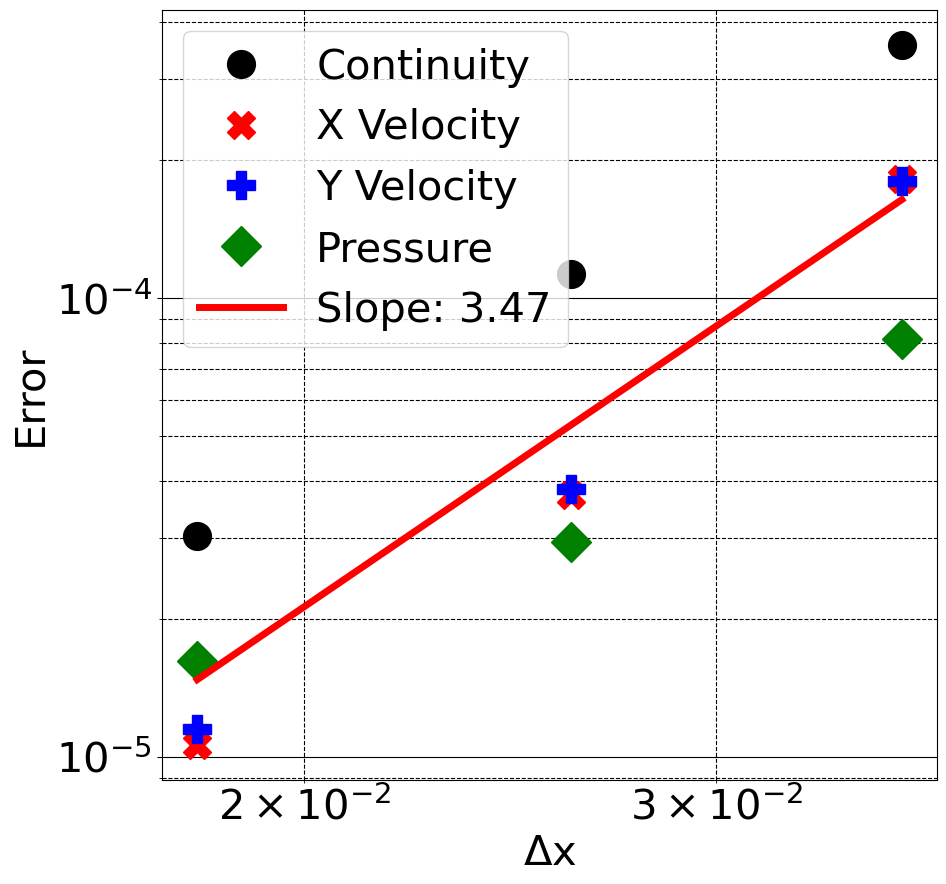}
		\caption{Polynomial Degree 3}
		\label{Fig:couette error courant 12 poldeg 3}
	\end{subfigure}
	\hspace{0.05\textwidth}
	\begin{subfigure}[t]{0.45\textwidth}
		\includegraphics[width=\textwidth]{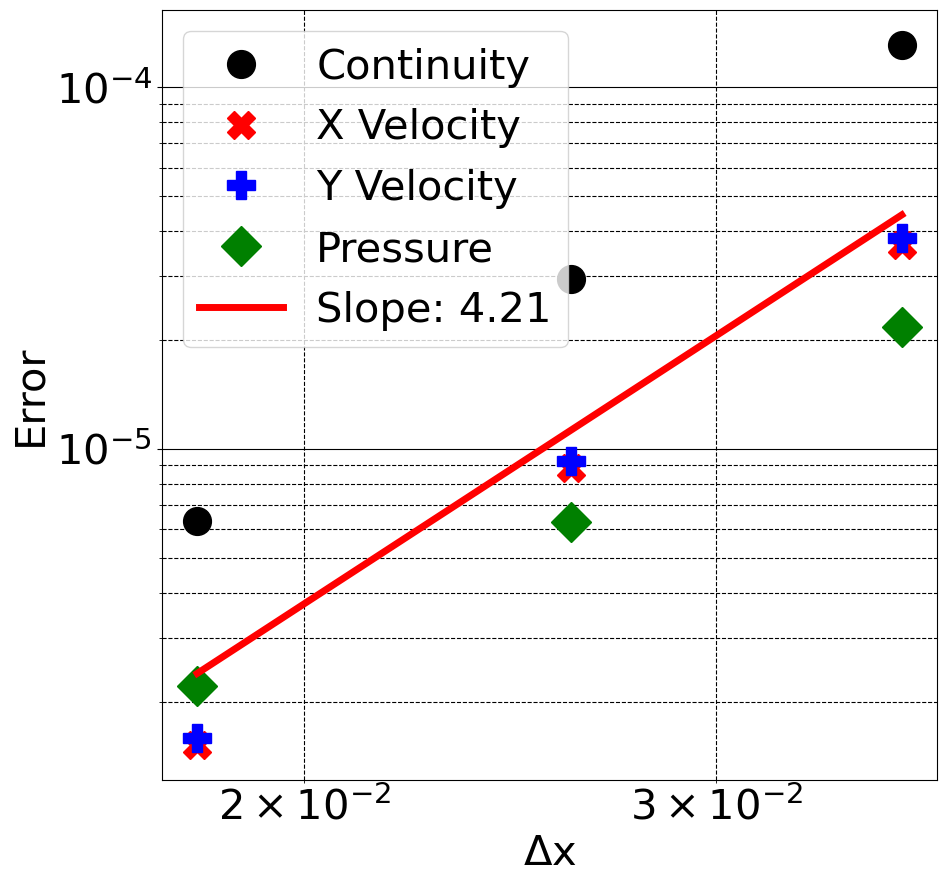}
		\caption{Polynomial Degree 4} \vspace{0.25cm}
	\end{subfigure}
	\begin{subfigure}[t]{0.45\textwidth}
		\includegraphics[width=\textwidth]{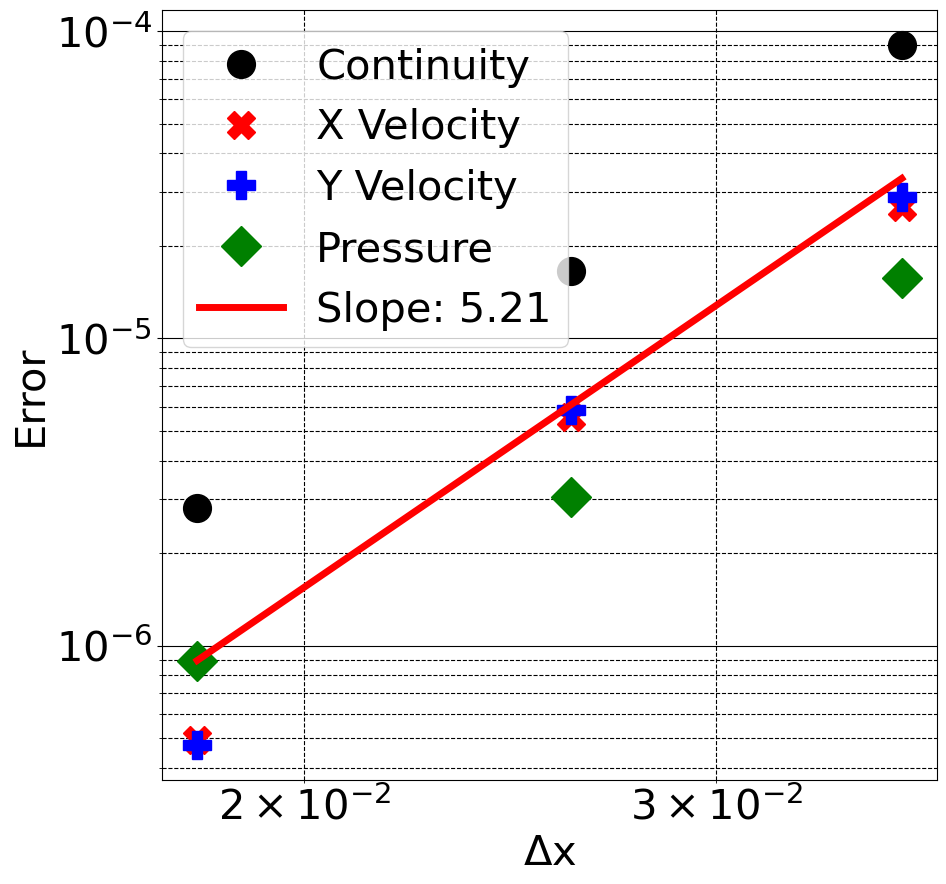}
		\caption{Polynomial Degree 5}
	\end{subfigure}
	\hspace{0.05\textwidth}
	\begin{subfigure}[t]{0.45\textwidth}
		\includegraphics[width=\textwidth]{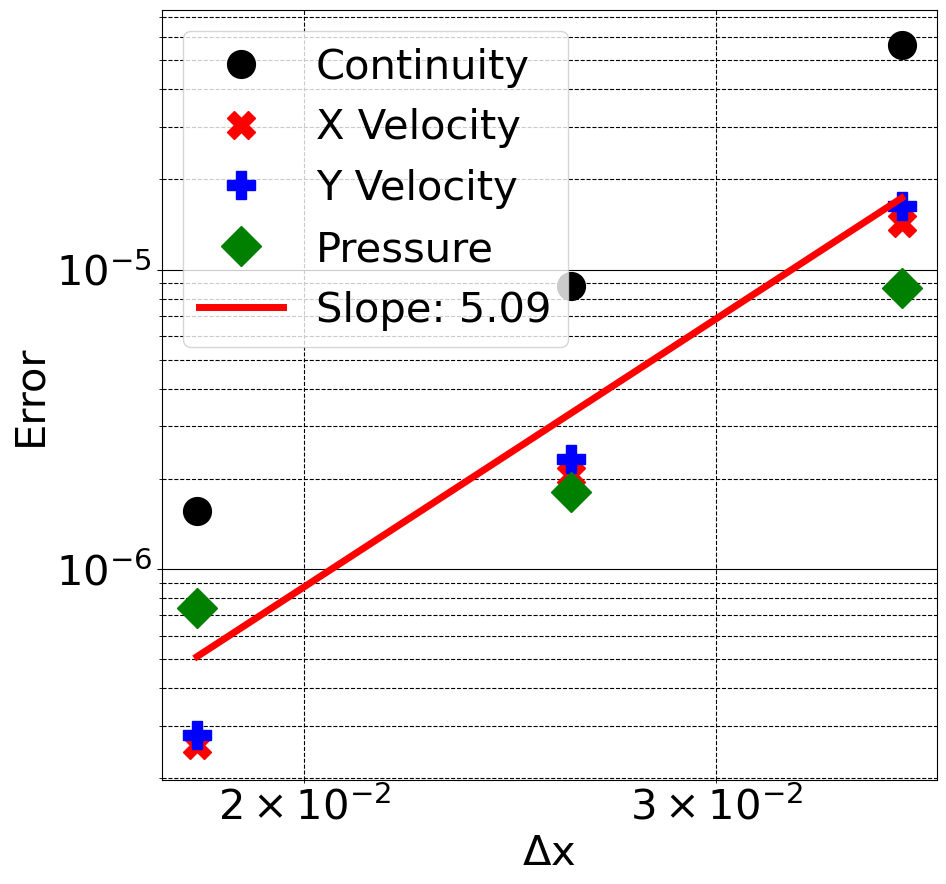}
		\caption{Polynomial Degree 6}
	\end{subfigure}
	\caption{Errors for Semi-Implicit Algorithm: Courant No.: 12}
	\label{Fig:couette error courant 12}
\end{figure}

\begin{figure}[H]
	\centering
	\includegraphics[width=0.7\textwidth]{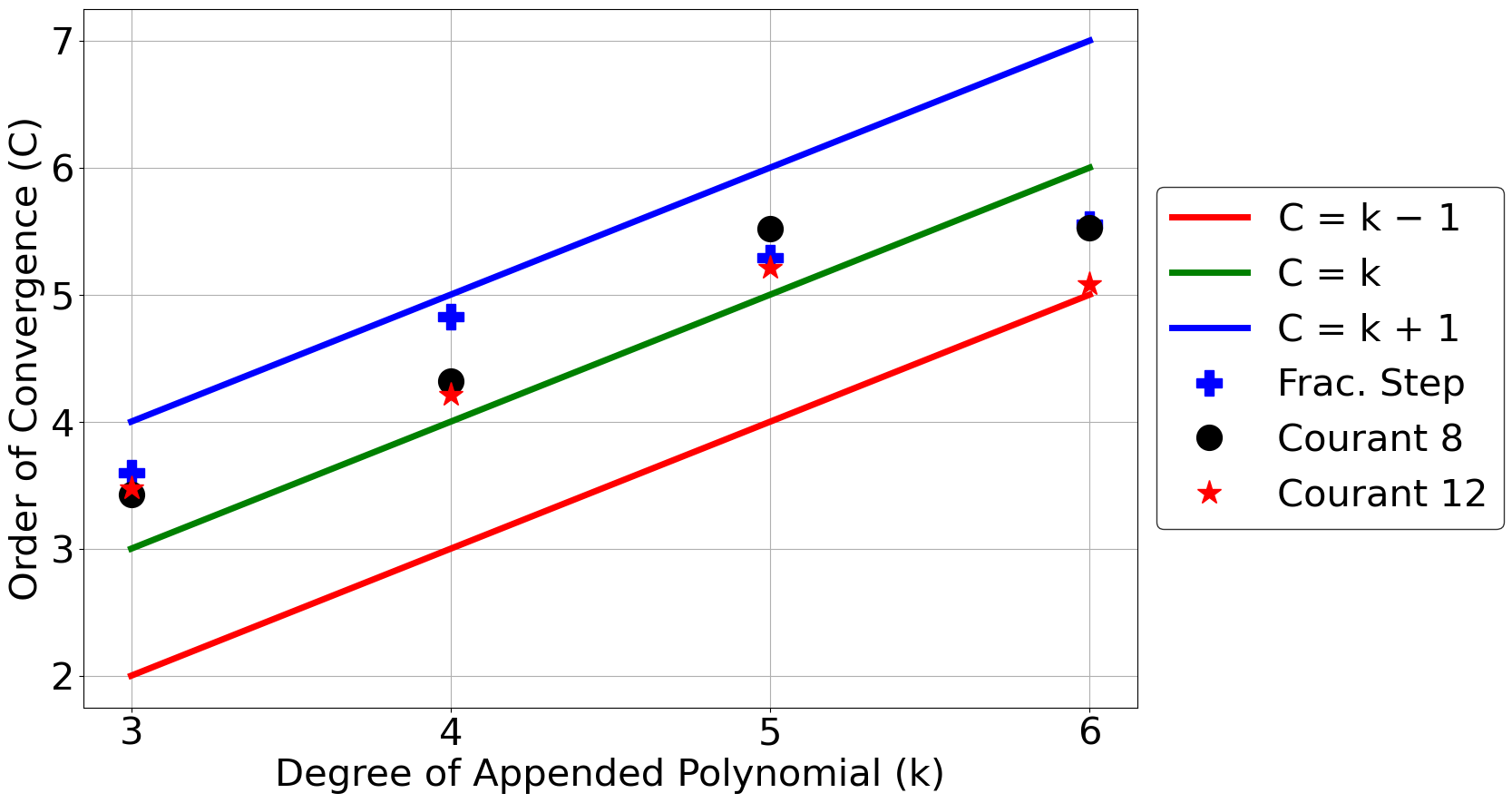}
	\caption{Order of Convergence with Polynomial Degree for Fractional Step \cite{shahane2020high} and Semi-Implicit Algorithms}
	\label{Fig:couette Order of Convergence}
\end{figure}

 As shown previously by \citet{shahane2020high} for the fractional step algorithm, the PHS basis functions and appended polynomials give high order convergence of the discretization errors. We expect the same behavior also for the semi-implicit algorithm. This is demonstrated below by calculating the $L_1$ norm of the differences in the velocity components and pressure field from the analytical solution. The divergence in steady state velocity field is defined as the error in the continuity equation. \Cref{Fig:couette error courant 12 poldeg 3} plots these errors as a function of $\Delta x$ on logarithmic scales for the three point distributions and a polynomial degree of 3. The best fit line for these 12 points is plotted. The slope of this line given by 3.47 signifies the order of convergence. Similarly, the orders of convergence are estimated for polynomial degrees of 4, 5 and 6 for the Courant number of 12 and presented in \cref{Fig:couette error courant 12}. As expected, we observe that the slopes increase monotonically with the polynomial degree.
For other Courant numbers, the slopes of the best fit lines are estimated similarly (not shown here). \Cref{Fig:couette Order of Convergence} plots these orders of convergence for the fractional step method with Courant number of 0.8 and the semi-implicit method for the Courant numbers of 8 and 12. The expected orders of convergence for interpolation, gradient and Laplacian operators are $k+1$, $k$ and $k-1$ respectively for a polynomial degree of $k$ \cite{flyer2016onrole_I,shahane2020high}. It can be seen that the order of convergence for the Navier-Stokes method, which has both gradient and Laplacian operators is between the three lines. As expected, the semi-implicit algorithm also has the same discretization error convergence properties of the fractional step method \cite{shahane2020high}.

\subsection{Problem with a Temporal Decay of the Velocity Field}
\begin{figure}[H]
	\centering
	\includegraphics[width=0.7\textwidth]{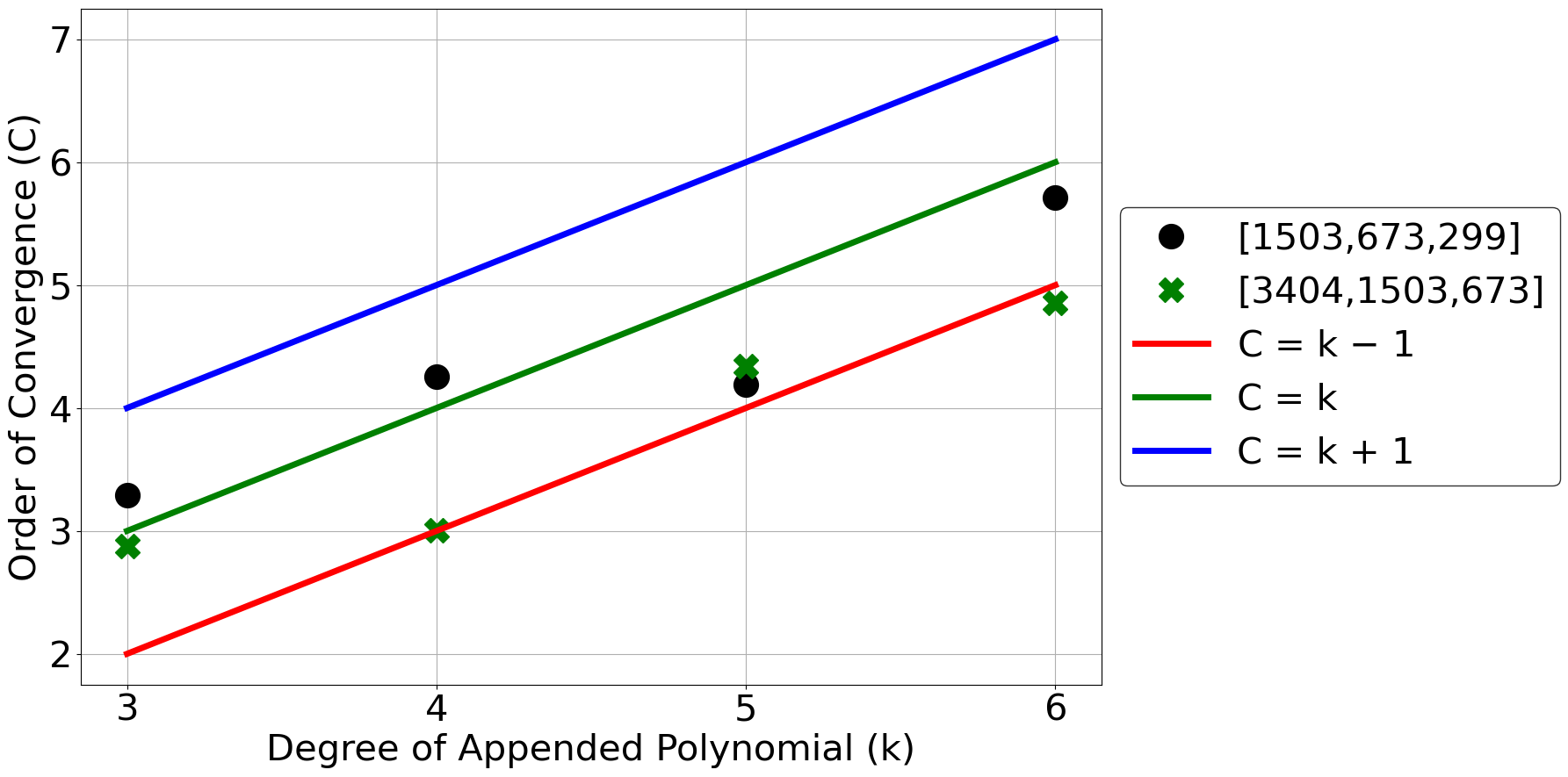}
	\caption{Order of Convergence with Polynomial Degree for Semi-Implicit Algorithm}
	\label{Fig:Order of Convergence Transient Bell-Coulela}
\end{figure}

In this section, we analyze the discretization errors for a flow which decays with time to a null solution. The initial condition on a unit square domain is given by the following stream function \cite{bell1989second}:
\begin{equation}
\Psi(x,y) = \frac{sin^2(\pi x) sin^2(\pi y)}{\pi}
\end{equation}
The initial velocity field corresponding to this stream function is:
\begin{equation}
\begin{aligned}
u(x,y)=&\frac{\partial \Psi}{\partial y} = sin^2(\pi x) sin(2 \pi y)\\
v(x,y)=&-\frac{\partial \Psi}{\partial x} = -sin(2 \pi x) sin^2(\pi y)\\
\end{aligned}
\end{equation}
The initial velocity field is divergence free. We set homogeneous Dirichlet boundary conditions on all four walls of the unit square. The Reynolds number is defined as $Re=\rho u_m L / \mu$ where, density ($\rho=1$), maximum velocity ($u_m=1$) and length of the square cavity ($L=1$) are prescribed. Based on the prescribed Reynolds number of 100, the dynamic viscosity ($\mu$) is computed. The initial velocity field is integrated in time till 0.5 seconds. Four different point distributions are used: [299, 673, 1503, 3404] points which correspond to an average $\Delta x$ of [0.0614, 0.0409, 0.0274, 0.0183] respectively. Simulations with all these points distributions and polynomial degrees from 3 to 6 are performed with a fixed $\Delta t=0.01$ seconds. After 0.5 seconds, the velocity components are interpolated to 100 uniform points along the horizontal and vertical center lines using the same PHS-RBF used for the solution. Let [$f^h$, $f^{2h}$, $f^{4h}$] denote these interpolated velocities for three successive point distributions: [1503, 673, 299] and [3404, 1503, 673]. We assume that subtracting $f^h$ from $f^{2h}$ and $f^{2h}$ from $f^{4h}$ cancels the temporal error since a fixed $\Delta t$ is used. Thus, using Richardson extrapolation, the spatial order of convergence is calculated as:
\begin{equation}
C = \log_2 \left(\frac{f^{4h} - f^{2h}}{f^{2h} - f^h}\right)
\label{Eq:richardson order expression}
\end{equation}
\Cref{Eq:richardson order expression} is applied to the $L_1$ norm of the differences between the interpolated velocities to estimate the order of convergence. These values are plotted in \cref{Fig:Order of Convergence Transient Bell-Coulela} for both the triplets of point distributions. We again see that all the orders of convergence lie in the expected range of $C=k-1$ to $C=k+1$. As a structured Cartesian point layout is not used, the $\Delta x$ is non-uniform. The slight deviations from the expected values may be attributed to the irregular point distributions. Moreover, the temporal errors may not exactly cancel when the differences are computed as in \cref{Eq:richardson order expression}.

\section{Applications to Complex Flows}
\subsection{Wake of a Circular Cylinder in a Uniform Flow}
\begin{figure}[H]
	\centering
	\includegraphics[width=0.75\textwidth]{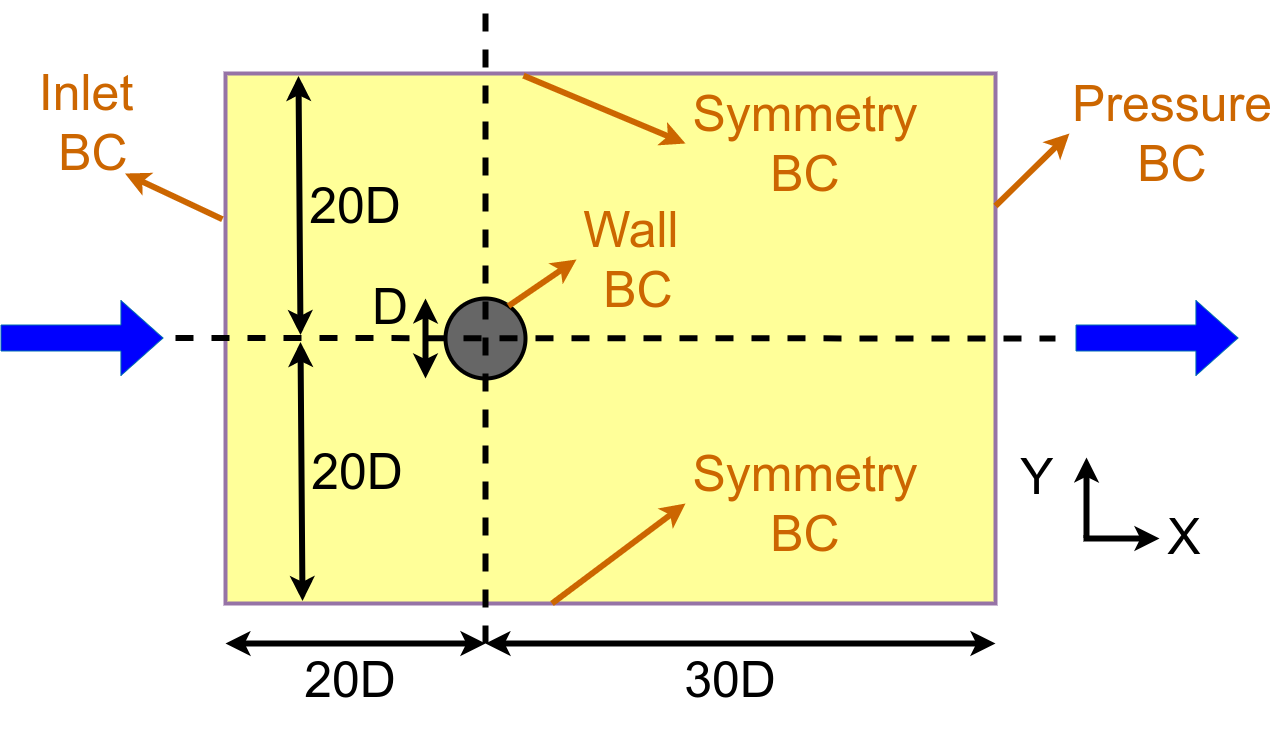}
	\caption{Domain with Diameter D $=1$ (not to scale)}
	\label{Fig:Cylinder Domain}
\end{figure}

Fluid flows over bluff bodies have been studied over several decades since they have many engineering applications. Here, we apply the semi-implicit algorithm to simulate flow over a circular cylinder. A schematic of the flow domain is shown in \cref{Fig:Cylinder Domain}. A uniform inlet flow of unit velocity is set at the left boundary. Homogeneous pressure boundary condition is applied on the right boundary with top and bottom boundaries set to symmetry conditions. No-penetration and no-slip boundary conditions are prescribed at the cylinder surface. Details of our implementation of these boundary conditions are described in \cref{Sec:Boundary Conditions}. The Reynolds number is defined as $Re=\rho U_i D /\mu$ where, density $\rho$, inlet velocity $U_i$ and diameter $D$ are all set to unity. The dynamic viscosity is estimated ($\mu$) based on the prescribed value of the Reynolds number. The inlet is placed at an upstream distance of 20D from the cylinder center. The downstream distance is taken as 30D.

\begin{figure}[H]
	\centering
	\begin{subfigure}[t]{0.45\textwidth}
		\includegraphics[width=\textwidth]{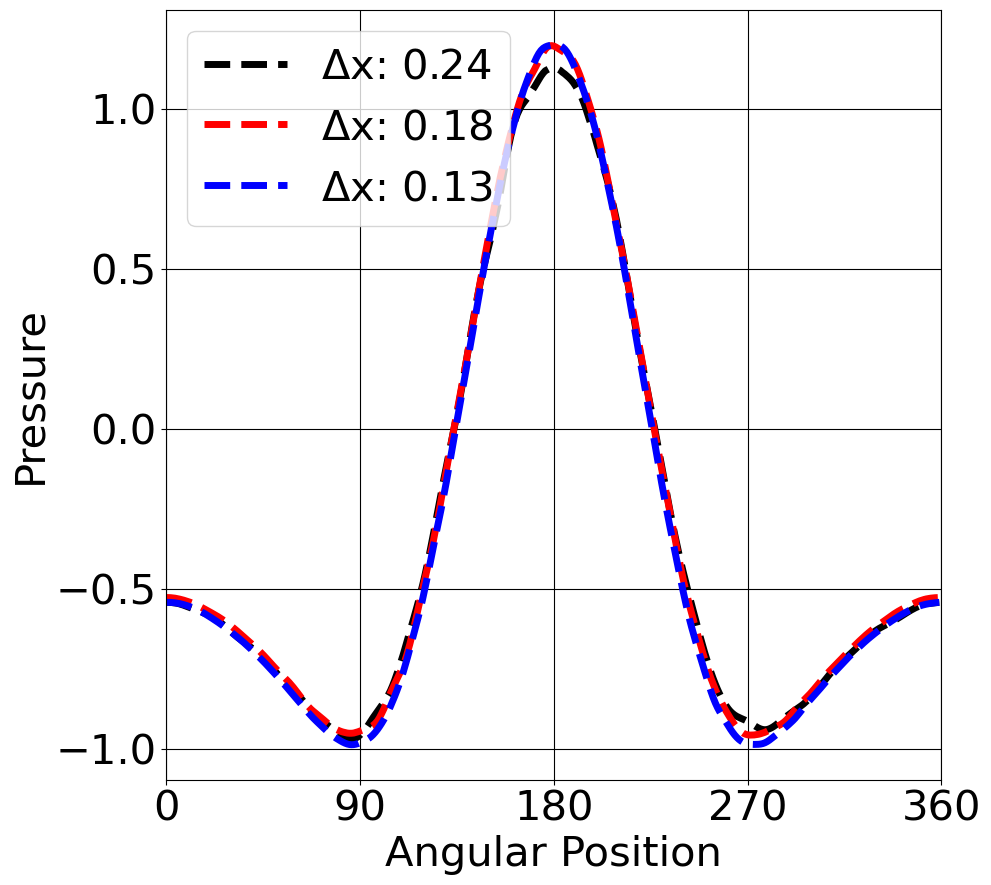}
		\caption{Polynomial Degree 5}
	\end{subfigure}
	\hspace{0.05\textwidth}
	\begin{subfigure}[t]{0.45\textwidth}
		\includegraphics[width=\textwidth]{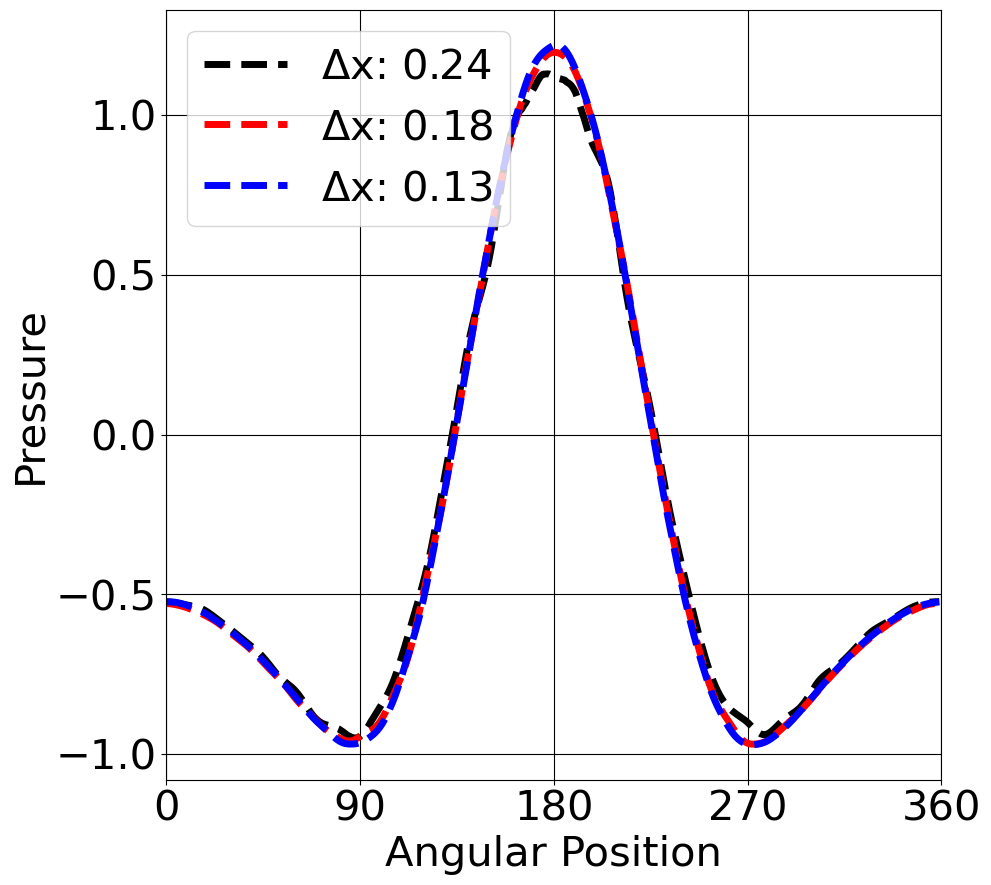}
		\caption{Polynomial Degree 6}
	\end{subfigure}
	\caption{Grid Independence for Re $=30$: Pressure Variation over Cylinder Surface}
	\label{Fig:Cylinder Grid Independence}
\end{figure}

We first consider Reynolds numbers of 10, 20 and 30 which have a steady state solution \cite{takami1969steady, tuann1978numerical, ding2004simulation, fornberg1980numerical, nieuwstadt1973viscous, gushchin1974numerical, dennis1970numerical}. Simulations are started with a uniform flow field in X direction. Time integration is performed using the BDF2 method as described in \cref{Sec:Semi-Implicit Algorithm} till the $L_1$ norm of the steady state error given by $\frac{\phi^{n+1} - \phi^n}{\Delta t}$ is less than 1E--6. $\phi^{n+1}$ and $\phi^n$ denote the velocity fields at consecutive timesteps. Three different point distributions with [33704, 59573, 105534] points is considered which corresponds to average $\Delta x$ of [0.24, 0.18, 0.13] respectively. Since the gradients are stronger near the cylinder, a non-uniform point distribution is generated using Gmsh \cite{geuzaine2009gmsh} such that the points on and around the cylinder are 5 times more refined than those in the far upstream and downstream regions. This approach is computationally more efficient than a uniform point distribution. All the simulations are performed with polynomial degrees of 5 and 6 and a Courant number of 10. \Cref{Fig:Cylinder Grid Independence} plots the variation of pressure over the cylindrical surface. It can be seen that grid independent solutions are obtained from the two finer point sets with $\Delta x$ of 0.18 and 0.13. For the remaining steady state simulations, we have used $\Delta x$ of 0.13 and the polynomial degree of 6.

Total forces per unit length experienced by the cylinder in the X and Y directions are estimated by integrating the stress components over the surface:
\begin{equation}
\begin{aligned}
F_x &= \frac{D}{2}\int_{0}^{2\pi} \left(\sigma_{xx} \cos\theta + \sigma_{yx} \sin\theta \right) d \theta \\
F_y &= \frac{D}{2}\int_{0}^{2\pi} \left(\sigma_{yy} \sin\theta + \sigma_{xy} \cos\theta \right) d \theta
\label{Eq:cylinder forces}
\end{aligned}
\end{equation}
where, components of the stress tensor ($\bm{\sigma}$) are given by:
\begin{equation}
\sigma_{xx} = -p + 2\mu\frac{\partial u}{\partial x} \hspace{0.5cm}
\sigma_{yy} = -p + 2\mu\frac{\partial v}{\partial y} \hspace{0.5cm}
\sigma_{xy}=\sigma_{yx} = \mu\left(\frac{\partial u}{\partial y} + \frac{\partial v}{\partial x}\right)
\label{Eq:cylinder stress}
\end{equation}
The pressure and strain rates are first interpolated to 360 uniform points along the cylindrical surface. Integrals in \cref{Eq:cylinder forces} are estimated using the Simpson's rule \cite{suli2003introduction} over these 360 points. The forces are non-dimensionalized to compute the drag ($C_D$) and lift ($C_L$) coefficients:
\begin{equation}
C_D = \frac{F_x}{\rho U_i^2 D/2} \hspace{1cm}
C_L = \frac{F_y}{\rho U_i^2 D/2}
\label{Eq:cylinder lift drag}
\end{equation}
\Cref{Tab:Cylinder Steady Comparison} lists the drag coefficient and pressure values at the leading and lagging stagnation regions for three Reynolds numbers in the steady state regime. Comparison with seven previously published studies shows that the current estimates are within the expected ranges. It must be emphasized that despite extensive numerical simulations of the flow over a circular cylinder, there is no clear agreement among the various solutions, primarily because of the different computational domain sizes used, upstream distances and outflow boundary conditions.  Our results are in close agreement with results of \citet{dennis1970numerical} for Re of 10 and 20, and with results of all three listed previous works for Re of 30.
\begin{table}[H]
	\centering
	\resizebox{\textwidth}{!}{%
		\begin{tabular}{|c|c|c|c|c|}
			\hline
			\begin{tabular}[c]{@{}c@{}}Reynolds \\ Number\end{tabular} &
			Reference &
			\begin{tabular}[c]{@{}c@{}}Drag \\ Coefficient\end{tabular} &
			\begin{tabular}[c]{@{}c@{}}Pressure \\ at theta=0\end{tabular} &
			\begin{tabular}[c]{@{}c@{}}Pressure \\ at theta=pi\end{tabular} \\ \hline
			\multirow{5}{*}{10} & \citet{dennis1970numerical}       & 2.846  & --0.742  & 1.489  \\
			& \citet{takami1969steady}       & 2.7541 & --0.6702 & 1.4744 \\
			& \citet{tuann1978numerical}        & 3.177  & --0.773  & 1.744  \\
			& \citet{ding2004simulation}         & 3.07   &   --    &  --    \\
			& \textbf{Present Work} & \textbf{2.877}  & \textbf{--0.7121} & \textbf{1.516}  \\ \hline
			\multirow{6}{*}{20} & \citet{dennis1970numerical}       & 2.045  & --0.589  & 1.269  \\
			& \citet{takami1969steady}       & 2.0027 & --0.537  & 1.2612 \\
			& \citet{tuann1978numerical}        & 2.253  & --0.614  & 1.457  \\
			& \citet{fornberg1980numerical}     & 2.0001 & --0.57   & 1.28   \\
			& \citet{ding2004simulation}         & 2.18   & --      &  --    \\
			& \textbf{Present Work} & \textbf{2.072}  & \textbf{--0.574 } & \textbf{1.283}  \\ \hline
			\multirow{4}{*}{30} & \citet{takami1969steady}       & 1.7167 & --0.5304 & 1.1836 \\
			& \citet{nieuwstadt1973viscous}   & 1.7329 & --0.5556 & 1.1765 \\
			& \citet{gushchin1974numerical}     & 1.808  & --      & 1.196  \\
			& \textbf{Present Work} & \textbf{1.739}  & \textbf{--0.5285} & \textbf{1.196}  \\ \hline
		\end{tabular}%
	}
	\caption{Drag Coefficient and Stagnation Pressures for Steady State Flows}
	\label{Tab:Cylinder Steady Comparison}
\end{table}

\begin{figure}[H]
	\centering
	\begin{subfigure}[t]{0.45\textwidth}
		\includegraphics[width=\textwidth]{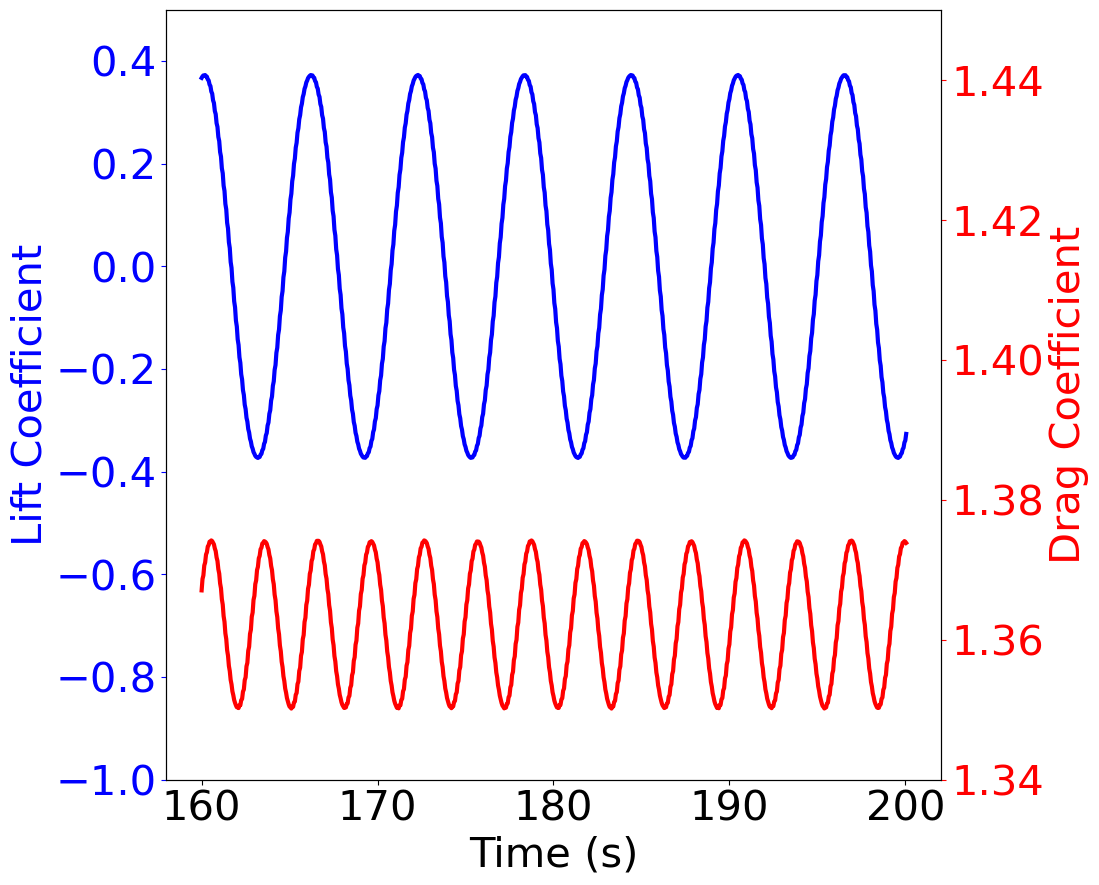}
		\caption{Reynolds Number: 100}
	\end{subfigure}
	\hspace{0.05\textwidth}
	\begin{subfigure}[t]{0.45\textwidth}
		\includegraphics[width=\textwidth]{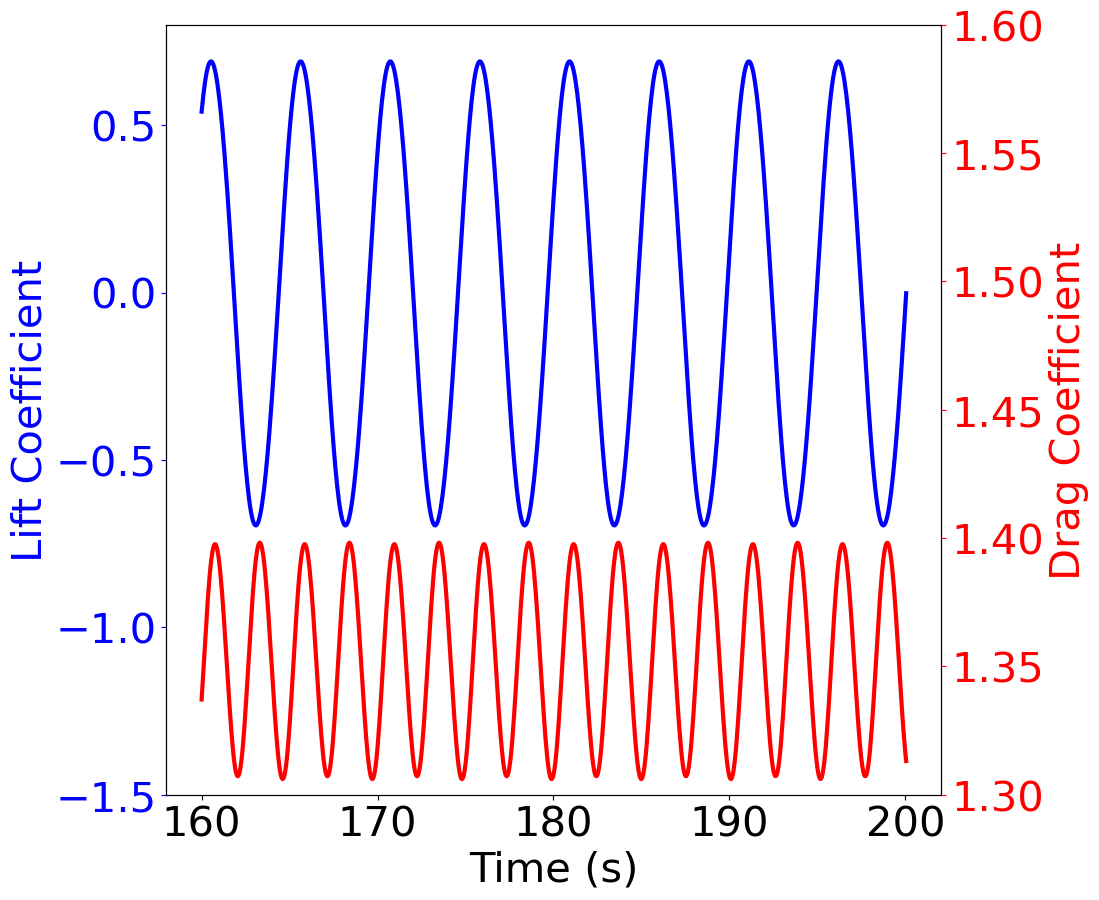}
		\caption{Reynolds Number: 200}
	\end{subfigure}
	\caption{Lift and Drag Coefficients (163798 points, Average $\Delta x = 0.1105$, Degree of Appended Polynomial: 5)}
	\label{Fig:Cylinder Lift and Drag Coefficients}
\end{figure}

Next we present the results for Reynolds numbers of 100 and 200 for which there is temporal variation due to vortex shedding. For these cases, a point distribution with 163798 points corresponding to an average $\Delta x$ of 0.10 with a polynomial degree of 6 is used. As described before, non-uniform points distribution is used for computational efficiency. Starting from a uniform velocity field, the equations are integrated in time till a stationary solution is obtained. \Cref{Fig:Cylinder Lift and Drag Coefficients} plots the temporal variations of lift and drag coefficients for both the Reynolds numbers. Sinusoidal variations are observed as expected \cite{ding2004simulation}. The Strouhal number is computed from the frequency of the lift coefficient: $St=fD/U_i$. \Cref{Tab:Cylinder Unsteady Comparison} documents the ranges of lift and drag coefficients together with the Strouhal numbers. As for the steady cases, present estimates fall within the documented ranges in the literature.

\begin{table}[H]
	\centering
	\resizebox{\textwidth}{!}{%
		\begin{tabular}{|c|c|c|c|c|}
			\hline
			\begin{tabular}[c]{@{}c@{}}Reynolds \\ Number\end{tabular} &
			Reference &
			\begin{tabular}[c]{@{}c@{}}Drag \\ Coefficient\end{tabular} &
			\begin{tabular}[c]{@{}c@{}}Lift \\ Coefficient\end{tabular} &
			\begin{tabular}[c]{@{}c@{}}Strouhal \\ Number\end{tabular} \\ \hline
			\multirow{4}{*}{100} & \citet{braza1986numerical}        & 1.364 $\pm$ 0.015  & $\pm$ 0.25  & 0.16  \\
			& \citet{liu1998preconditioned}          & 1.35 $\pm$ 0.012   & $\pm$ 0.339 & 0.165 \\
			& \citet{ding2004simulation}         & 1.325 $\pm$ 0.008  & $\pm$ 0.28  & 0.164 \\
			& \textbf{Present Work} & \textbf{1.362 $\pm$ 0.0118} & $\pm$ \textbf{0.3331}  & \textbf{0.1649} \\ \hline
			\multirow{5}{*}{200} & \citet{belov1995new}        & 1.19 $\pm$ 0.042   & $\pm$ 0.64  &       \\
			& \citet{braza1986numerical}         & 1.4 $\pm$ 0.05     & $\pm$ 0.75  & 0.193 \\
			& \citet{liu1998preconditioned}          & 1.31 $\pm$ 0.049   & $\pm$ 0.69  & 0.192 \\
			& \citet{ding2004simulation}         & 1.327 $\pm$ 0.045  & $\pm$ 0.60  & 0.196 \\
			& \textbf{Present Work} & \textbf{1.353 $\pm$ 0.0448}   & $\pm$ \textbf{0.6931}  & \textbf{0.1963} \\ \hline
		\end{tabular}%
	}
	\caption{Strouhal Number, Drag and Lift Coefficients for Unsteady Flows}
	\label{Tab:Cylinder Unsteady Comparison}
\end{table}

\subsection{Couette Flow in an Ellipse with Two Rotating Inner Cylinders}
\begin{figure}[H]
	\begin{subfigure}[t]{0.45\textwidth}
		\includegraphics[width=\textwidth]{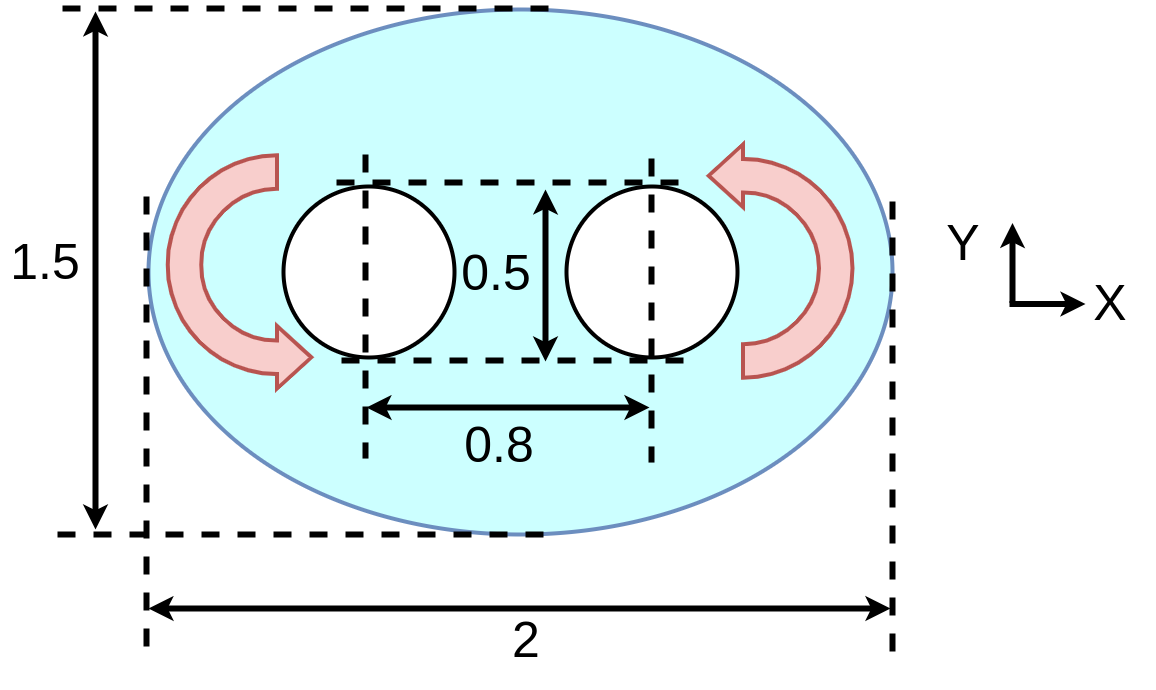}
		\caption{Unidirectional Rotation}
		\label{Fig:Double Couette Geometry unidirectional}
	\end{subfigure}
	\hspace{0.05\textwidth}
	\begin{subfigure}[t]{0.45\textwidth}
		\includegraphics[width=\textwidth]{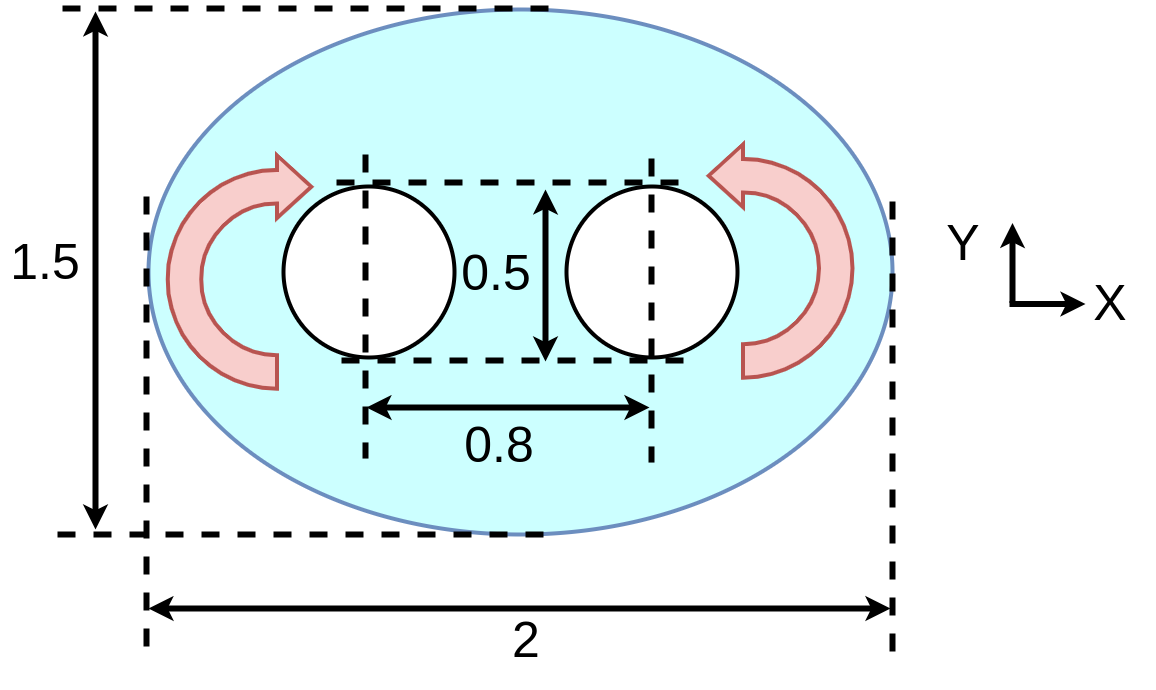}
		\caption{Counter Rotation}
		\label{Fig:Double Couette Geometry counter}
	\end{subfigure}
	\caption{Geometry}
	\label{Fig:Double Couette Geometry}
\end{figure}
We next present results of a Couette flow generated by two rotating circular cylinders inside a stationary elliptical enclosure. The fluid is initially at rest in the annular region between the cylinders. \Cref{Fig:Double Couette Geometry unidirectional} shows the domain with dimensions in which both the cylinders rotate in the counterclockwise direction and \cref{Fig:Double Couette Geometry counter} shows the case with cylinders rotating in opposite directions. A rotational velocity ($\omega$) of 1 radians per second is set at both the inner cylinders. The Reynolds number is defined with respect to the radius of the inner cylinder and its velocity: $Re=\rho \omega r_i^2 /\mu$ where, $\omega=1$, $\rho=1$, $r_i=0.25$ and $\mu$ is computed based on the Reynolds number. We have prescribed a Reynolds number of 40 for which we see a steady two dimensional flow.
\begin{figure}[H]
	\begin{subfigure}[t]{0.45\textwidth}
		\includegraphics[width=\textwidth]{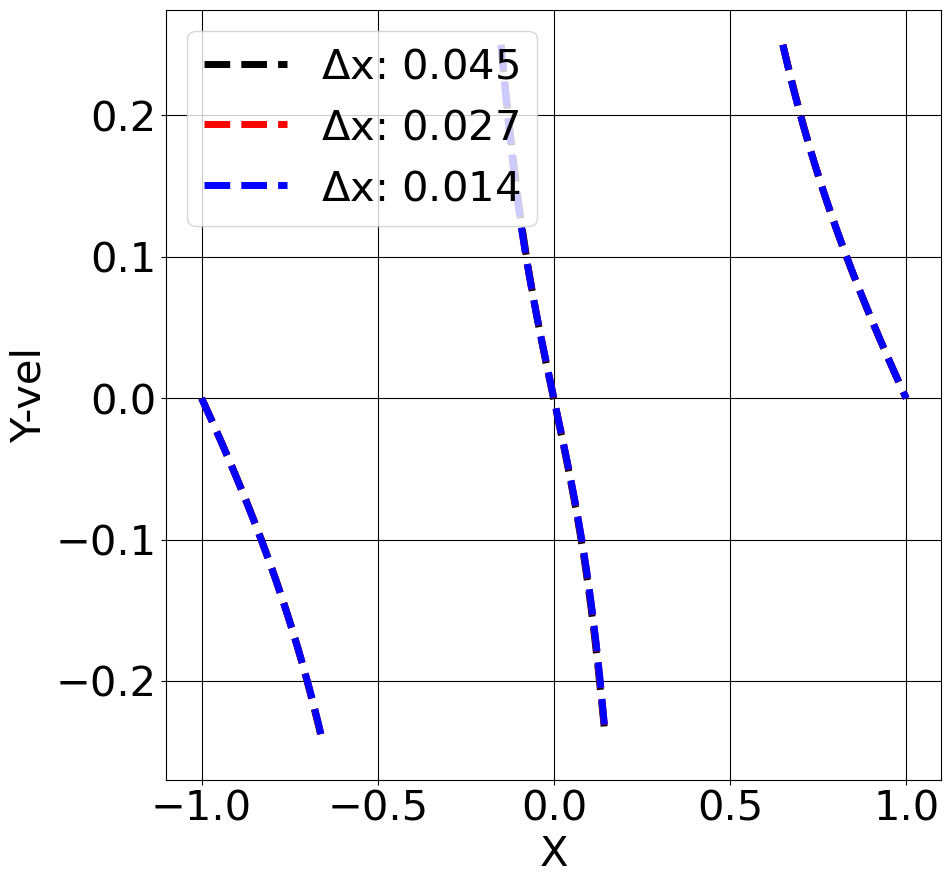}
		\caption{Unidirectional Rotation}
	\end{subfigure}
	\hspace{0.05\textwidth}
	\begin{subfigure}[t]{0.45\textwidth}
		\includegraphics[width=\textwidth]{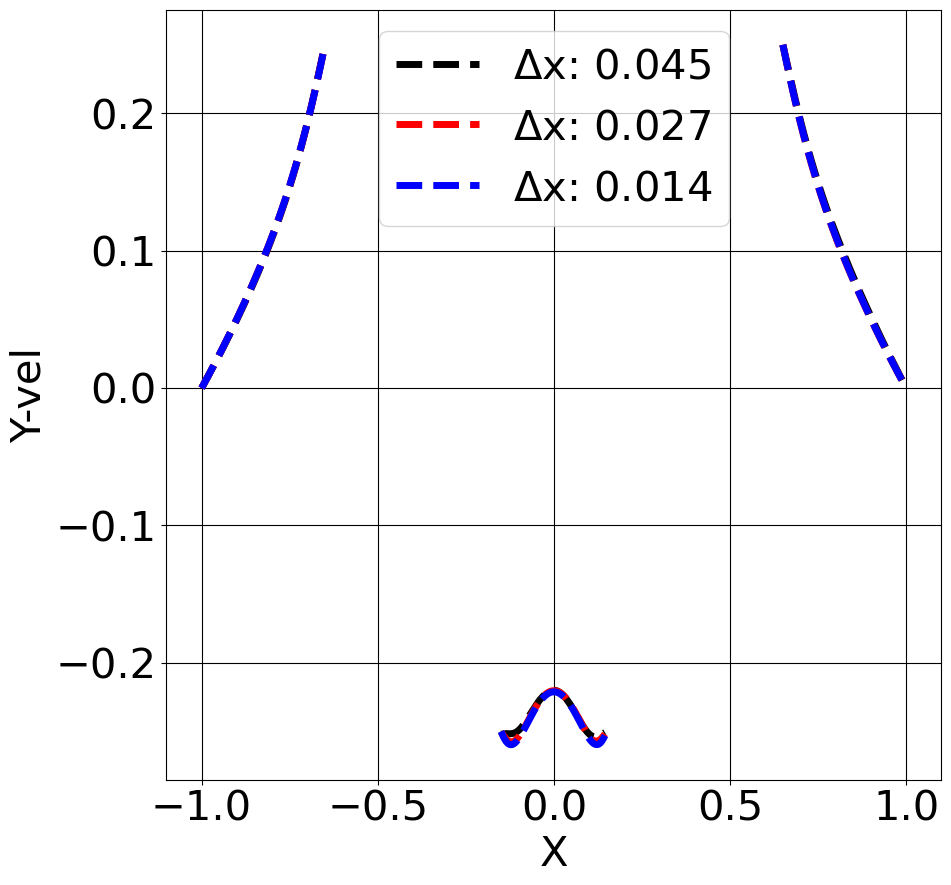}
		\caption{Counter Rotation}
	\end{subfigure}
	\caption{Grid Independence for Polynomial Degree 5: Y-Vel along Horizontal Center Line}
	\label{Fig:Double Couette Geometry Grid Independence}
\end{figure}

\begin{figure}[H]
	\begin{subfigure}[t]{0.45\textwidth}
		\includegraphics[width=\textwidth]{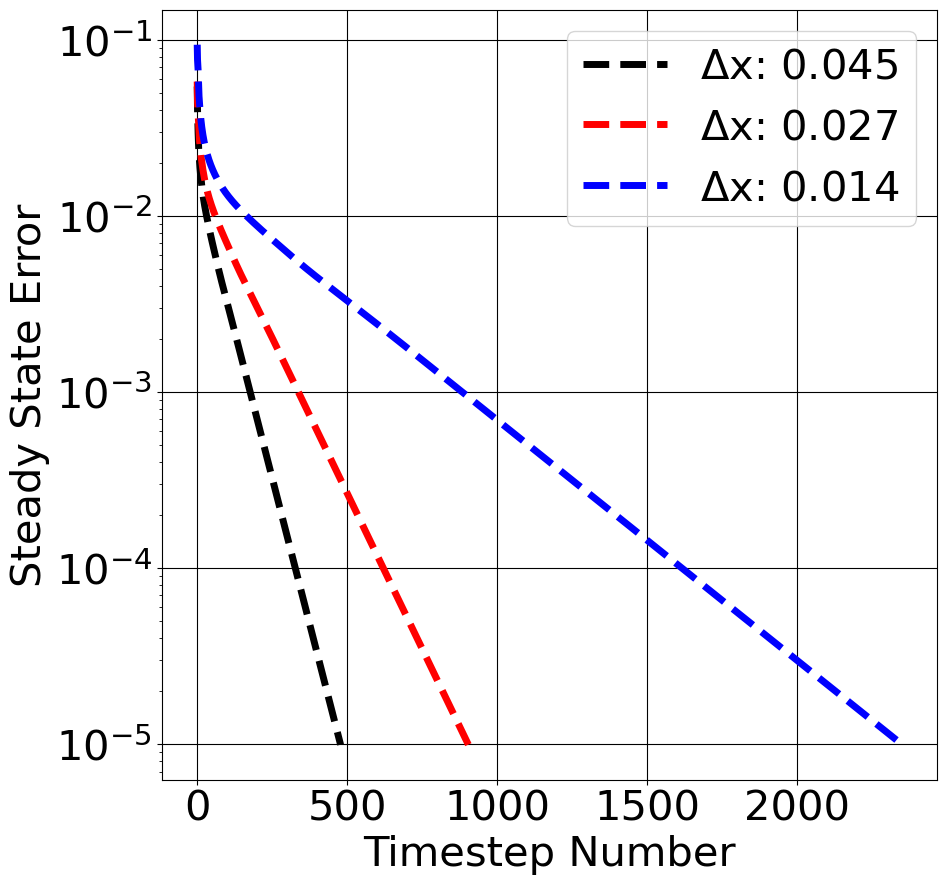}
		\caption{Unidirectional Rotation}
	\end{subfigure}
	\hspace{0.05\textwidth}
	\begin{subfigure}[t]{0.45\textwidth}
		\includegraphics[width=\textwidth]{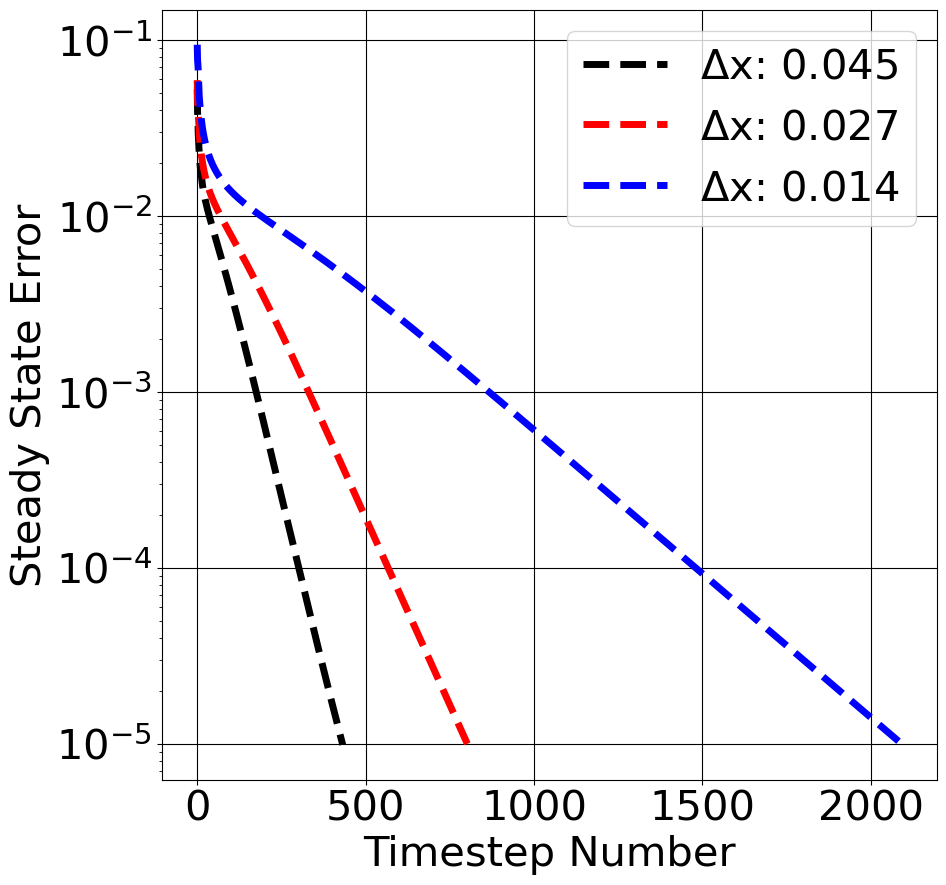}
		\caption{Counter Rotation}
	\end{subfigure}
	\caption{Convergence of the Semi-Implicit Algorithm with Polynomial Degree 5 and Courant No. 10}
	\label{Fig:Double Couette steady error}
\end{figure}

We have simulated the flows using three point distributions with 1021, 2748 and 10843 points which correspond to average $\Delta x$ values of 0.045, 0.027 and 0.014. We have used a polynomial degree of 5 for all the simulations. \Cref{Fig:Double Couette Geometry Grid Independence} shows that the Y component of the velocity along the horizontal center lines overlap. The finest point set is used for all the subsequent simulations. \Cref{Fig:Double Couette steady error} plots convergence of the steady state error with timestep for the polynomial degree of 5. Similar to the plots shown in \cref{Sec:couette flow}, the case with higher $\Delta x$ converges to steady state in less number of timesteps for the same Courant number.

\begin{figure}[H]
	\begin{subfigure}[t]{0.45\textwidth}
		\includegraphics[width=\textwidth]{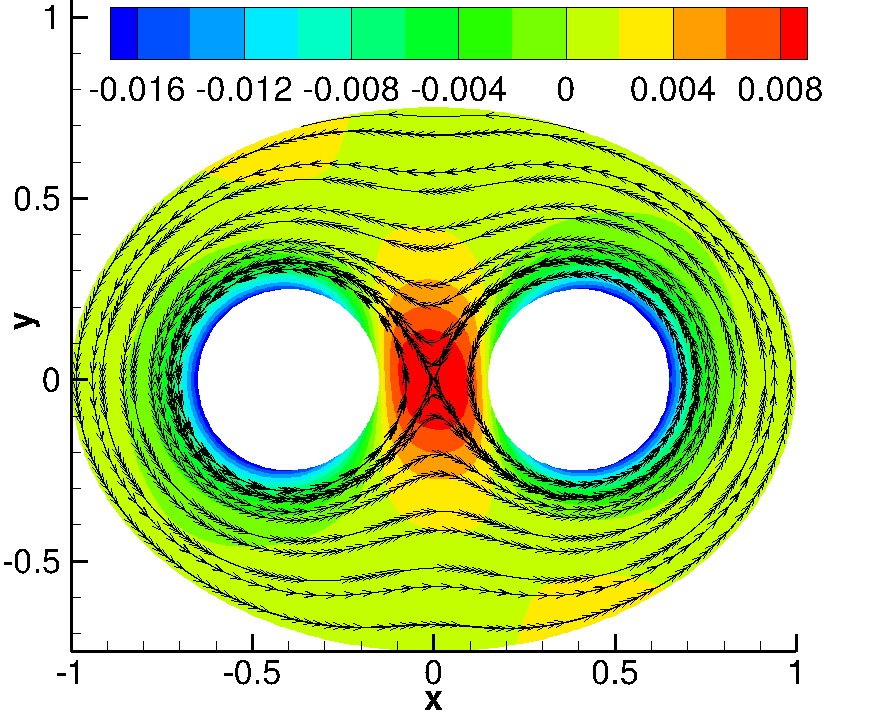}
		\caption{Unidirectional Rotation}
	\end{subfigure}
	\hspace{0.05\textwidth}
	\begin{subfigure}[t]{0.45\textwidth}
		\includegraphics[width=\textwidth]{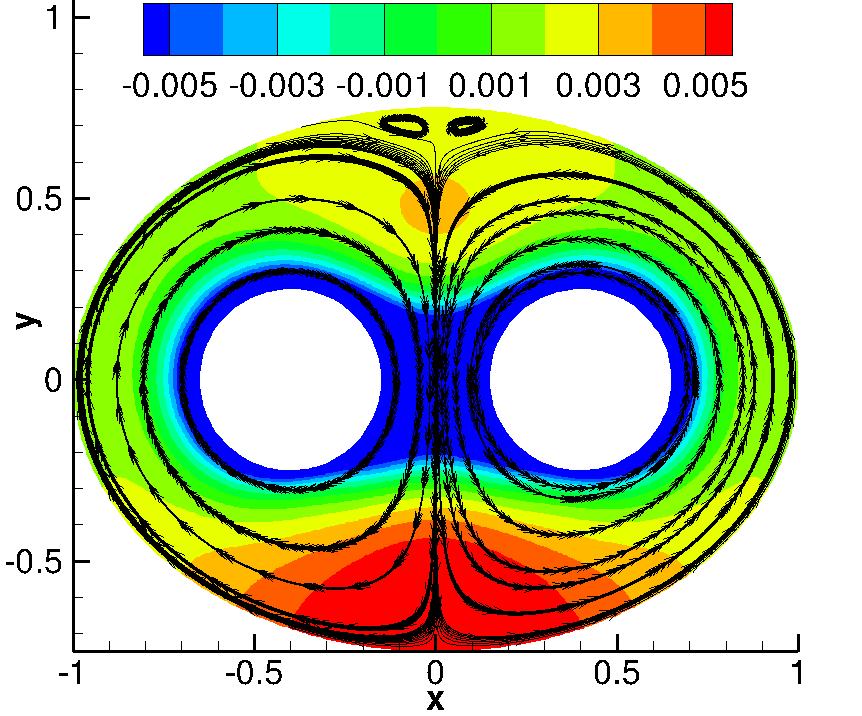}
		\caption{Counter Rotation}
	\end{subfigure}
	\caption{Contours of Pressure with Streamlines}
	\label{Fig:Double Couette Pressure with Streamlines}
\end{figure}

\Cref{Fig:Double Couette Pressure with Streamlines} plots the contours of pressure superposed with streamlines. For the case of unidirectional rotation, a stagnation point is formed between the cylinders as the velocities are in opposite directions in that region. On the other hand, in the case of counter rotation, the velocities add in the region between the two cylinders. Thus, we see a region of low pressure between the the cylinders and a stagnation high pressure region on the bottom wall.

\subsection{Flow in a Periodic Bellowed Channel}
\begin{figure}[H]
	\centering
	\includegraphics[width=0.6\textwidth]{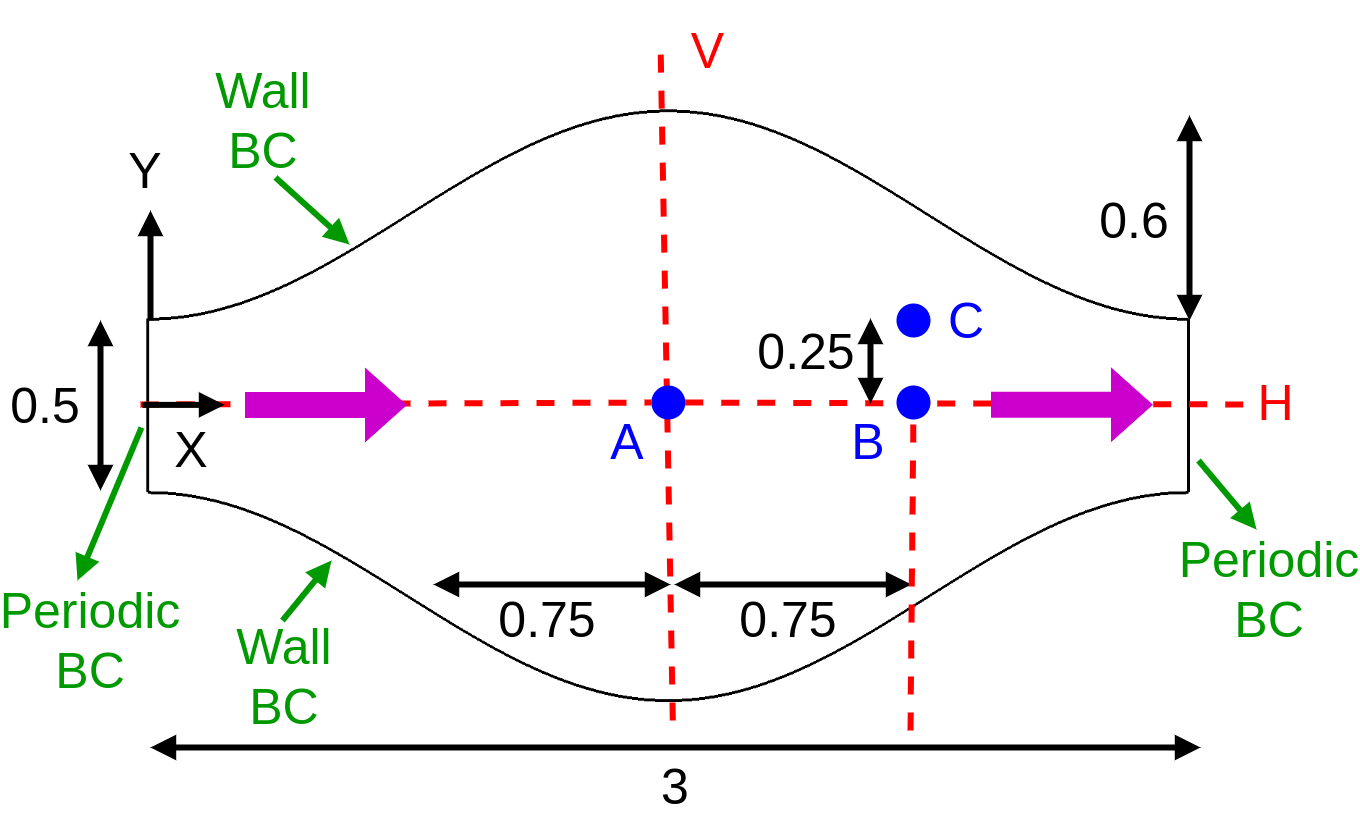}
	\caption{Geometry of the Bellowed Channel}
	\label{Fig:bellow schematic}
\end{figure}
\Cref{Fig:bellow schematic} shows the dimensions of a sinusoidal bellow with fluid flowing from left to right \cite{wang1995convective,stephanoff1980flow}. The top and bottom boundaries are walls with no-slip and no-penetration conditions. Periodicity is prescribed on the left and right boundaries as described in \cref{Sec:Boundary Conditions}. Velocity fields are interpolated on uniform points along a horizontal and vertical center lines for analysis. Similarly, temporal variations of the fields at three reference points A, B and C are also plotted. The minimum width of the channel used as the characteristic length ($L$) is set to 0.5 and density ($\rho$) is set to unity. The periodic flow is driven with a mean pressure gradient in the X direction:
\begin{equation}
\frac{\partial \bar{p}}{\partial x} = \frac{12 \mu}{L^2}
\label{Eq:bellow dp_dx}
\end{equation}
Velocity in the X direction averaged over the minimum width at the left boundary is used to calculate the flow Reynolds number: $Re=\frac{\rho \bar{u} L}{\mu}$. Simulations are performed for two different viscosity values: 0.01 and 0.002 which correspond to Reynolds numbers of 55.6 and 172.9 respectively. We find that 55.6 gives a steady state flow whereas, 172.9 generates an unsteady flow. Both the simulations are started with a zero initial condition and integrated in time till stationary fields are obtained.
\par
\begin{figure}[H]
	\centering
	\begin{subfigure}[t]{0.45\textwidth}
		\includegraphics[width=\textwidth]{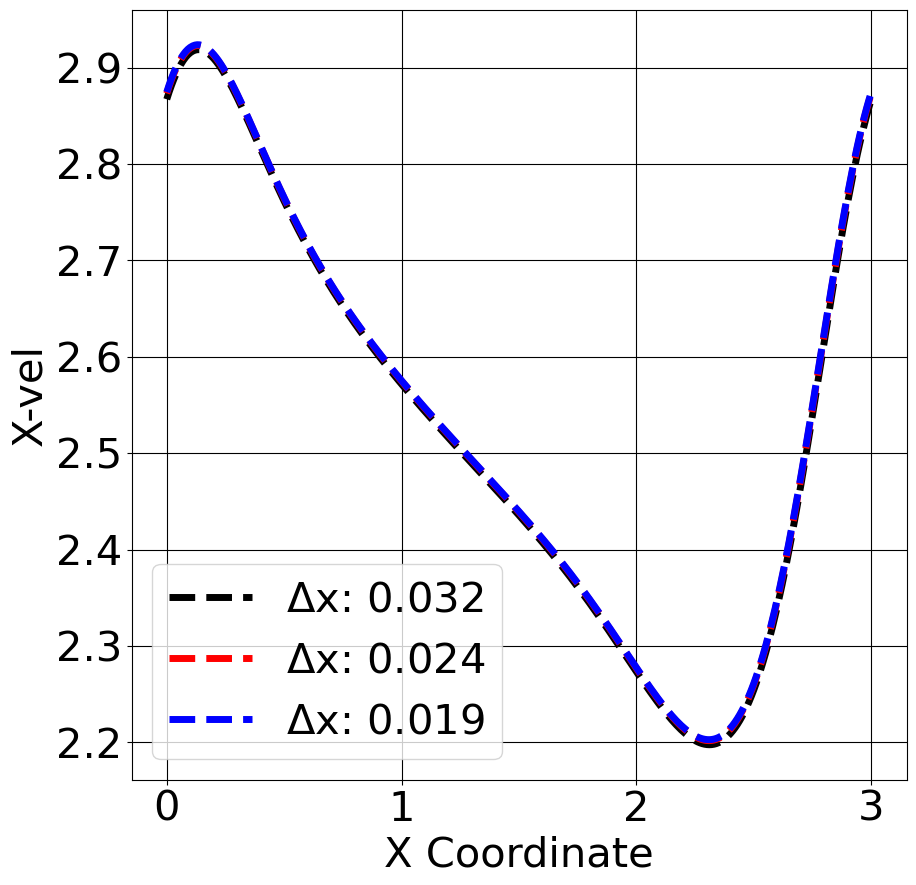}
		\caption{X-Vel: Horizontal Center Line (H)}
	\end{subfigure}
	\hspace{0.05\textwidth}
	\begin{subfigure}[t]{0.45\textwidth}
		\includegraphics[width=\textwidth]{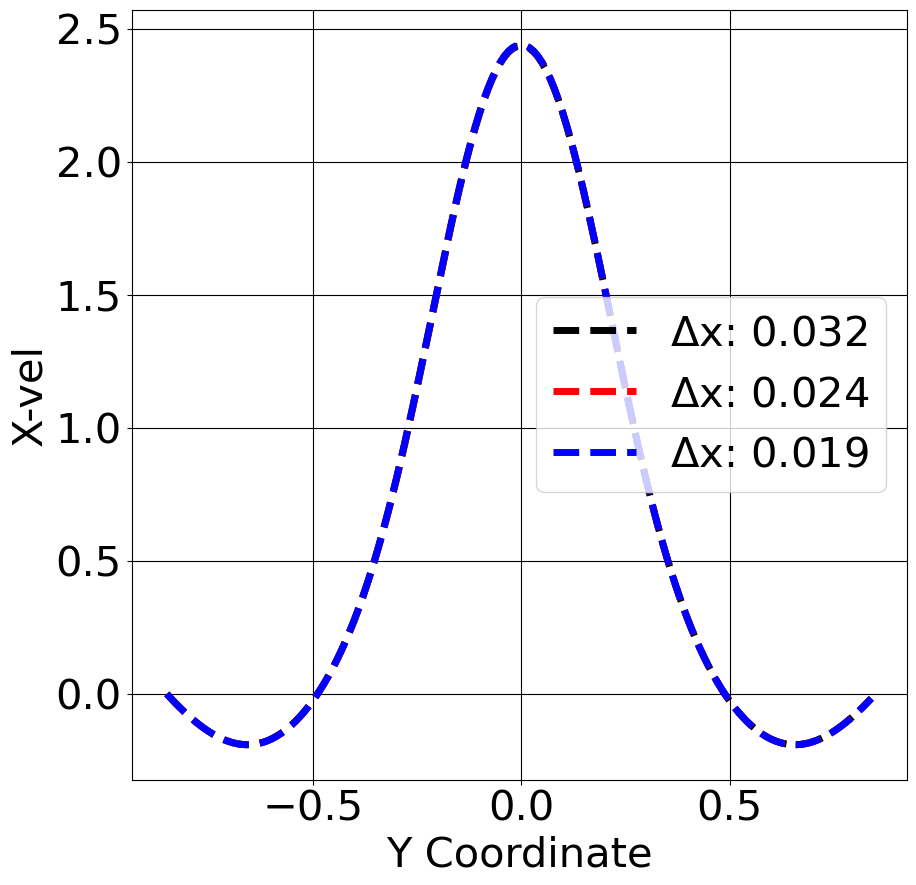}
		\caption{X-Vel: Vertical Center Line (V)}
	\end{subfigure}\vspace{0.25cm}
	\begin{subfigure}[t]{0.45\textwidth}
		\includegraphics[width=\textwidth]{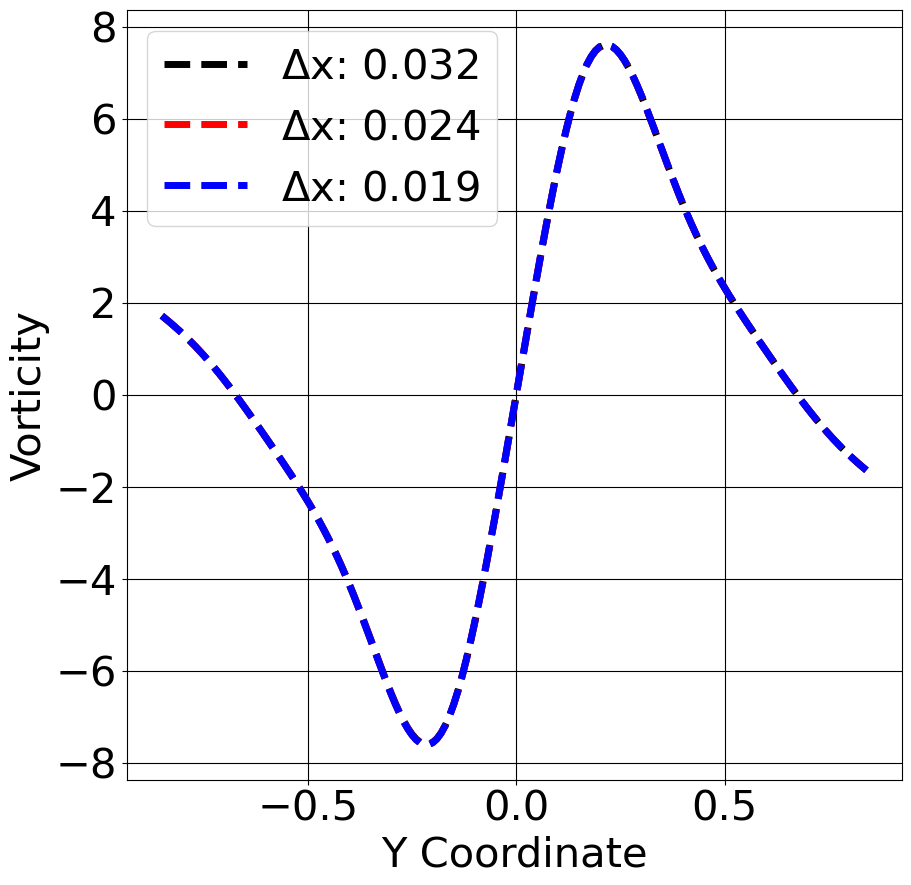}
		\caption{Vorticity: Vertical Center Line (V)}
	\end{subfigure}
	\caption{Grid Independence for Re $=55.6$ with Polynomial Degree 5}
	\label{Fig:bellow Grid Independence}
\end{figure}

\begin{figure}[H]
	\centering
	\begin{subfigure}[t]{0.45\textwidth}
		\includegraphics[width=\textwidth]{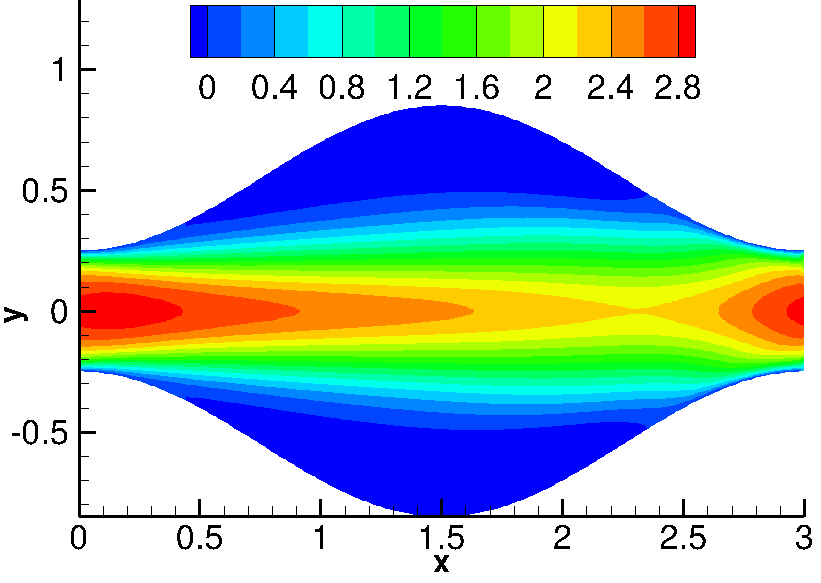}
		\caption{X-Vel}
	\end{subfigure}
	\hspace{0.05\textwidth}
	\begin{subfigure}[t]{0.45\textwidth}
		\includegraphics[width=\textwidth]{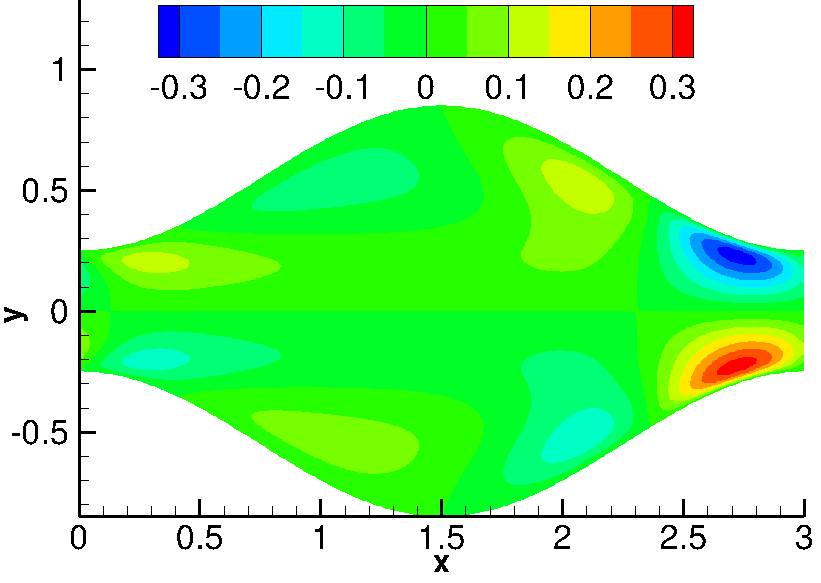}
		\caption{Y-Vel}
	\end{subfigure}\vspace{0.25cm}
	\begin{subfigure}[t]{0.45\textwidth}
		\includegraphics[width=\textwidth]{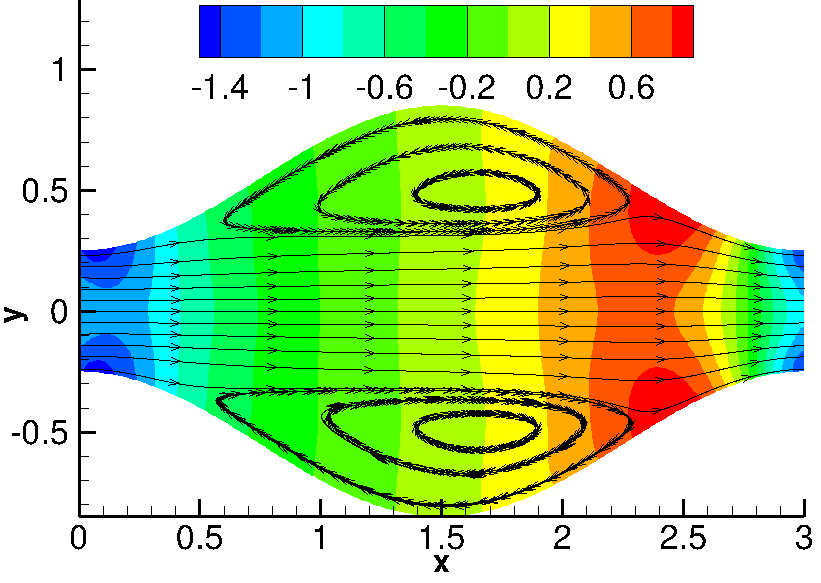}
		\caption{Pressure Variations and Streamlines}
		\label{Fig:bellow Re 55.6 contours p streamline}
	\end{subfigure}
	\caption{Contour Plots for Re $=55.6$ ($\Delta x=0.019$ with Polynomial Degree 5)}
	\label{Fig:bellow Re 55.6 contours}
\end{figure}

 Three different point distributions with 3626, 6367 and 9778 points which correspond to average $\Delta x$ of 0.032, 0.024 and 0.019 respectively are used with a polynomial degree of 5. The X component of velocity and vorticity are plotted along the horizontal and vertical center lines in \cref{Fig:bellow Grid Independence}. It can be seen that the line plots for all the point distributions overlap. For all the simulations in this section, we use the finest distribution with 9778 points and average $\Delta x=0.019$.
\Cref{Fig:bellow Re 55.6 contours} plots the contours for the Reynolds number of 55.6. A steady state solution is obtained for this case. Note that since a mean pressure gradient (\cref{Eq:bellow dp_dx}) is applied as a force in the X-momentum equation, pressure perturbations on top of the mean pressure are plotted in \cref{Fig:bellow Re 55.6 contours p streamline} together with streamlines. As expected, the steady flow is symmetric with respect to the horizontal center line.

\begin{figure}[H]
	\centering
	\begin{subfigure}[t]{0.45\textwidth}
		\includegraphics[width=\textwidth]{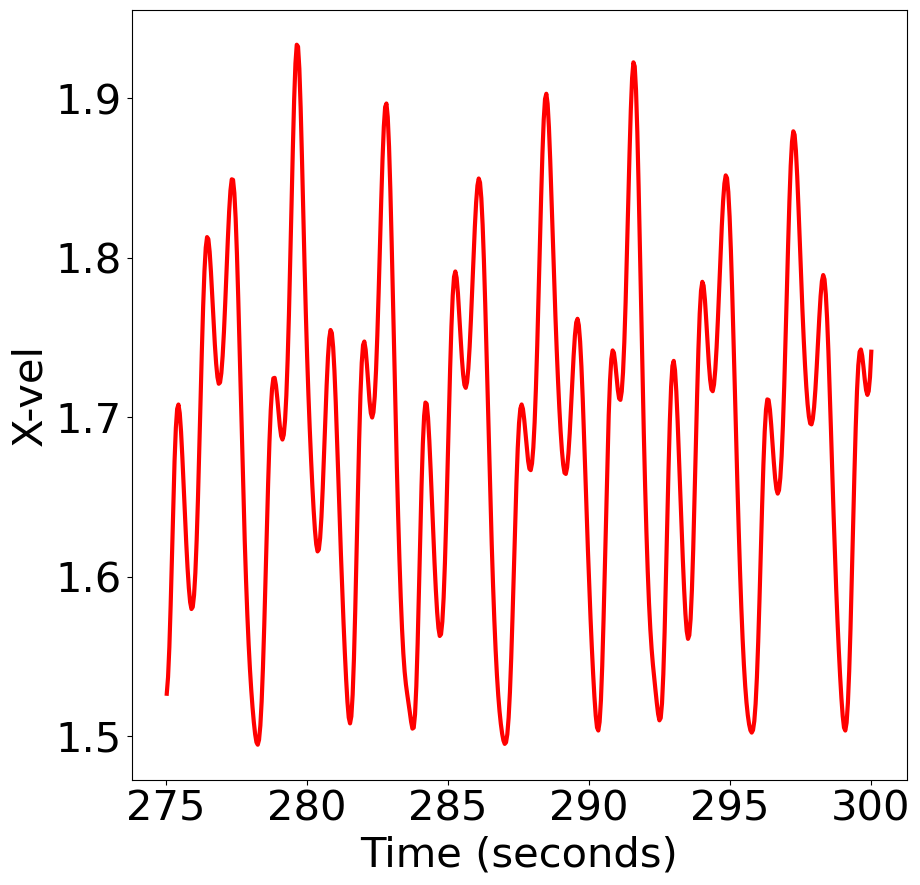}
		\caption{Point A: X$=1.5$, Y$=0$}
	\end{subfigure}
	\hspace{0.05\textwidth}
	\begin{subfigure}[t]{0.45\textwidth}
		\includegraphics[width=\textwidth]{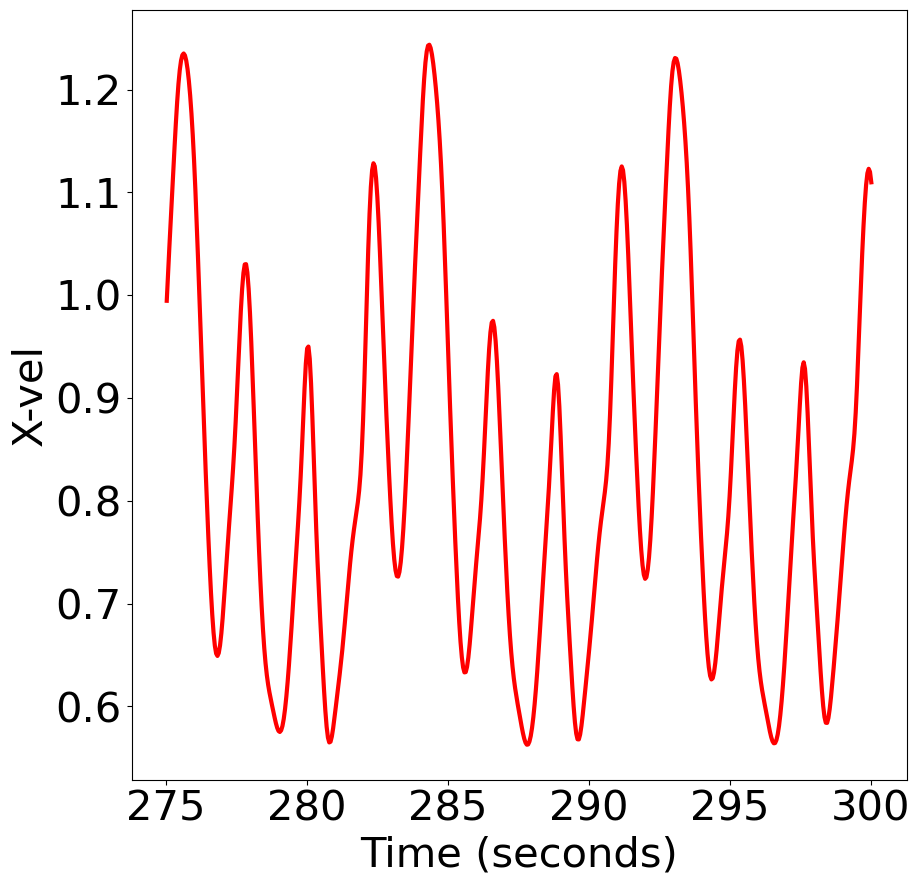}
		\caption{Point C: X$=2.25$, Y$=0.25$}
	\end{subfigure}
	\caption{Temporal Variation of X-Vel for Re$=172.9$}
	\label{Fig:bellow X-vel temporal}
\end{figure}

\begin{figure}[H]
	\centering
	\begin{subfigure}[t]{0.45\textwidth}
		\includegraphics[width=\textwidth]{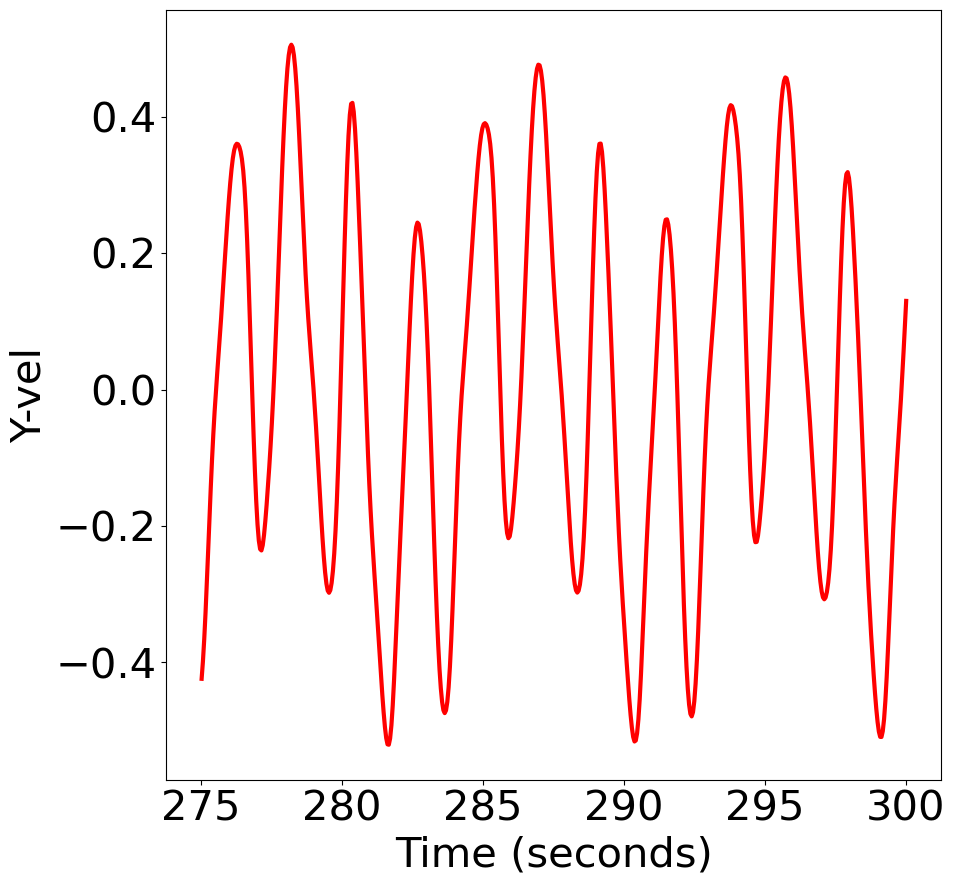}
		\caption{Point B: X$=2.25$, Y$=0$}
	\end{subfigure}
	\hspace{0.05\textwidth}
	\begin{subfigure}[t]{0.45\textwidth}
		\includegraphics[width=\textwidth]{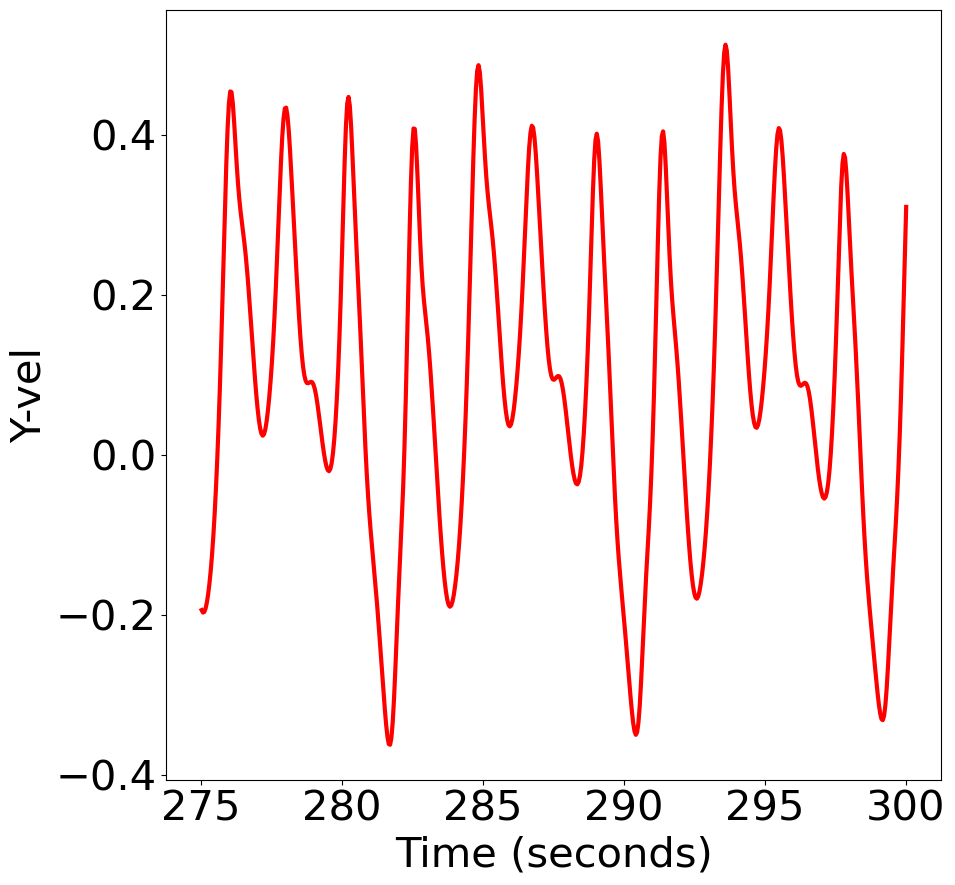}
		\caption{Point C: X$=2.25$, Y$=0.25$}
	\end{subfigure}
	\caption{Temporal Variation of Y-Vel for Re $=172.9$}
	\label{Fig:bellow Y-vel temporal}
\end{figure}

\begin{figure}[H]
	\centering
	\begin{subfigure}[t]{0.45\textwidth}
		\includegraphics[width=\textwidth]{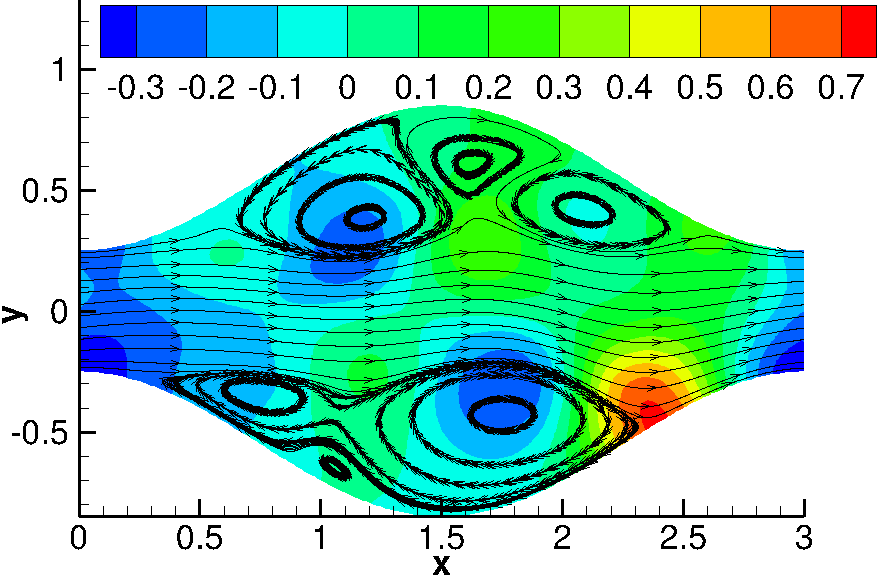}
		\caption{Time: 220 seconds}
	\end{subfigure}
	\hspace{0.05\textwidth}
	\begin{subfigure}[t]{0.45\textwidth}
		\includegraphics[width=\textwidth]{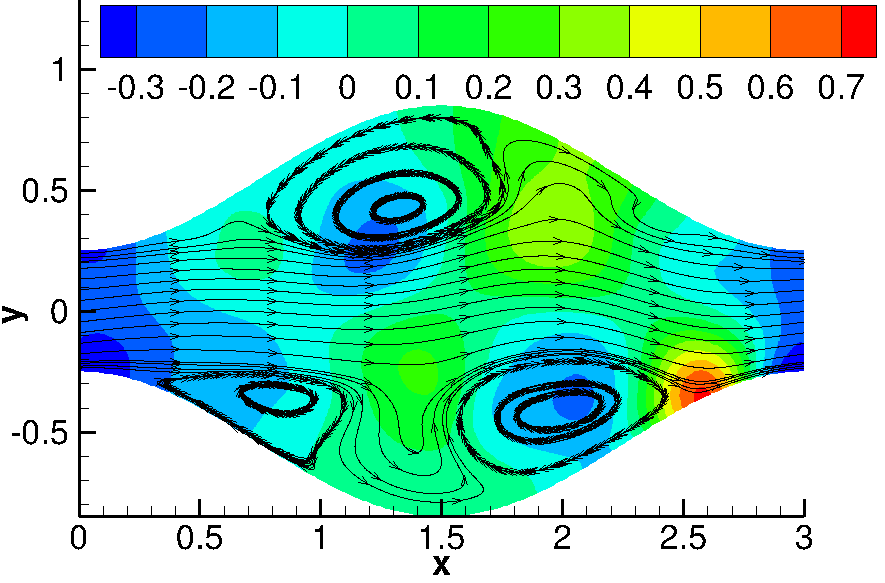}
		\caption{Time: 240 seconds}
	\end{subfigure}
	\begin{subfigure}[t]{0.45\textwidth}
		\includegraphics[width=\textwidth]{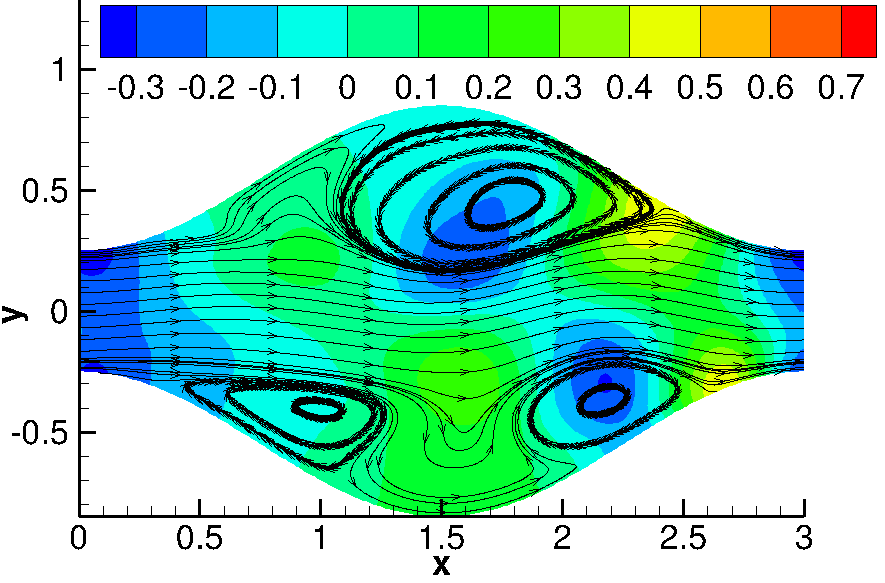}
		\caption{Time: 260 seconds}
	\end{subfigure}
	\hspace{0.05\textwidth}
	\begin{subfigure}[t]{0.45\textwidth}
		\includegraphics[width=\textwidth]{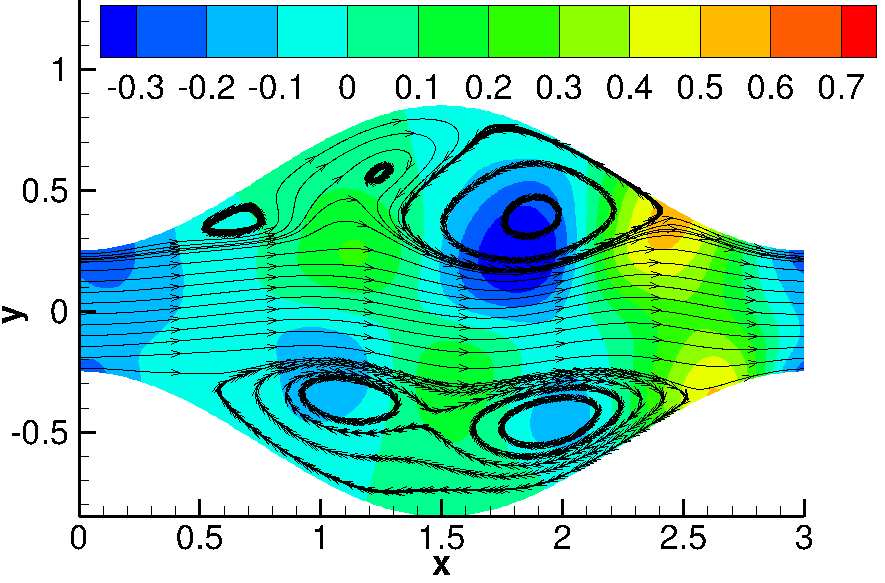}
		\caption{Time: 280 seconds}
	\end{subfigure}
	\caption{Contours of Excess Pressure over Mean Gradient with Streamlines for Re$=172.9$ at Various Time Instants}
	\label{Fig:bellow pressure streamlines Re 172.9}
\end{figure}

\Cref{Fig:bellow X-vel temporal,Fig:bellow Y-vel temporal} plot the temporal variation of the velocity components at the three points A, B and C marked in \cref{Fig:bellow schematic} for the Reynolds number of 172.9. We observe that the X-velocity oscillates around a mean positive value since the flow direction is from left to right. The mean value at point C is lower than the point A which is along the horizontal center line. The Y-velocity on the other hand oscillates around a mean of zero at both points B and C. Multiple frequencies are observed in the velocity signals at the points A, B and C. \Cref{Fig:bellow pressure streamlines Re 172.9} plots the contours of the excess pressure over the mean gradient superposed with streamlines at various time instants after a stationary solution is obtained. We see multiple vortices oscillating with time.

\section{Summary}
A meshless semi-implicit iterative algorithm for the solution of incompressible Navier-Stokes equations is presented. The PHS-RBFs with appended polynomials give exponential order of convergence for the degree of the appended polynomial. The BDF2 scheme is used for temporal integration of the momentum equations. A Poisson equation for pressure corrections is formulated by imposing divergence free condition on the iterated velocity field. At each time step, the momentum and pressure correction equations are repeatedly solved until the velocities and pressure converge to a pre-specified tolerance. We demonstrate that this semi-implicit procedure is stable and convergent for high Courant numbers.
\par In this paper, we have demonstrated the convergence and discretization accuracy of the algorithm for two model problems and simulated three other complex problems. In all cases, the algorithm is stable for Courant numbers in excess of ten. No additional under-relaxation is needed to converge the nonlinearities and coupling between the equations. The model problems include a circular Couette flow with exact analytical solution, and the decay of an inital flow field to a null solution. In both cases, we demonstrate the exponential convergence of the discretization error per the degree of the appended polynomial. Comparison of the CPU runtimes with the previously published fractional step method \cite{shahane2020high} shows that this algorithm is around 5 to 6 times faster for the Courant number of 12. Subsequently, both steady and unsteady flows over a circular cylinder are simulated and the lift and drag coefficients are compared with values in the literature. The method is also applied to simulate two practical problems. The first problem deals with flow in the gap formed between two rotating circular cylinders and a stationary outer ellipse. Results are presented for the cases of unidirectional and counter rotations. Next we have simulated flow inside a sinusoidal bellow with periodic boundary conditions. Two Reynolds numbers are considered showing steady and unsteady flows. The unsteady case shows complex flow patterns with multiple oscillating vortices. For all these problems, a systematic grid independence study followed by line and contour plots are presented. These are in qualitative agreement with published results in literature.
\par The algorithm has a potential to accurately and efficiently solve many fluid flow and heat transfer problems in complex domains. An open source code Meshless Multi-Physics Software (MeMPhyS) \cite{shahanememphys} is available for interested users of the algorithm. Future extensions will include a parallel implementation of the software on high-performance CPUs and GPUs and implementation of different physical models for turbulence and multiple phases.

\bibliography{References}

\begin{thebibliography}{70}
\expandafter\ifx\csname natexlab\endcsname\relax\def\natexlab#1{#1}\fi
\providecommand{\url}[1]{\texttt{#1}}
\providecommand{\href}[2]{#2}
\providecommand{\path}[1]{#1}
\providecommand{\DOIprefix}{doi:}
\providecommand{\ArXivprefix}{arXiv:}
\providecommand{\URLprefix}{URL: }
\providecommand{\Pubmedprefix}{pmid:}
\providecommand{\doi}[1]{\href{http://dx.doi.org/#1}{\path{#1}}}
\providecommand{\Pubmed}[1]{\href{pmid:#1}{\path{#1}}}
\providecommand{\bibinfo}[2]{#2}
\ifx\xfnm\relax \def\xfnm[#1]{\unskip,\space#1}\fi
\bibitem[{Shahane and Vanka(2021)}]{shahanememphys}
\bibinfo{author}{S.~Shahane}, \bibinfo{author}{S.~P. Vanka},
  \bibinfo{title}{{MeMPhyS: Meshless Multi-Physics Software}},
  \bibinfo{year}{2021}. \URLprefix
  \url{https://github.com/shahaneshantanu/memphys}.
\bibitem[{Harlow and Welch(1965)}]{harlow1965numerical}
\bibinfo{author}{F.~Harlow}, \bibinfo{author}{J.~Welch},
\newblock \bibinfo{title}{Numerical calculation of time-dependent viscous
  incompressible flow of fluid with free surface},
\newblock \bibinfo{journal}{The Physics of Fluids} \bibinfo{volume}{8}
  (\bibinfo{year}{1965}) \bibinfo{pages}{2182--2189}.
\bibitem[{Chorin(1968)}]{chorin1968numerical}
\bibinfo{author}{A.~J. Chorin},
\newblock \bibinfo{title}{Numerical solution of the navier-stokes equations},
\newblock \bibinfo{journal}{Mathematics of computation} \bibinfo{volume}{22}
  (\bibinfo{year}{1968}) \bibinfo{pages}{745--762}.
\bibitem[{Temam(1969)}]{temam1969approximation}
\bibinfo{author}{R.~Temam},
\newblock \bibinfo{title}{Approximation of the solution of the navier- stokes
  equations by the fractional step method(approximate solution method for
  navier-stokes equations for incompressible viscous fluids)},
\newblock \bibinfo{journal}{Archive for Rational Mechanics and Analysis}
  \bibinfo{volume}{32} (\bibinfo{year}{1969}) \bibinfo{pages}{135--153}.
\bibitem[{Kim and Moin(1985)}]{kim1985application}
\bibinfo{author}{J.~Kim}, \bibinfo{author}{P.~Moin},
\newblock \bibinfo{title}{Application of a fractional-step method to
  incompressible navier-stokes equations},
\newblock \bibinfo{journal}{Journal of computational physics}
  \bibinfo{volume}{59} (\bibinfo{year}{1985}) \bibinfo{pages}{308--323}.
\bibitem[{Rosenfeld et~al.(1991)Rosenfeld, Kwak, and
  Vinokur}]{rosenfeld1991fractional}
\bibinfo{author}{M.~Rosenfeld}, \bibinfo{author}{D.~Kwak},
  \bibinfo{author}{M.~Vinokur},
\newblock \bibinfo{title}{A fractional step solution method for the unsteady
  incompressible navier-stokes equations in generalized coordinate systems},
\newblock \bibinfo{journal}{Journal of Computational Physics}
  \bibinfo{volume}{94} (\bibinfo{year}{1991}) \bibinfo{pages}{102--137}.
\bibitem[{Shahane et~al.(2019)Shahane, Aluru, Ferreira, Kapoor, and
  Vanka}]{shahane2019finite}
\bibinfo{author}{S.~Shahane}, \bibinfo{author}{N.~Aluru},
  \bibinfo{author}{P.~Ferreira}, \bibinfo{author}{S.~G. Kapoor},
  \bibinfo{author}{S.~P. Vanka},
\newblock \bibinfo{title}{Finite volume simulation framework for die casting
  with uncertainty quantification},
\newblock \bibinfo{journal}{Applied Mathematical Modelling}
  \bibinfo{volume}{74} (\bibinfo{year}{2019}) \bibinfo{pages}{132--150}.
\bibitem[{Issa(1986)}]{issa1986solution}
\bibinfo{author}{R.~I. Issa},
\newblock \bibinfo{title}{Solution of the implicitly discretised fluid flow
  equations by operator-splitting},
\newblock \bibinfo{journal}{Journal of computational physics}
  \bibinfo{volume}{62} (\bibinfo{year}{1986}) \bibinfo{pages}{40--65}.
\bibitem[{Van~Doormaal and Raithby(1984)}]{van1984enhancements}
\bibinfo{author}{J.~P. Van~Doormaal}, \bibinfo{author}{G.~D. Raithby},
\newblock \bibinfo{title}{Enhancements of the simple method for predicting
  incompressible fluid flows},
\newblock \bibinfo{journal}{Numerical heat transfer} \bibinfo{volume}{7}
  (\bibinfo{year}{1984}) \bibinfo{pages}{147--163}.
\bibitem[{Patankar and Spalding(1983)}]{patankar1983calculation}
\bibinfo{author}{S.~V. Patankar}, \bibinfo{author}{D.~B. Spalding},
\newblock \bibinfo{title}{A calculation procedure for heat, mass and momentum
  transfer in three-dimensional parabolic flows},
\newblock in: \bibinfo{booktitle}{Numerical prediction of flow, heat transfer,
  turbulence and combustion}, \bibinfo{publisher}{Elsevier},
  \bibinfo{year}{1983}, pp. \bibinfo{pages}{54--73}.
\bibitem[{Patankar et~al.(1983)Patankar, Pratap, and
  Spalding}]{patankar1983prediction}
\bibinfo{author}{S.~V. Patankar}, \bibinfo{author}{V.~S. Pratap},
  \bibinfo{author}{D.~B. Spalding},
\newblock \bibinfo{title}{Prediction of laminar flow and heat transfer in
  helically coiled pipes},
\newblock in: \bibinfo{booktitle}{Numerical Prediction of Flow, Heat Transfer,
  Turbulence and Combustion}, \bibinfo{publisher}{Elsevier},
  \bibinfo{year}{1983}, pp. \bibinfo{pages}{117--129}.
\bibitem[{Patankar(2018)}]{patankar2018numerical}
\bibinfo{author}{S.~V. Patankar}, \bibinfo{title}{Numerical heat transfer and
  fluid flow}, \bibinfo{publisher}{Taylor \& Francis}, \bibinfo{year}{2018}.
\bibitem[{Vanka(1985)}]{vanka1985block}
\bibinfo{author}{S.~P. Vanka},
\newblock \bibinfo{title}{Block-implicit calculation of steady turbulent
  recirculating flows},
\newblock \bibinfo{journal}{International journal of heat and mass transfer}
  \bibinfo{volume}{28} (\bibinfo{year}{1985}) \bibinfo{pages}{2093--2103}.
\bibitem[{Braaten and Patankar(1990)}]{braaten1990fully}
\bibinfo{author}{M.~Braaten}, \bibinfo{author}{S.~V. Patankar},
\newblock \bibinfo{title}{Fully coupled solution of the equations for
  incompressible recirculating flows using a penalty-function finite-difference
  formulation},
\newblock \bibinfo{journal}{Computational mechanics} \bibinfo{volume}{6}
  (\bibinfo{year}{1990}) \bibinfo{pages}{143--155}.
\bibitem[{Darwish et~al.(2009)Darwish, Sraj, and
  Moukalled}]{darwish2009coupled}
\bibinfo{author}{M.~Darwish}, \bibinfo{author}{I.~Sraj},
  \bibinfo{author}{F.~Moukalled},
\newblock \bibinfo{title}{A coupled finite volume solver for the solution of
  incompressible flows on unstructured grids},
\newblock \bibinfo{journal}{Journal of Computational Physics}
  \bibinfo{volume}{228} (\bibinfo{year}{2009}) \bibinfo{pages}{180--201}.
\bibitem[{Vanka and Leaf(1983)}]{vanka1983fully}
\bibinfo{author}{S.~P. Vanka}, \bibinfo{author}{G.~Leaf},
  \bibinfo{title}{Fully-coupled solution of pressure-linked fluid flow
  equations}, \bibinfo{type}{Technical Report}, Argonne National Lab., IL
  (USA), \bibinfo{year}{1983}.
\bibitem[{Chen and Przekwas(2010)}]{chen2010coupled}
\bibinfo{author}{Z.~Chen}, \bibinfo{author}{A.~Przekwas},
\newblock \bibinfo{title}{A coupled pressure-based computational method for
  incompressible/compressible flows},
\newblock \bibinfo{journal}{Journal of Computational Physics}
  \bibinfo{volume}{229} (\bibinfo{year}{2010}) \bibinfo{pages}{9150--9165}.
\bibitem[{Sohankar(2004)}]{sohankar2004and}
\bibinfo{author}{A.~Sohankar},
\newblock \bibinfo{title}{The les and dns simulations of heat transfer and
  fluid flow in a plate-fin heat exchanger with vortex generators}
  (\bibinfo{year}{2004}).
\bibitem[{Manhart(2004)}]{manhart2004zonal}
\bibinfo{author}{M.~Manhart},
\newblock \bibinfo{title}{A zonal grid algorithm for dns of turbulent boundary
  layers},
\newblock \bibinfo{journal}{Computers \& Fluids} \bibinfo{volume}{33}
  (\bibinfo{year}{2004}) \bibinfo{pages}{435--461}.
\bibitem[{Ha et~al.(2018)Ha, Park, and You}]{ha2018gpu}
\bibinfo{author}{S.~Ha}, \bibinfo{author}{J.~Park}, \bibinfo{author}{D.~You},
\newblock \bibinfo{title}{A gpu-accelerated semi-implicit fractional-step
  method for numerical solutions of incompressible navier--stokes equations},
\newblock \bibinfo{journal}{Journal of Computational Physics}
  \bibinfo{volume}{352} (\bibinfo{year}{2018}) \bibinfo{pages}{246--264}.
\bibitem[{Kalitzin et~al.(2003)Kalitzin, Wu, and Durbin}]{kalitzin2003dns}
\bibinfo{author}{G.~Kalitzin}, \bibinfo{author}{X.~Wu}, \bibinfo{author}{P.~A.
  Durbin},
\newblock \bibinfo{title}{Dns of fully turbulent flow in a lpt passage},
\newblock \bibinfo{journal}{International Journal of Heat and Fluid Flow}
  \bibinfo{volume}{24} (\bibinfo{year}{2003}) \bibinfo{pages}{636--644}.
\bibitem[{You et~al.(2000)You, Choi, and Yoo}]{you2000modified}
\bibinfo{author}{J.~You}, \bibinfo{author}{H.~Choi}, \bibinfo{author}{J.~Y.
  Yoo},
\newblock \bibinfo{title}{A modified fractional step method of keeping a
  constant mass flow rate in fully developed channel and pipe flows},
\newblock \bibinfo{journal}{KSME international journal} \bibinfo{volume}{14}
  (\bibinfo{year}{2000}) \bibinfo{pages}{547--552}.
\bibitem[{Robichaux et~al.(1999)Robichaux, Balachandar, and
  Vanka}]{robichaux1999three}
\bibinfo{author}{J.~Robichaux}, \bibinfo{author}{S.~Balachandar},
  \bibinfo{author}{S.~P. Vanka},
\newblock \bibinfo{title}{Three-dimensional floquet instability of the wake of
  square cylinder},
\newblock \bibinfo{journal}{Physics of Fluids} \bibinfo{volume}{11}
  (\bibinfo{year}{1999}) \bibinfo{pages}{560--578}.
\bibitem[{nek(2021)}]{nek5000}
\bibinfo{title}{{NEK: a fast and scalable high-order solver for computational
  fluid dynamics}}, \bibinfo{year}{2021}. \URLprefix
  \url{https://github.com/Nek5000}.
\bibitem[{Spalart et~al.(1991)Spalart, Moser, and Rogers}]{spalart1991spectral}
\bibinfo{author}{P.~R. Spalart}, \bibinfo{author}{R.~D. Moser},
  \bibinfo{author}{M.~M. Rogers},
\newblock \bibinfo{title}{Spectral methods for the navier-stokes equations with
  one infinite and two periodic directions},
\newblock \bibinfo{journal}{Journal of Computational Physics}
  \bibinfo{volume}{96} (\bibinfo{year}{1991}) \bibinfo{pages}{297--324}.
\bibitem[{Tafti and Vanka(1991{\natexlab{a}})}]{tafti1991numerical}
\bibinfo{author}{D.~Tafti}, \bibinfo{author}{S.~P. Vanka},
\newblock \bibinfo{title}{A numerical study of the effects of spanwise rotation
  on turbulent channel flow},
\newblock \bibinfo{journal}{Physics of Fluids A: Fluid Dynamics}
  \bibinfo{volume}{3} (\bibinfo{year}{1991}{\natexlab{a}})
  \bibinfo{pages}{642--656}.
\bibitem[{Tafti and Vanka(1991{\natexlab{b}})}]{tafti1991three}
\bibinfo{author}{D.~Tafti}, \bibinfo{author}{S.~P. Vanka},
\newblock \bibinfo{title}{A three-dimensional numerical study of flow
  separation and reattachment on a blunt plate},
\newblock \bibinfo{journal}{Physics of Fluids A: Fluid Dynamics}
  \bibinfo{volume}{3} (\bibinfo{year}{1991}{\natexlab{b}})
  \bibinfo{pages}{2887--2909}.
\bibitem[{Pawlowski et~al.(2006)Pawlowski, Shadid, Simonis, and
  Walker}]{pawlowski2006globalization}
\bibinfo{author}{R.~P. Pawlowski}, \bibinfo{author}{J.~N. Shadid},
  \bibinfo{author}{J.~P. Simonis}, \bibinfo{author}{H.~F. Walker},
\newblock \bibinfo{title}{Globalization techniques for newton--krylov methods
  and applications to the fully coupled solution of the navier--stokes
  equations},
\newblock \bibinfo{journal}{SIAM review} \bibinfo{volume}{48}
  (\bibinfo{year}{2006}) \bibinfo{pages}{700--721}.
\bibitem[{Knoll and Keyes(2004)}]{knoll2004jacobian}
\bibinfo{author}{D.~A. Knoll}, \bibinfo{author}{D.~E. Keyes},
\newblock \bibinfo{title}{Jacobian-free newton--krylov methods: a survey of
  approaches and applications},
\newblock \bibinfo{journal}{Journal of Computational Physics}
  \bibinfo{volume}{193} (\bibinfo{year}{2004}) \bibinfo{pages}{357--397}.
\bibitem[{Vanka(1986{\natexlab{a}})}]{vanka1986blockcmame}
\bibinfo{author}{S.~P. Vanka},
\newblock \bibinfo{title}{Block-implicit multigrid calculation of
  two-dimensional recirculating flows},
\newblock \bibinfo{journal}{Computer Methods in Applied Mechanics and
  Engineering} \bibinfo{volume}{59} (\bibinfo{year}{1986}{\natexlab{a}})
  \bibinfo{pages}{29--48}.
\bibitem[{Vanka(1986{\natexlab{b}})}]{vanka1986blockjcp}
\bibinfo{author}{S.~P. Vanka},
\newblock \bibinfo{title}{Block-implicit multigrid solution of navier-stokes
  equations in primitive variables},
\newblock \bibinfo{journal}{Journal of Computational Physics}
  \bibinfo{volume}{65} (\bibinfo{year}{1986}{\natexlab{b}})
  \bibinfo{pages}{138--158}.
\bibitem[{John and Tobiska(2000)}]{john2000numerical}
\bibinfo{author}{V.~John}, \bibinfo{author}{L.~Tobiska},
\newblock \bibinfo{title}{Numerical performance of smoothers in coupled
  multigrid methods for the parallel solution of the incompressible
  navier--stokes equations},
\newblock \bibinfo{journal}{International Journal for Numerical Methods in
  Fluids} \bibinfo{volume}{33} (\bibinfo{year}{2000})
  \bibinfo{pages}{453--473}.
\bibitem[{Monaghan(2012)}]{monaghan2012smoothed}
\bibinfo{author}{J.~Monaghan},
\newblock \bibinfo{title}{Smoothed particle hydrodynamics and its diverse
  applications},
\newblock \bibinfo{journal}{Annual Review of Fluid Mechanics}
  \bibinfo{volume}{44} (\bibinfo{year}{2012}) \bibinfo{pages}{323--346}.
\bibitem[{Perrone and Kao(1975)}]{perrone1975general}
\bibinfo{author}{N.~Perrone}, \bibinfo{author}{R.~Kao},
\newblock \bibinfo{title}{A general finite difference method for arbitrary
  meshes},
\newblock \bibinfo{journal}{Computers \& Structures} \bibinfo{volume}{5}
  (\bibinfo{year}{1975}) \bibinfo{pages}{45--57}.
\bibitem[{Liu et~al.(1995)Liu, Jun, and Zhang}]{liu1995reproducing}
\bibinfo{author}{W.~Liu}, \bibinfo{author}{S.~Jun}, \bibinfo{author}{Y.~Zhang},
\newblock \bibinfo{title}{Reproducing kernel particle methods},
\newblock \bibinfo{journal}{International journal for numerical methods in
  fluids} \bibinfo{volume}{20} (\bibinfo{year}{1995})
  \bibinfo{pages}{1081--1106}.
\bibitem[{Belytschko et~al.(1994)Belytschko, Lu, and
  Gu}]{belytschko1994element}
\bibinfo{author}{T.~Belytschko}, \bibinfo{author}{Y.~Lu},
  \bibinfo{author}{L.~Gu},
\newblock \bibinfo{title}{Element-free galerkin methods},
\newblock \bibinfo{journal}{International journal for numerical methods in
  engineering} \bibinfo{volume}{37} (\bibinfo{year}{1994})
  \bibinfo{pages}{229--256}.
\bibitem[{Liszka et~al.(1996)Liszka, Duarte, and Tworzydlo}]{liszka1996hp}
\bibinfo{author}{T.~Liszka}, \bibinfo{author}{C.~Duarte},
  \bibinfo{author}{W.~Tworzydlo},
\newblock \bibinfo{title}{hp-meshless cloud method},
\newblock \bibinfo{journal}{Computer Methods in Applied Mechanics and
  Engineering} \bibinfo{volume}{139} (\bibinfo{year}{1996})
  \bibinfo{pages}{263--288}.
\bibitem[{Babu{\v{s}}ka and Melenk(1997)}]{babuvska1997partition}
\bibinfo{author}{I.~Babu{\v{s}}ka}, \bibinfo{author}{J.~Melenk},
\newblock \bibinfo{title}{The partition of unity method},
\newblock \bibinfo{journal}{International journal for numerical methods in
  engineering} \bibinfo{volume}{40} (\bibinfo{year}{1997})
  \bibinfo{pages}{727--758}.
\bibitem[{O{\~n}ate et~al.(1996)O{\~n}ate, Idelsohn, Zienkiewicz, and
  Taylor}]{onate1996finite}
\bibinfo{author}{E.~O{\~n}ate}, \bibinfo{author}{S.~Idelsohn},
  \bibinfo{author}{O.~Zienkiewicz}, \bibinfo{author}{R.~Taylor},
\newblock \bibinfo{title}{A finite point method in computational mechanics.
  applications to convective transport and fluid flow},
\newblock \bibinfo{journal}{International journal for numerical methods in
  engineering} \bibinfo{volume}{39} (\bibinfo{year}{1996})
  \bibinfo{pages}{3839--3866}.
\bibitem[{Shu et~al.(2003)Shu, Ding, and Yeo}]{shu2003local}
\bibinfo{author}{C.~Shu}, \bibinfo{author}{H.~Ding}, \bibinfo{author}{K.~Yeo},
\newblock \bibinfo{title}{Local radial basis function--based differential
  quadrature method and its application to solve two--dimensional
  incompressible navier--stokes equations},
\newblock \bibinfo{journal}{Computer methods in applied mechanics and
  engineering} \bibinfo{volume}{192} (\bibinfo{year}{2003})
  \bibinfo{pages}{941--954}.
\bibitem[{Ding et~al.(2006)Ding, Shu, Yeo, and Xu}]{ding2006numerical}
\bibinfo{author}{H.~Ding}, \bibinfo{author}{C.~Shu}, \bibinfo{author}{K.~Yeo},
  \bibinfo{author}{D.~Xu},
\newblock \bibinfo{title}{Numerical computation of three--dimensional
  incompressible viscous flows in the primitive variable form by local
  multiquadric differential quadrature method},
\newblock \bibinfo{journal}{Computer Methods in Applied Mechanics and
  Engineering} \bibinfo{volume}{195} (\bibinfo{year}{2006})
  \bibinfo{pages}{516--533}.
\bibitem[{Sanyasiraju and Chandhini(2008)}]{sanyasiraju2008local}
\bibinfo{author}{Y.~Sanyasiraju}, \bibinfo{author}{G.~Chandhini},
\newblock \bibinfo{title}{Local radial basis function based gridfree scheme for
  unsteady incompressible viscous flows},
\newblock \bibinfo{journal}{Journal of Computational Physics}
  \bibinfo{volume}{227} (\bibinfo{year}{2008}) \bibinfo{pages}{8922--8948}.
\bibitem[{Vidal et~al.(2016)Vidal, Kassab, and Divo}]{vidal2016direct}
\bibinfo{author}{A.~Vidal}, \bibinfo{author}{A.~Kassab},
  \bibinfo{author}{E.~Divo},
\newblock \bibinfo{title}{{A direct velocity--pressure coupling Meshless
  algorithm for incompressible fluid flow simulations}},
\newblock \bibinfo{journal}{Engineering Analysis with Boundary Elements}
  \bibinfo{volume}{72} (\bibinfo{year}{2016}) \bibinfo{pages}{1--10}.
\bibitem[{Zamolo and Nobile(2019)}]{zamolo2019solution}
\bibinfo{author}{R.~Zamolo}, \bibinfo{author}{E.~Nobile},
\newblock \bibinfo{title}{{Solution of incompressible fluid flow problems with
  heat transfer by means of an efficient RBF--FD meshless approach}},
\newblock \bibinfo{journal}{Numerical Heat Transfer, Part B: Fundamentals}
  (\bibinfo{year}{2019}) \bibinfo{pages}{1--24}.
\bibitem[{Shahane et~al.(2020)Shahane, Radhakrishnan, and
  Vanka}]{shahane2020high}
\bibinfo{author}{S.~Shahane}, \bibinfo{author}{A.~Radhakrishnan},
  \bibinfo{author}{S.~P. Vanka},
\newblock \bibinfo{title}{A high-order accurate meshless method for solution of
  incompressible fluid flow problems},
\newblock \bibinfo{journal}{arXiv preprint arXiv:2010.01702}
  (\bibinfo{year}{2020}).
\bibitem[{Bayona et~al.(2017)Bayona, Flyer, Fornberg, and
  Barnett}]{bayona2017role_II}
\bibinfo{author}{V.~Bayona}, \bibinfo{author}{N.~Flyer},
  \bibinfo{author}{B.~Fornberg}, \bibinfo{author}{G.~A. Barnett},
\newblock \bibinfo{title}{On the role of polynomials in rbf-fd approximations:
  Ii. numerical solution of elliptic pdes},
\newblock \bibinfo{journal}{Journal of Computational Physics}
  \bibinfo{volume}{332} (\bibinfo{year}{2017}) \bibinfo{pages}{257--273}.
\bibitem[{Kosec and {\v{S}}arler(2008)}]{kosec2008solution}
\bibinfo{author}{G.~Kosec}, \bibinfo{author}{B.~{\v{S}}arler},
\newblock \bibinfo{title}{Solution of thermo-fluid problems by collocation with
  local pressure correction},
\newblock \bibinfo{journal}{International Journal of Numerical Methods for Heat
  \& Fluid Flow}  (\bibinfo{year}{2008}).
\bibitem[{Kosec and Slak(2020)}]{kosec2020radial}
\bibinfo{author}{G.~Kosec}, \bibinfo{author}{J.~Slak},
\newblock \bibinfo{title}{Radial basis function-generated finite differences
  solution of natural convection problem in 3d},
\newblock in: \bibinfo{booktitle}{AIP Conference Proceedings}, volume
  \bibinfo{volume}{2293}, \bibinfo{organization}{AIP Publishing LLC},
  \bibinfo{year}{2020}, p. \bibinfo{pages}{420094}.
\bibitem[{Flyer et~al.(2016)Flyer, Fornberg, Bayona, and
  Barnett}]{flyer2016onrole_I}
\bibinfo{author}{N.~Flyer}, \bibinfo{author}{B.~Fornberg},
  \bibinfo{author}{V.~Bayona}, \bibinfo{author}{G.~Barnett},
\newblock \bibinfo{title}{{On the role of polynomials in RBF--FD
  approximations: I. Interpolation and accuracy}},
\newblock \bibinfo{journal}{Journal of Computational Physics}
  \bibinfo{volume}{321} (\bibinfo{year}{2016}) \bibinfo{pages}{21--38}.
\bibitem[{Radhakrishnan et~al.(2021)Radhakrishnan, Xu, Shahane, and
  Vanka}]{radhakrishnan2021non}
\bibinfo{author}{A.~Radhakrishnan}, \bibinfo{author}{M.~Xu},
  \bibinfo{author}{S.~Shahane}, \bibinfo{author}{S.~P. Vanka},
\newblock \bibinfo{title}{A non-nested multilevel method for meshless solution
  of the poisson equation in heat transfer and fluid flow},
\newblock \bibinfo{journal}{arXiv preprint arXiv:2104.13758}
  (\bibinfo{year}{2021}).
\bibitem[{Hardy(1971)}]{hardy1971multiquadric}
\bibinfo{author}{R.~Hardy},
\newblock \bibinfo{title}{Multiquadric equations of topography and other
  irregular surfaces},
\newblock \bibinfo{journal}{Journal of geophysical research}
  \bibinfo{volume}{76} (\bibinfo{year}{1971}) \bibinfo{pages}{1905--1915}.
\bibitem[{Kansa(1990{\natexlab{a}})}]{kansa1990multiquadrics_I}
\bibinfo{author}{E.~Kansa},
\newblock \bibinfo{title}{{Multiquadrics—A scattered data approximation
  scheme with applications to computational fluid-dynamics—I surface
  approximations and partial derivative estimates}},
\newblock \bibinfo{journal}{Computers \& Mathematics with applications}
  \bibinfo{volume}{19} (\bibinfo{year}{1990}{\natexlab{a}})
  \bibinfo{pages}{127--145}.
\bibitem[{Kansa(1990{\natexlab{b}})}]{kansa1990multiquadrics_II}
\bibinfo{author}{E.~Kansa},
\newblock \bibinfo{title}{{Multiquadrics—A scattered data approximation
  scheme with applications to computational fluid-dynamics—II solutions to
  parabolic, hyperbolic and elliptic partial differential equations}},
\newblock \bibinfo{journal}{Computers \& mathematics with applications}
  \bibinfo{volume}{19} (\bibinfo{year}{1990}{\natexlab{b}})
  \bibinfo{pages}{147--161}.
\bibitem[{Kansa and Hon(2000)}]{kansa2000circumventing}
\bibinfo{author}{E.~Kansa}, \bibinfo{author}{Y.~Hon},
\newblock \bibinfo{title}{Circumventing the ill-conditioning problem with
  multiquadric radial basis functions: applications to elliptic partial
  differential equations},
\newblock \bibinfo{journal}{Computers \& Mathematics with applications}
  \bibinfo{volume}{39} (\bibinfo{year}{2000}) \bibinfo{pages}{123--137}.
\bibitem[{S{\"u}li and Mayers(2003)}]{suli2003introduction}
\bibinfo{author}{E.~S{\"u}li}, \bibinfo{author}{D.~F. Mayers},
  \bibinfo{title}{An introduction to numerical analysis},
  \bibinfo{publisher}{Cambridge university press}, \bibinfo{year}{2003}.
\bibitem[{LeVeque(2007)}]{leveque2007finite}
\bibinfo{author}{R.~J. LeVeque}, \bibinfo{title}{Finite difference methods for
  ordinary and partial differential equations: steady-state and time-dependent
  problems}, \bibinfo{publisher}{SIAM}, \bibinfo{year}{2007}.
\bibitem[{Bell et~al.(1989)Bell, Colella, and Glaz}]{bell1989second}
\bibinfo{author}{J.~B. Bell}, \bibinfo{author}{P.~Colella},
  \bibinfo{author}{H.~M. Glaz},
\newblock \bibinfo{title}{A second-order projection method for the
  incompressible navier-stokes equations},
\newblock \bibinfo{journal}{Journal of Computational Physics}
  \bibinfo{volume}{85} (\bibinfo{year}{1989}) \bibinfo{pages}{257--283}.
\bibitem[{Takami and Keller(1969)}]{takami1969steady}
\bibinfo{author}{H.~Takami}, \bibinfo{author}{H.~B. Keller},
\newblock \bibinfo{title}{Steady two-dimensional viscous flow of an
  incompressible fluid past a circular cylinder},
\newblock \bibinfo{journal}{The Physics of Fluids} \bibinfo{volume}{12}
  (\bibinfo{year}{1969}) \bibinfo{pages}{II--51}.
\bibitem[{Tuann and Olson(1978)}]{tuann1978numerical}
\bibinfo{author}{S.-Y. Tuann}, \bibinfo{author}{M.~D. Olson},
\newblock \bibinfo{title}{Numerical studies of the flow around a circular
  cylinder by a finite element method},
\newblock \bibinfo{journal}{Computers \& Fluids} \bibinfo{volume}{6}
  (\bibinfo{year}{1978}) \bibinfo{pages}{219--240}.
\bibitem[{Ding et~al.(2004)Ding, Shu, Yeo, and Xu}]{ding2004simulation}
\bibinfo{author}{H.~Ding}, \bibinfo{author}{C.~Shu}, \bibinfo{author}{K.~Yeo},
  \bibinfo{author}{D.~Xu},
\newblock \bibinfo{title}{Simulation of incompressible viscous flows past a
  circular cylinder by hybrid fd scheme and meshless least square-based finite
  difference method},
\newblock \bibinfo{journal}{Computer Methods in Applied Mechanics and
  Engineering} \bibinfo{volume}{193} (\bibinfo{year}{2004})
  \bibinfo{pages}{727--744}.
\bibitem[{Fornberg(1980)}]{fornberg1980numerical}
\bibinfo{author}{B.~Fornberg},
\newblock \bibinfo{title}{A numerical study of steady viscous flow past a
  circular cylinder},
\newblock \bibinfo{journal}{Journal of Fluid Mechanics} \bibinfo{volume}{98}
  (\bibinfo{year}{1980}) \bibinfo{pages}{819--855}.
\bibitem[{Nieuwstadt and Keller(1973)}]{nieuwstadt1973viscous}
\bibinfo{author}{F.~Nieuwstadt}, \bibinfo{author}{H.~Keller},
\newblock \bibinfo{title}{Viscous flow past circular cylinders},
\newblock \bibinfo{journal}{Computers \& Fluids} \bibinfo{volume}{1}
  (\bibinfo{year}{1973}) \bibinfo{pages}{59--71}.
\bibitem[{Gushchin and Shchennikov(1974)}]{gushchin1974numerical}
\bibinfo{author}{V.~A. Gushchin}, \bibinfo{author}{V.~Shchennikov},
\newblock \bibinfo{title}{A numerical method of solving the navier-stokes
  equations},
\newblock \bibinfo{journal}{USSR Computational Mathematics and Mathematical
  Physics} \bibinfo{volume}{14} (\bibinfo{year}{1974})
  \bibinfo{pages}{242--250}.
\bibitem[{Dennis and Chang(1970)}]{dennis1970numerical}
\bibinfo{author}{S.~Dennis}, \bibinfo{author}{G.-Z. Chang},
\newblock \bibinfo{title}{Numerical solutions for steady flow past a circular
  cylinder at reynolds numbers up to 100},
\newblock \bibinfo{journal}{Journal of Fluid Mechanics} \bibinfo{volume}{42}
  (\bibinfo{year}{1970}) \bibinfo{pages}{471--489}.
\bibitem[{Geuzaine and Remacle(2009)}]{geuzaine2009gmsh}
\bibinfo{author}{C.~Geuzaine}, \bibinfo{author}{J.~Remacle},
\newblock \bibinfo{title}{Gmsh: A 3-d finite element mesh generator with
  built-in pre-and post-processing facilities},
\newblock \bibinfo{journal}{International Journal for Numerical Methods in
  Engineering} \bibinfo{volume}{79} (\bibinfo{year}{2009})
  \bibinfo{pages}{1309--1331}.
\bibitem[{Braza et~al.(1986)Braza, Chassaing, and Minh}]{braza1986numerical}
\bibinfo{author}{M.~Braza}, \bibinfo{author}{P.~Chassaing},
  \bibinfo{author}{H.~H. Minh},
\newblock \bibinfo{title}{Numerical study and physical analysis of the pressure
  and velocity fields in the near wake of a circular cylinder},
\newblock \bibinfo{journal}{Journal of fluid mechanics} \bibinfo{volume}{165}
  (\bibinfo{year}{1986}) \bibinfo{pages}{79--130}.
\bibitem[{Liu et~al.(1998)Liu, Zheng, and Sung}]{liu1998preconditioned}
\bibinfo{author}{C.~Liu}, \bibinfo{author}{X.~Zheng},
  \bibinfo{author}{C.~Sung},
\newblock \bibinfo{title}{Preconditioned multigrid methods for unsteady
  incompressible flows},
\newblock \bibinfo{journal}{Journal of Computational physics}
  \bibinfo{volume}{139} (\bibinfo{year}{1998}) \bibinfo{pages}{35--57}.
\bibitem[{Belov et~al.(1995)Belov, Martinelli, and Jameson}]{belov1995new}
\bibinfo{author}{A.~Belov}, \bibinfo{author}{L.~Martinelli},
  \bibinfo{author}{A.~Jameson},
\newblock \bibinfo{title}{A new implicit algorithm with multigrid for unsteady
  incompressible flow calculations},
\newblock in: \bibinfo{booktitle}{33rd Aerospace sciences meeting and exhibit},
  \bibinfo{year}{1995}, p.~\bibinfo{pages}{49}.
\bibitem[{Wang and Vanka(1995)}]{wang1995convective}
\bibinfo{author}{G.~Wang}, \bibinfo{author}{S.~P. Vanka},
\newblock \bibinfo{title}{Convective heat transfer in periodic wavy passages},
\newblock \bibinfo{journal}{International Journal of Heat and Mass Transfer}
  \bibinfo{volume}{38} (\bibinfo{year}{1995}) \bibinfo{pages}{3219--3230}.
\bibitem[{Stephanoff et~al.(1980)Stephanoff, Sobey, and
  Bellhouse}]{stephanoff1980flow}
\bibinfo{author}{K.~Stephanoff}, \bibinfo{author}{I.~J. Sobey},
  \bibinfo{author}{B.~Bellhouse},
\newblock \bibinfo{title}{On flow through furrowed channels. part 2. observed
  flow patterns},
\newblock \bibinfo{journal}{Journal of Fluid Mechanics} \bibinfo{volume}{96}
  (\bibinfo{year}{1980}) \bibinfo{pages}{27--32}.

\end{thebibliography}

\end{document}